\pdfoutput=1
\documentclass{article}

\usepackage{ifthen}

\newboolean{arxiv}
\newboolean{sisc}
\newboolean{toms}

\newboolean{markNewContent}

\setboolean{arxiv}{true}
\setboolean{sisc}{false}
\setboolean{toms}{false}

\setboolean{markNewContent}{false}

\ifthenelse{\boolean{arxiv}}{
  \usepackage{algorithm}  
\newtheorem{definition}{Definition} 

\usepackage{geometry}
\usepackage{amsfonts}

\newcommand{\TheTitle}{Quasi-matrix-free hybrid multigrid on dynamically adaptive Cartesian grids} 

\title{\TheTitle}

\author{
  Marion Weinzierl
  \thanks{
    Department of Mathematical Sciences, Durham University,
    UK.
   (marion.weinzierl@durham.ac.uk)
  }
  \and
  Tobias Weinzierl
  \thanks{
    School of Engineering and Computing Sciences, Durham University,
    UK.
    (tobias.weinzierl@durham.ac.uk)
  }
}

\newtheorem{observation}{Observation}
\newcommand{\tbl}{\caption }

}{} 

\ifthenelse{\boolean{sisc}}{
  
\newcommand{\TheTitle}{Quasi-matrix-free hybrid multigrid on dynamically adaptive Cartesian grids}

\newcommand{\TheAuthors}{M. Weinzierl and T. Weinzierl}

\headers{Quasi-matrix-free hybrid multigrid}{\TheAuthors}

\title{\TheTitle}

\Crefname{ALC@unique}{Line}{Lines}

\author{
  Marion Weinzierl
  \thanks{
    Department of Mathematical Sciences, 
    Durham University
    (\email{marion.weinzierl@dur.ac.uk})
  }
  \and
  Tobias Weinzierl
  \thanks{
    School of Engineering and Computing Sciences, Durham University
    (\email{tobias.weinzierl@dur.ac.uk})
  }
}

\slugger{sisc}{xxxx}{xx}{x}{x--x}

\setlength\floatsep{0.4\baselineskip}
\setlength\textfloatsep{0.4\baselineskip}

}{} 

\ifthenelse{\boolean{toms}}{
\doi{0000001.0000001}

\issn{1234-56789}

\newcommand{\TheTitle}{Quasi-matrix-free hybrid multigrid on dynamically adaptive Cartesian grids}

\newcommand{\TheAuthors}{M. Weinzierl and T. Weinzierl}

\title{\TheTitle}

\markboth{\TheAuthors}{\TheTitle}

\author{
  MARION WEINZIERL
  \affil{
    Department of Mathematical Sciences, 
    Durham University
  }
  TOBIAS WEINZIERL
  \affil{
    School of Engineering and Computing Sciences, 
    Durham University
  }
}

\acmformat{
  Marion Weinzierl, Tobias Weinzierl, 2016. \TheTitle
 }

\category{}{Mathematics of computing}{Mathematical software}[Solvers]
\category{}{Computing methodologies}{ Symbolic and algebraic algorithms}[Linear algebra algorithms]
\category{}{Computing methodologies}{Modeling and simulation}[Multiscale systems]

\newtheorem{observation}[theorem]{Observation}

}{}

\usepackage{graphicx}
\usepackage[numbers,sort&compress]{natbib}
\usepackage{xcolor}
\usepackage{pdflscape}
\usepackage{algorithm}
\usepackage{algpseudocode}
\usepackage{amssymb}

\usepackage{enumitem}

\newcommand{\SIf}[1]{\State \algorithmicif\ {#1} \algorithmicthen}
\newcommand{\EndSIf}{ \unskip \algorithmicend\ \algorithmicif}

\ifthenelse{\boolean{markNewContent}}{
  \newcommand{\new}[1]{\textcolor{blue}{#1}}
}{
  \newcommand{\new}[1]{#1}
}

\begin{document}

%
%
\ifthenelse{\boolean{toms}}{
  \begin{abstract}
We \new{present} a family of spacetree-based multigrid realizations using the
tree's multiscale nature to derive coarse grids.
They align with matrix-free geometric multigrid solvers as they
never assemble the system matrices which is cumbersome for
dynamically adaptive grids and full multigrid.
The most sophisticated realizations use BoxMG \new{to construct}
operator-dependent prolongation and restriction in combination with Galerkin/Petrov-Galerkin coarse-grid operators.
This yields robust solvers for nontrivial elliptic problems.
\new{We embed} the algebraic, problem- and grid-dependent
multigrid operators as stencils into the grid and evaluate all
matrix-vector products in-situ throughout the grid traversals.
While such \new{ an approach} is not literally
matrix-free---the grid carries the matrix---we propose to switch to a hierarchical representation of all
operators.
Only differences of algebraic operators to their geometric
counterparts are held.
These hierarchical differences can be stored and exchanged with small memory
footprint.
Our realizations support arbitrary dynamically adaptive grids while they
vertically integrate the multilevel operations through spacetree linearization.
This yields good memory access characteristics, while standard colouring of mesh
entities with domain decomposition \new{allows} us to use parallel manycore
clusters.
%
\new{
All realization ingredients are detailed such that they can be
used by other codes.
}
\end{abstract}

  \maketitle
}{
  \maketitle
  \begin{abstract}
We \new{present} a family of spacetree-based multigrid realizations using the
tree's multiscale nature to derive coarse grids.
They align with matrix-free geometric multigrid solvers as they
never assemble the system matrices which is cumbersome for
dynamically adaptive grids and full multigrid.
The most sophisticated realizations use BoxMG \new{to construct}
operator-dependent prolongation and restriction in combination with Galerkin/Petrov-Galerkin coarse-grid operators.
This yields robust solvers for nontrivial elliptic problems.
\new{We embed} the algebraic, problem- and grid-dependent
multigrid operators as stencils into the grid and evaluate all
matrix-vector products in-situ throughout the grid traversals.
While such \new{ an approach} is not literally
matrix-free---the grid carries the matrix---we propose to switch to a hierarchical representation of all
operators.
Only differences of algebraic operators to their geometric
counterparts are held.
These hierarchical differences can be stored and exchanged with small memory
footprint.
Our realizations support arbitrary dynamically adaptive grids while they
vertically integrate the multilevel operations through spacetree linearization.
This yields good memory access characteristics, while standard colouring of mesh
entities with domain decomposition \new{allows} us to use parallel manycore
clusters.
%
\new{
All realization ingredients are detailed such that they can be
used by other codes.
}
\end{abstract}

}

\ifthenelse{\boolean{sisc}}{

}{}

\ifthenelse{\boolean{toms}}{
 \begin{bottomstuff}
%
   Author's addresses: 
   M. Weinzierl, Department of Mathematical Sciences, Durham University,
   DH1 3LE Durham, United Kingdom;
   T. Weinzierl, School of Engineering and Computing Sciences, Durham University,
   DH1 3LE Durham, United Kingdom
 \end{bottomstuff}
}{}

\section{Introduction}

%
%
Quadtrees, octrees and their generalization, all covered by the
term spacetree, are popular meshing paradigms in scientific computing, either
for plain meshing or as building blocks within forests of trees
\cite{Rudi:15:GordonBell,Burstedde:11:p4est,Mehl:06:MG,Reps:15:Helmholtz,Griebel:99:SFCAndMultigrid,Sundar:08:BalancedOctrees,Sundar:12:ParallelMultigrid,Weinzierl:11:Peano}.
They offer a geometric multiscale hierarchy, they support in-situ
meshing, they offer the opportunity to erase and refine parts of the grid
easily, their structuredness facilitates efficient implementations, and fast
domain decomposition methods are known for them---notably in combination with
space-filling curves (SFCs) \cite{Bader:13:SFCs}.
Therefore, they are well-suited for multigrid solvers for elliptic partial
differential equations (PDEs) discretized by finite elements.
Stagnating or decreasing memory per core, a widening memory gap and 
communication constraints hereby encourage the use of
matrix-free multigrid \cite{ascac10exa}.  Instead of assembling an equation system, all required matrix-vector multiplications
 (matvecs) are computed directly on the grid.
No matrix allocation or maintenance costs arise.

Matvecs without assembly or storage of matrix entries require that all operators can efficiently be determined on-the-fly
\new{from} the grid.
In algebraic multigrid (e.g. \cite{stueben01amg} and references therein),
inter-grid transfer (prolongation and restriction) and coarse-grid operators 
however depend recursively over several mesh levels on the PDE discretization.
\new{An on-the-fly computation requires the partial,
repeated and redundant re-assembly of the system matrix entries in every solver
sweep.
Furthermore, algebraic methods derive the coarse grid from the operator
characteristics, i.e.~a good coarse grid anticipates the problem's behaviour.
The grid and its data structures are not given a priori.
} 
Geometric multigrid (e.g. \cite{brandt84MGguide} and references therein)
in contrast \new{is better suited for matrix-free, fast implementations as 
it fixes the grid hierarchy and hardcodes or prescribes the operators. Unfortunately, it} quickly becomes unstable.

\new{We use the convection-diffusion equation
\begin{equation}
 - \nabla \new{\cdot} \left( \epsilon \nabla \right) u + \left( v \nabla \right) u = f
 \qquad \mbox{on\ } \Omega = (0,1)^d, u: \Omega \mapsto \mathbf{R} 
 \label{equation:pde}
\end{equation}
\noindent as demonstrator, where the diffusion coeffient $\epsilon $ is a
\new{rank-2 tensor with diagonal entries} $\epsilon _1, \ldots, \epsilon _d$ and the convection velocities $v=(v_1,\ldots,v_d)^T$ for
$d\in\{2,3\}$. It describes, for example, the transport of chemicals $u$ in a
fluid given by $v$ or yields a precursor for the Navier-Stokes equations
where $v\mapsto u$ and $f \mapsto f(u)$.
This equation is challenging for geometric multigrid if $\epsilon$, $v$ 
or---not studied in the present work---the shape of $\Omega$ are non-trivial and
if no particular/manual effort is invested to construct well-suited coarse
grids.
Algebraic multigrid, in contrast, can robustly and efficiently solve it as
black box:
Classic coarse grid identification places coarse grid points along
material-/$\epsilon $-transitions, makes coarse grids anticipate anisotropies
and convection, and it anticipates the domain layout.
Classic multigrid operator construction takes the same triad into account.
Both complement each other, i.e.~sophisticated AMG makes the operators
compensate weaker coarse grid choices and the other way round.}

%
%
\new{
We see three approaches for combining the robustness of algebraic multigrid
with the advantages of geometric multigrid.
First, we can use the spacetree concept to define the grid hierarchy and the fine-grid discretization (stencils), 
but compute the coarse-grid operators and inter-grid transfer operators algebraically (not on-the-fly).
An expensive coarse-grid/c-point identification phase, as necessary in algebraic
multigrid through matching on graphs, e.g., is thus omitted. At the same time,
the linear algebra subroutines and the memory allocation can take advantage of the structured grid.} As a second possibility, we may couple a geometric and an algebraic multigrid code.
On fine meshes, diffusion processes typically
dominate while material parameters are reasonably smooth. 
Matrix-free geometric multigrid based on rediscretization and fixed
inter-grid operators yield reasonable convergence.
On coarse meshes, one can employ algebraic multigrid, \new{alternative}
iterative schemes or a direct solver.
Such an approach
\cite{Gmeiner:15:HHG,Lu:14:HybridMG,Sundar:12:ParallelMultigrid,Rudi:15:GordonBell}
is robust and fast at modest memory requirements as long as the finest algebraic
problem with a matrix setup remains reasonably coarse within the grid hierarchy.
Finally, we may also decide to tailor all operators to the
problem manually and to hard-code it into the solver software.
Inter-grid transfer operators, for example, can anticipate the $\epsilon $ and
$v$ behaviour on finer meshes through homogenization.
\new{Such a strategy mirrors the fact that many engineers put significant
effort into the creation of high-quality simulation grids. We fix the grid but
invest into the operator.
Such an approach} does not work on-the-fly as black box.

%
%
\new{In this paper, we present a fourth technique that merges the advantages of
spacetrees plus rediscretization with the robustness of algebraic multigrid: space-tree based
multigrid using Petrov-Galerkin coarse-grid and operator-dependent 
inter-grid transfer operators.} 
\new{This} allows us to solve significantly harder elliptic problems without matrix
assembly than a classic geometric approach while we retain the advantages of
geometric multigrid.
\new{On the one hand,}
our solvers do not increase the memory footprint
significantly and thus are a convenient building block for extreme scale
simulations where the memory per core is a limiting factor.
On the other hand, our solvers preserve the geometric structuredness of
spacetrees which is an important characteristic for many optimisation techniques
and domain decomposition approaches.

Our description consists of three parts.
We start from an outline of the used tree data
structures, grid traversal and terminology
\new{(Sect.~\ref{sections:spacetrees})}. 
\new{In Sect.~\ref{sections:solver} we then revisit} three multigrid solver
variants that fit to spacetrees and spacetree traversals:
additive, BPX-type and multiplicative solvers \cite{Bastian:98:AdditiveVsMultiplicativeMG}
on hierarchical generating systems \cite{Griebel:94:Multilevelmethoden}.
All work in a strictly element-wise multiscale sense, that is, they only require
one cell or vertex record, respectively, plus their parents (next coarser
entities) at a time, \new{and they read each individual spacetree cell only
once per multigrid smoothing step.
This makes them well-prepared for future architecture with a widening memory
gap \cite{Dongarra:14:ApplMathExascaleComputing}.
We next introduce block Jacobi-type
smoothers that preserve the tree's single run-through policy, and we enhance
the geometric multigrid implementation with element-wise Galerkin coarse-grid operators.}
These operators are embedded into the tree, that is, the
grid acts both as an organizational and as a compute data structure. 
Galerkin coarse-grid operators improve the convergence rates, but the
solver becomes really robust only once we use operator-dependent prolongation
and restriction through Black Box Multigrid (BoxMG)  
\cite{dendy82blackbox, dendy83blackboxNonsym, dendy10blackboxCoarseBy3}.
BoxMG on locally refined grids for bipartitioning is known
\cite{shapira04matrixMG}.
\new{Our contribution is that we simplify its realisation on (dynamically)
adaptive grids through} the
fusion of BoxMG with full approximation storage (FAS)
\cite{trottenberg01multigrid} \new{realized through HTMG \cite{Griebel:90:HTMM}}.
Our code furthermore reduces the computation of BoxMG operators through
mirroring \new{all} computations \new{onto} one reference configuration.
\new{This renders the programming of BoxMG, notably for $d=3$, simpler than
published by any other author.}

Embedding all operators as stencils into the grid is not literally 
matrix-free.
In the second part of the manuscript (Sect.~\ref{sections:stencil-compression}),
we thus replace the storage of the
operators in the tree with a compressed encoding.
We determine the difference of all operators
to geometric rediscretization or $d$-linear operators, respectively, and store
only the hierarchical differences of algebraic to geometric stencils. 
In the best case, the difference is negligible and no floating point data is to
be held at all.
This reduces the memory footprint.
We work almost matrix-free but preserve BoxMG's robustness.
\new{To the best of our knowledge, such an approach to realize algebraic
operators almost at the memory cost of rediscretization-based
multigrid is new.}

The final part of the paper (Sect.~\ref{sections:parallelisation}) sketches
a shared and distributed memory parallelization.
The shared memory discussion derives data dependency graphs that can be
fed into a task-based system \cite{AlOnazi:17:TaskBasedAlgebraicMG}.
The distributed memory discussion studies data flow characteristics for
non-overlapping domain decompositions well-suited for MPI.
For both parallelization variants we show that our approach is by
construction well-suited for parallel machines; a property stemming from the
strict element-wise data access.

\new{Our contribution benefits from the fusion of three ingredients:
Operator-dependent prolongation and restriction on dynamically adaptive spacetrees, stencil compression and concurrency
analysis. These three parts} introduce two notions of hybrid
algebraic-geometric multigrid that are orthogonal to the classic notion of
hybrid where coarser grids are tackled by algebraic multigrid while fine grids
benefit from geometric multigrid with rediscretization: our approach is hybrid as we
(i) stick to the geometric multigrid structure but have algebraic operators and
(ii) determine algebraic operators but store only their difference to geometric
operators.
The present manuscript's \new{idioms perform} on reasonably well-posed problems. 
\new{They} enlarge the applicability of geometric multigrid with all its
geometric advantages.
This enlargement is made possible by an integration of known multigrid
techniques.
\new{Their elegant integration is,} to the best of our knowledge, new.

\section{Previous work and shortcomings of present approach}

\begin{table}
  \tbl{
  \new{
   Overview of discussed implementation techniques (first row) with references
   to corresponding manuscript sections (second row).
   Where appropriate, we add in the bottom part in which previous work we
   have studied related or similar concepts. Dissertations, as they have not
   appeared in peer-reviewed journals, are put in brackets.
  }
   \label{table:solver-variants}
  }{
  {\tiny
   \begin{tabular}{l|p{1.54cm}p{1.54cm}p{1.54cm}|p{1.54cm}|p{1.54cm}p{1.54cm}}
   Techniques: 
     & Matrix-free on spacetree with rediscretisation and HTMG
     & Element-wise tree block smoothers
     & Galerkin and BoxMG
     & Operator compression
     & Shared memory
     & Distributed memory
     \\
   \hline
   Section
     & Sect.~\ref{sections:solver:geometric}
     & Sect.~\ref{sections:blocksmoothers}
     & Sect.~\ref{sections:galerkin},\ref{sections:boxmg}
     & Sect.~\ref{sections:stencil-compression}
     & Sect.~\ref{section:parallelisation:sharedmem}
     & Sect.~\ref{section:parallelisation:distributedmem}
     \\ 
   \hline
   Additive  
     & \cite{Mehl:06:MG} without HTMG, \cite{Reps:15:Helmholtz}
     &&
     & 
     \\
   BPX 
     & \cite{Reps:15:Helmholtz}
     &&&      
     \\ 
   Multiplicative 
     & (\cite{Weinzierl:09:Diss})
     & (\cite{mweinzierl13diss}))
     & (\cite{mweinzierl13diss}))
     & \cite{Bungartz:10:Precompiler} for unknowns instead of stencils 
     & (\cite{mweinzierl13diss})
     & (\cite{Weinzierl:09:Diss}),(\cite{mweinzierl13diss})
     \\ 
     & 
     & 
     & 
     & \cite{Eckhardt:16:SPH} for unknowns in SPH 
     & 
     & 
   \end{tabular}
 }
}
\end{table}

\new{
Peano \cite{Software:Peano,Weinzierl:17:PeanoSoftware} serves as code base to
realize our single-touch tree traversals.
All implementation ideas however apply to other spacetree software, too.
Single-touch additive multigrid solvers for spacetrees with rediscretization are
subject of discussion in \cite{Mehl:06:MG}, though the discussion lacks details on the handling of
dynamic adaptivity.
The FAS and HTMG combination is explored in \cite{Weinzierl:09:Diss}, and
detailed for additive multigrid and BPX in \cite{Reps:15:Helmholtz}.
A combined, concise presentation for additive, BPX and multiplicative solvers
is new (cmp.~Table~\ref{table:solver-variants}). }

\new{Our augmentation of tree solvers with block smoothers stems from the dissertation
\cite{mweinzierl13diss} where block smoothers are solely applied to
multiplicative solvers.
The present manuscript generalizes them to
additive and BPX solvers, too.
\citet{mweinzierl13diss} also introduces the fusion of the tree with
Galerkin operators and BoxMG---again solely for the multiplicative case and
lacking our mirroring/reference element idea which simplifies the coding.
The compression of solution data through a hierarchical transform is explored in
\cite{Bungartz:10:Precompiler} for multigrid and in \cite{Eckhardt:16:SPH} for
SPH codes. Its application to the operators is, to the best of our knowledge, new, 
though we reuse some the aforementioned mechanisms.
Fragments of our shared memory or distributed memory concepts
 first can be found in
\cite{Weinzierl:09:Diss} or \cite{mweinzierl13diss},
respectively.
We focus on the correlation with multigrid and a theoretical concurrency
analysis here, while technical details are subject of other publications 
\cite{Eckhardt:10:Blocking,Weinzierl:17:PeanoSoftware}.
}

%
%
Our studies concentrate on operators with the sparsity
pattern of d-linear discretization and inter-grid transfer operators.
Wider operators, which are reasonable from an HPC point of view
\cite{Ghysels:13:ModelMG} or mandatory if stronger 
solver ingredients are required, induce different memory access patterns.
The present ideas continue to apply but require additional work.
Besides that, our experimental data stem from a multigrid prototype which
is not tuned.
For real-world computations, real scalability and performance engineering is
mandatory, and it might turn out that it is reasonable to compromise between our
academic approach and severely optimised strategies and existing libraries.
\new{Manycore vectorization for example has been successfully demonstrated
\cite{Reps:15:Helmholtz} where we fuse our tree paradigm with batched BLAS
concepts \cite{Dongarra:16:Batched}.
In general, performance and scalability engineering for particular
problems requires tailored solutions.}
Finally, our experiments restrict to academic setups \new{challenging enough to}
uncover the approach's potential.
It is obvious that the application to real-world experiments introduces further
challenges such as more sophisticated boundary conditions.

Besides experimental and implementational limitations due to the manuscript's
scope, there are conceptual limitations to the presented family of solvers:
Our approach reduces the memory footprint and data exchange but sacrifices
compute power.
Roadmaps predict that this is a reasonable strategy for
next generation linear algebra \cite{ascac10exa}. 
However, the flops are not for free yet.
The most severe restriction is that our approach cannot tackle problems
where geometric coarsening cannot resolve coarse-scale effects (as convection) anymore or where the
discretization requires very strong smoothers.
As such phenomena typically arise for coarser discretizations, 
it will remain obligatory for these applications to apply algebraic or
direct solvers on coarse levels, even though our techniques are applied.
The proposed approach can \new{widen} this notion of coarseness:
Our solvers remain stable for way coarser convection operators than a pure
geometric approach. 

\section{Spacetrees and FAS on generating systems}
\label{sections:spacetrees}

\begin{figure}
  \begin{center}
    \includegraphics[width=0.3\textwidth]{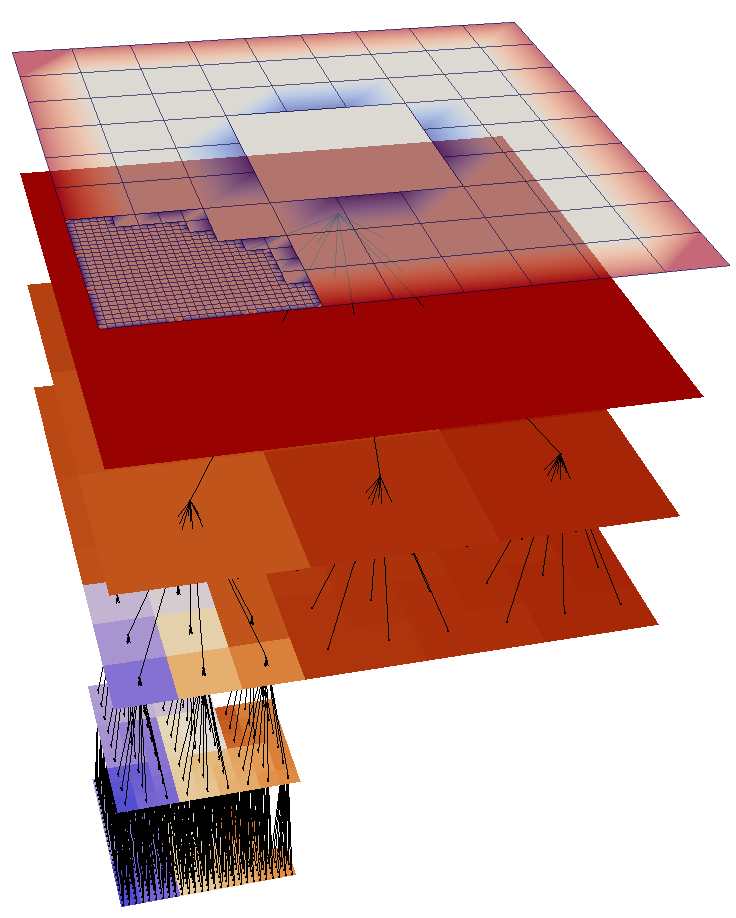}
    \hspace{1.2cm}
    \includegraphics[width=0.35\textwidth]{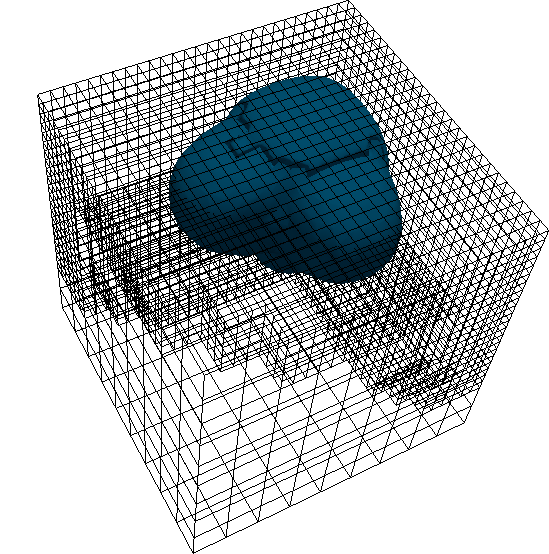}
  \end{center}
  \caption{
    Left: Spacetree for $d=2$ from \cite{Weinzierl:14:BlockFusion}.
    The top layer shows an adaptive Cartesian grid that results from a
    union of the individual levels of the tree (below). 
    Right: Grid for the three-dimensional convection-diffusion equation
    \new{in the $\texttt{checkerboard}$ setup} with an isosurface of the solution at $u =
    0.1$.
    \label{figure:spacetree}
  }
\end{figure}

Let a spacetree be a generalization of octrees or
quadtrees w.r.t.~the dimension $d$ and the cut cardinality $k\geq 2$ (Fig.~\ref{figure:spacetree}):
\new{
We embed $\Omega $ into a $d$-dimensional hypercube and cut it
equidistantly into $k$ pieces along each coordinate axis.
The original hypercube is the {\em root}.
It acts as parent node to the newly created $k^d$ smaller hypercubes that are
{\em children} of the {\em root}.
Root in turn is the {\em parent} of its children.
This process is repeated recursively.
Per cube we decide independently
whether to refine. 
The construction scheme yields a spacetree.
Each node of the tree graph represents one {\em cell}, i.e., a square ($d=2$)
or cube ($d=3$).
Let the {\em level} of a cell be the minimal number of construction steps we
need to create it. 
Root has level 0.
This usage of the term level results from a graph language.
}

For the present paper, we use $k=3$ as we base our experiments on the Peano
software. 
\new{We traverse the trees starting from the root which represents the
coarsest grid and run through them cell by cell.
Efficient storage concepts for such traversals are known (cmp.~Appendix
\ref{sections:appendix:spacetrees} which details the construction and traversal
of the grid and the underlying data structure). 
}

\new{Our present codes exploit}
the fact that a spacetree yields a cascade of ragged Cartesian grids, i.e.~each grid level defines vertices, 
\new{but} each level might cover a smaller part of the domain than the next
coarser one.
As a result, a vertex is unique through its position in space plus its level.
We distinguish three different vertex types: 
a vertex is {\em hanging} if it has less than $2^d$ adjacent cells on the same
level; a vertex is {\em refined} if there exists another non-hanging vertex at
the same position in space on a finer level, which implies that all adjacent
cells are refined further; and all other vertices are {\em unrefined}.

Our operators stem from a finite element discretization of
(\ref{equation:pde}) with $d$-linear shape functions.
This yields $3^d$-point stencils on regular grids.
Let each non-hanging vertex induce a shape function that spans all
adjacent cells of the same level.
The spacetree's shape functions then form a hierarchical generating system
\cite{Griebel:94:Multilevelmethoden}.
We combine the Full Approximation Storage
(FAS) \cite{brandt77mg,trottenberg01multigrid} scheme with the
hierarchical generating system and the idea of the Fast Adaptive Composite-Grid
Method (FAC) \cite{hart86fac, mccormick89fac}:

\begin{enumerate}
  \item Let a refined vertex hold the nodal value of all the vertices at the
  same spatial position with a higher level: $ u_\ell = I u_{\ell +1} $
  with $I$ being the point-wise injection. We copy the $u$ values of every vertex
  onto the $u$ values of coarser vertices for each vertex pair sharing
  the same spatial position.
  \item Let a hanging vertex's value be the $d$-linear interpolant
  from the coarser meshes.
  \item \new{Rely on} the same discretization technique on every level.
\end{enumerate}

\begin{figure}
  \begin{center}
    \includegraphics[width=0.65\textwidth]{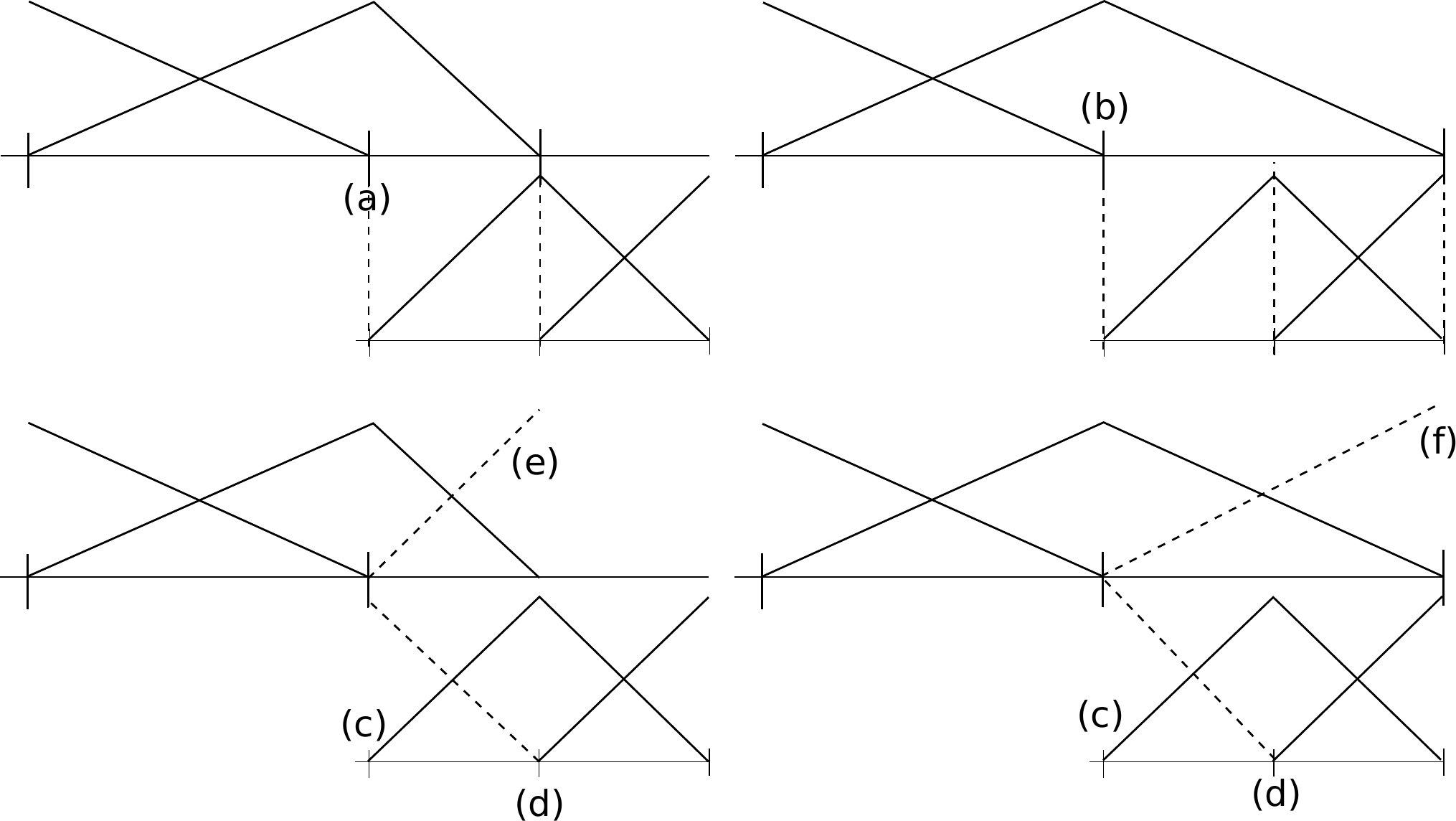}
  \end{center}
  \caption{
    Classic nodal shape functions yield discretization stencils for
    non-equidistant grids along the resolution boundary (in point (a) of
    left/top sketch), even if $\epsilon $ and $v$ are invariant. 
    We use a marginally hierarchical basis (top right) where the stencil in (b)
    can be computed such as any other stencil on this level. 
    No source code modification is required.
    We read the adaptive grid as overlapping domain decomposition with
    Dirichlet-Dirichlet coupling:
    In both variants, the hanging node acts on the
    fine grid as if there were an additional shape function halve centred in
    (c).
    No stencil is to be altered w.r.t.~coarse-fine transitions in point (d).
    The values in (e) or (f), respectively result from injection, i.e.~these
    points can be read as halved shape functions placed on the coarse grid,
    acting as Dirichlet points there, and thus coupling fine grid to coarse grid
    solution.
    In a multigrid context, (f) plays two roles: it acts as coarse system
    Dirichlet point and as left-most correction point on the coarser level.
    This induces a transition dilemma: the point carries the solution while it
    should be zero for the exact solution, e.g.
    With FAS, all correction system points carry a nodal representation, too,
    i.e.~this dilemma is resolved.
    \label{figure:fas}
  }
\end{figure}

\noindent
Smoothers then
can be read from a domain decomposition point of view, where
coarser grids prescribe the values at hanging nodes while fine-grid values
yield Dirichlet values in regions of the coarse grid overlapped by finer discretizations.
This renders the handling of hanging nodes and, more general, adaptivity
straightforward.
Notably, it implies that any discretization stencil can be unaware of
resolution transitions. 
Otherwise, the region where a fine grid transitions into a correction grid as it
is refined further requires special attention and additional coding (Figure
\ref{figure:fas}):
The semantics of the degrees of freedom change from
discretization weights into correction weights.
While the degrees of freedom at the transition have to carry the solution as
they act as boundary to the PDE solve on coarse mesh parts, their degree of
freedom weights have to approach zero once the solution starts to converge.
A FAS, i.e.~a scheme that starts to change coarse grid values from an injection
of fine grid solutions, removes this contradiction:
$u_\ell$ plays two different roles in adaptive meshes.
In unrefined regions, it carries a discretization of the PDE, while it
carries solution plus correction otherwise.
$A_\ell$ plays two roles, too.
In unrefined regions, it represents a discretization of the PDE, while it
encodes discretization plus correction term in refined regions.
The operator does not have to be changed at resolution transitions, 
which otherwise yields a large number of modified stencils even if we apply tree
balancing \cite{Sampath:08:Dendro} to reduce the number of possibilities of
coarse to fine cell configurations. 
While a coarsened/injected data representation is advantageous for 
non-linear problems, we consequently highlight a different advantage:


\begin{observation}
FAS allows us to ignore that some unknowns of one level carry a solution while
others have to carry a multigrid correction, as the latters' degrees of freedom
encode a correction term plus the coarsened solution.
No case distinction w.r.t.~the \new{semantics of the stencil or the unknown} is
required.
\end{observation}

\noindent
All discretization operators have strict local support,
as a stencil spans $3^d$ neighbouring vertices.
Let inter-grid transfer operators have local support, too.
A {\em parent} vertex $b$ of a vertex $a$ is any vertex that has at least
one adjacent cell that is a parent of an adjacent cell of $b$.
Interpolation and restriction couple vertices and parents only.
The notion of parent vertices here is slightly too generous.
The $d$-linear interpolation and restriction operators are actually even
sparser.

The spacetree's cascade of grids embeds function spaces into each other.
We propose to exploit this \new{beyond sole FAS} through the
HTMG idea \cite{Griebel:90:HTMM}. Let $A_\ell u_\ell  = b_\ell$
be the linear equation system for (\ref{equation:pde}) on a level $\ell
$.
The $u_\ell$ here are weights of the shape functions in $\Omega _{h,\ell}$ and
 comprise degrees of freedom on the fine grid mesh of level $\ell$.
 They also contain correction degrees of freedom if there are unrefined
 vertices on $\ell$.
The $b_\ell$ result from a discretization of $f$ if we are on the finest
level.
Otherwise, they comprise a multigrid correction component.
Instead of a correction equation, we use the Petrov-Galerkin coarse-grid
operator definition, switch to FAS and solve
\begin{eqnarray*}
   A_\ell (u_\ell+e_\ell) & = & R\left( b_{\ell+1}-A_{\ell+1}(id - P I) 
   u_{\ell+1}\right) 
   =: R\left( b_{\ell+1}-A_{\ell+1}\hat u_{\ell+1} \right) 
   =: R \hat r_{\ell+1}.
\end{eqnarray*}

\noindent
Here, $A_\ell$ is the level operator, $e_\ell$ or $r_\ell$ are the nodal
shape function correction weights and residuals, respectively, $P$ and $R$ are prolongation and
restriction, $id$ is the identity, $\hat r$ is called the hierarchical residual,
and $\hat u$ the hierarchical value.
Given an injected (a coarse) representation of the solution on a level $\ell $
through $I$, a multigrid solver computes a coarse-grid
update, keeps track of this coarse-grid solution modification $e_\ell$ and
prolongs this value back to the fine grid later all with one set of unknowns.

Our solver family relies on the following spacetree traversal paradigm:

\begin{definition}
 A {\em multiscale element-wise traversal} of the spacetree is a traversal of
 the cascade of grids $\Omega _{h,\ell}$ where
 \begin{enumerate}[topsep=0pt,itemsep=-1ex,partopsep=0.5ex,parsep=1ex,
leftmargin=0.5cm]
   \item each cell is processed only once, and the traversal offers throughout
   this processing access to all adjacent vertex data,
   \item each cell access allows for access to the cell's parent as well as its
   adjacent vertices,
   \item and a refined cell is processed after its children.
 \end{enumerate}
 \label{definition:multiscale-traversal}
\end{definition}

\noindent
Def.~\ref{definition:multiscale-traversal} gives a partial order on
the tree.
The order can be formalised as a set of operation applications $handleCell$. \new{These operation
applications (\textit{events}) are implemented as function calls that are
invoked at certain points during the spacetree traversal.} $handleCell$ implicitly introduces two further orders
$touchVertexFirstTime$
and $touchVertexLastTime$ on the tree's vertices that specify when a vertex is read for the first time and for the last time. 
Mirroring the statements from Def.~\ref{definition:multiscale-traversal},
the corresponding operations shall have access to their parent data, too.

We finally augment this set by a fourth transition:
$descend$ accepts a cell and its adjacent vertices as well as the $3^d$ children
of the cell with their vertices and precedes any $handleCell$ on any child. 
\new{
Such an event fits into a depth first traversal and does preserve the single
touch policy---every cell/vertex is read only once per multiscale grid
sweep---if we make the traversal code load within a refined cell first all $k^d$ children cells before is continues to recurse further (Figure
\ref{figure:events} in the appendix).
Traversals loading more than the direct children make some assumptions or
have to have knowledge about the grid structure.
$descend$ does not need this and thus fits to our strictly element-wise
mindset.
}

Definition~\ref{definition:multiscale-traversal} is a formalisation of our algorithmic ingredients
which directly fits to many standard ways to run through a spacetree, in particular
depth-first and breadth-first.
Depth-first traversals of spacetrees always resemble the
construction of space-filling curves as long as the run-through order of the
$k^d$ children of a parent is deterministic. 
Notably, Hilbert and Morton ($k=2$) and Peano ($k=3$) fall into place.

\section{Solver realizations}
\label{sections:solver}

A stencil of a vertex $v$ describes the entries of one row in $A_\ell $.
Such a row $(A_\ell)_v $ decomposes over the cells adjacent to $v$.
To compute $r_v=(A_\ell)_v u$, we can either compute $r_v$ by a sum over the
vertices, or we can split up $(A_\ell)_v $ additively over all the cells and
accumulate $r_v$ 
element-wisely. 
For $d=2$, the stencil 
\begin{equation}
 \left[
  \begin{array}{ccc}
   s_6 & s_7 & s_8 \\
   s_3 & s_4 & s_5 \\
   s_0 & s_1 & s_2
  \end{array}
 \right]
  \ \mbox{decomposes into}\ 
 \left[
  \begin{array}{ccc}
   s_6 & \frac{s_7}{2} & 0 \\
   \frac{s_3}{2} & \frac{s_4}{4} & 0 \\
   0 & 0 & 0
  \end{array}
 \right]
 +
 \left[
  \begin{array}{ccc}
   0 & \frac{s_7}{2} & s_8 \\
   0 & \frac{s_4}{4} & \frac{s_5}{2} \\
   0 & 0 & 0
  \end{array}
 \right]
 +
 \left[
  \begin{array}{ccc}
   0 & 0 & 0 \\
   0 & \frac{s_4}{4} & \frac{s_5}{2} \\
   0 & \frac{s_1}{2} & s_2
  \end{array}
 \right]
 +
 \left[
  \begin{array}{ccc}
   0 & 0 & 0 \\
   \frac{s_3}{2} & \frac{s_4}{4} & 0 \\
   s_0 & \frac{s_1}{2} & 0
  \end{array}
 \right].
 \label{equation:03_solver:stencils}
\end{equation}


\begin{table}
  \tbl{
    Overview of unknowns per vertex $v$ in the three geometric multigrid
    variants.
    \texttt{pers} indicates that the value is required
    in-between two solver iterations, i.e.~has to be stored persistently. \texttt{tmp} denotes
    that it is required only temporarily.
    \label{table:unknowns}
  }{
  {\footnotesize
  \begin{tabular}{l|p{8.8cm}|ccc}
    & Description  & Add & BPX & Mult \\
    \hline 
    $u$
     & Weight of shape function centred at vertex, i.e., function value
     in $v$.
     & pers & pers & pers \\
    $\hat u$
     & Hierarchical surplus in $v$, i.e.~difference to coarsened solution. 
     & tmp & tmp & tmp \\
    $r$
     & Residual.
     & tmp & tmp & tmp \\
    $\hat r$
     & Hierarchical residual.
     & tmp & tmp & tmp \\
    $d$
     & Helper variable transferring solver updates (deltas) between levels. 
     & pers & pers & pers \\
    $b$
     & Right-hand side.
     & pers & pers & pers \\
    $i$
     & Injected impact of the smoother. 
     & & pers \\ 
  \end{tabular}
  }
  }
\end{table}

\subsection{Geometric multigrid variants}
\label{sections:solver:geometric}

\noindent
The realization of a multiplicative\new{, geometric multigrid solver within 
strictly element-wise spacetree traversals} requires a state automaton $S$
steering the algorithm.
Let $S$ have two properties $current(S)$  and $old(S)$ that
identify the current smoothing level and the previous one. 
At startup, $current(S)=old(S)=\ell _{max}$. 
The automaton $S$ supports two transitions:
If we invoke $S \gets ascend(S)$, $current(S)$ is decreased. 
If we invoke $S \gets descend(S)$, $current(S)$ is increased.

\new{Both operations are controlled from a loop over the multigrid 
iterations. Within, the actual multigrid solver is a set of nested for loops
mirroring the well-known V- or W-cycle pattern.}
We study only $V(\mu _{pre},\mu _{post})$-cycles with $\mu _{pre} \geq 1$
and $\mu _{post} \geq 1$ here.
\new{From within these loops, }Algorithm
\ref{algorithm:geometric::multiplicative} is run once per multigrid cycle step.
If we realize a $V(\mu _{pre},\mu _{post})$-cycle, the algorithm's recursive
function is invoked $(\mu _{pre} + \mu _{post})(\ell _{max}-\ell _{min}+1)$
times.
The $S$ state transitions are invoked between two function calls.
We end up with three different state configurations.
For $current(S)=old(S)$, we perform either a pre- or a post-smoothing
step, which is not the first smoothing step on level $current(S)$.
For $current(S)=old(S)-1$, we run the first pre-smoothing step on
$current(S)$ and restrict the right-hand side from level
$old(S)$ to the next coarser level.
For $current(S)=old(S)+1$, we run the first post-smoothing step on
$current(S)$ and prolong corrections from level
$old(S)$ to the next finer level.

\begin{algorithm}[htb]
 \caption{
   Geometric multiplicative multigrid with rediscretization. It embeds into a
   multiscale element-wise spacetree traversal such as a depth-first ordering
   and is invoked by {\sc geomMult}($\ell _{max},S$) per smoothing step. $S$ is
   a state machine holding the current and previous active smoothing level. 
   The combination of these two levels distinguishes pre- from post-processing
   and triggers inter-grid data transfers.
 }
 \label{algorithm:geometric::multiplicative}
 {\footnotesize
  \begin{algorithmic}[1]
    \Function{geomMult}{$\ell, S$}
    
     \SIf{$\ell = current(S) \vee \left( current(S)>\ell \wedge
       vertex \ unrefined \right) $}
        $r_\ell \gets 0$
      \EndSIf 
      
     \SIf{$\ell = current(S) \vee \left( current(S)>\ell \wedge
       vertex \ unrefined \right) $}
         $\hat r_\ell \gets 0$
      \EndSIf 
      
     \SIf{$ current(S) < old(S) \wedge  old(S)=\ell $}
         $\hat u_\ell \gets u_\ell - Pu_{\ell -1} $
      \EndSIf 
      
     \SIf{$ current(S) < old(S) \wedge current(S)=\ell \wedge vertex \ refined $}
         $b_\ell \gets 0 $
      \EndSIf 
      
     \SIf{$ current(S) > old(S) \wedge current(S)=\ell $}
         $ u_\ell \gets  u_\ell + Pd_{\ell -1} $
      \EndSIf 
      
     \SIf{$ current(S) > old(S) \wedge current(S)=\ell $}
         $ d_\ell \gets d_\ell + Pd_{\ell -1}  $
      \EndSIf 
      
     \SIf{$\ell < \ell_{max} \wedge \ell < max\{ current(S), old(S) \}$}
        \Call{geomMult}{$\ell +1$,S}
     \EndSIf
     
     \SIf{$ current(S) = \ell \vee current(S) < \ell \wedge cell\ unrefined$}
         $ r_\ell \gets -A_\ell u_\ell  $
      \EndSIf 
      
     \SIf{$current(S) < old(S) \wedge \ell = old(S) $}
         $ \hat r_\ell \gets -A_\ell \hat u_\ell  $
      \EndSIf 
      
      \SIf{$ current(S) = \ell \vee current(S) < \ell  $}
         $ r_\ell \gets r_\ell + b_\ell $
      \EndSIf      
      
      \SIf{$ current(S) < old(S) \wedge \ell =  old(S) $}
         $\hat r_\ell \gets  \hat r_\ell + b_\ell $
      \EndSIf 
      
      \SIf{$  \ell = current(S) \vee ( current(S)>\ell \ \wedge vertex \
      unrefined $} 
       \State  $\phantom{if} u_\ell \gets  u_\ell + \omega  \
       diag^{-1}(A_\ell) \ r_\ell $
      \EndSIf 

     \State $d_\ell \gets \left\{
      \begin{array}{ll}
       0 & \mbox{if}\ \ell = current(S) \wedge vertex \ refined \\
       d_\ell + \omega  \ diag^{-1}(A_\ell) \  r_\ell & \mbox{if}\ \ell =
       current(S) \vee ( current(S)>\ell \ \wedge \\
         & vertex \ unrefined ) \\
       d_\ell & \mbox{otherwise}
      \end{array}
      \right.$
      
      \SIf{$ current(S) \geq \ell $}
         $u_{\ell-1} \gets  Iu_{\ell} $
      \EndSIf 
      
      \SIf{$ current(S)<old(S) \wedge old(S)=\ell  $}
         $b_{\ell-1} \gets R\hat r_{\ell} $
      \EndSIf 
      
     \EndFunction
  \end{algorithmic}
  }
\end{algorithm}

Cells contribute to the residual either if they are on the active level or if
they belong to a level that is coarser than the active level and are adjacent to
at least one unrefined vertex.
\new{Such} a multiscale smoothing can be read in a domain decomposition
way, where coarse-grid values are overwritten by overlapping fine-grid
values.
Many multiplicative multigrid codes coarsen all grid regions
when they ascend.
The present code, in contrast, coarsens the levels that are finer than the previous smoothing level.
This makes the code easier to understand. 
A fine-grid region then is subject to the more smoothing steps the coarser it
is.

\begin{observation}
We can implement matrix-free, {\em geometric multiplicative multigrid}
within an element-wise multiscale traversal 
with one smoothing step per grid sweep.
\end{observation}

\begin{algorithm}[htb]
  \caption{
   Geometric additive multigrid based upon rediscretization. It embeds into a
   multiscale element-wise spacetree traversal such as a depth-first ordering
   and is invoked by {\sc geomAdd}($\ell _{max}$). If all hanging
   nodes are made to hold the $d$-linear interpolant of the next coarser levels,
   the code works on arbitrary adaptive meshes.
 }
 \label{algorithm:geometric::additive}
 {\footnotesize
  \begin{algorithmic}[1]
    \Function{geomAdd}{$\ell $} 
     \State $d_\ell \gets d_\ell + Pd_{\ell-1}$; $u_\ell \gets u_\ell +
     Pd_{\ell-1}$; $\hat u_\ell \gets u_\ell - Pu_{\ell -1}$
       \Comment{Vertex-wise \new{correction}}
     \State $r_\ell \gets 0$; $\hat r_\ell \gets 0$
       \Comment{\new{(coarse-to-fine data flow)}}
     \State $b_\ell \gets 0$ for all $b$ associated to refined vertices
     \SIf{$\ell < \ell_{max}$}
        \Call{geomAdd}{$\ell +1$}
     \EndSIf 
       \Comment{Recursion into finer grid levels}
     \State $r_\ell \gets -A_\ell u_\ell $ ; $\hat r_\ell \gets -A_\ell \hat
     u_\ell $
       \Comment{Cell-wise residual accumulations} 
     \State $r_\ell \gets r_\ell + b_\ell $; $\hat r_\ell \gets \hat r_\ell +
     b_\ell $
     \State $d_\ell \gets \omega \ diag^{-1}(A_\ell) \  r_\ell $
       \Comment{Vertex-wise smoothing}
     \State $u_\ell \gets u_\ell + d_\ell $
     \SIf{$\ell > \ell _{min}$}
       $u_{\ell-1} \gets I u_{\ell}$; $b_{\ell-1} \gets R \hat r_{\ell}$
     \EndSIf
          \Comment{Vertex-wise \new{restriction and} injection}
     \EndFunction
  \end{algorithmic}
 }
\end{algorithm}

\noindent
\new{ In \cite{Reps:15:Helmholtz}, two additive solver alternatives are
presented that work in an element-wise multiscale traversal and both require, amortized, one multiscale grid sweep per additive cycle. The versions of these solvers that we use in this work
are Algorithm \ref{algorithm:geometric::additive} and Algorithm \ref{algorithm:geometric::bpx}.
Algorithm \ref{algorithm:geometric::additive} is a classical additive scheme, and Algorithm \ref{algorithm:geometric::bpx} is a BPX variant. 
A detailed describtion these algorithms and their implementation can be found in Appendix \ref{appendix:geoaddMG}. 
}
\begin{observation}
We can implement a matrix-free, {\em geometric additive multigrid
\linebreak solver} (a solver based upon $d$-linear inter-grid operators and
rediscretization of coarse grid matrices) within an element-wise multiscale traversal that requires,
amortized, one multiscale grid sweep per additive cycle.
\end{observation} 

\begin{observation}
We can implement a matrix-free, {\em geometric BPX solver} 
within an element-wise multiscale traversal that requires, amortized,
one multiscale grid sweep per additive cycle.
\end{observation}

\begin{algorithm}[htb]
 \caption{
   Geometric BPX variant with rediscretization. It embeds into a
   multiscale element-wise spacetree traversal such as a depth-first ordering
   and is invoked by {\sc geomBPX}($\ell _{max}$).
 }
 \label{algorithm:geometric::bpx}
 {\footnotesize
  \begin{algorithmic}[1]
    \Function{geomBPX}{$\ell $} 
     \State $d_\ell \gets d_\ell + \left\{ 
      \begin{array}{ll}
       Pd_{\ell-1} & \mbox{if c-point} \\
       Pd_{\ell-1} - Pi_{\ell-1} & \mbox{else}
      \end{array}
      \right.$
       \Comment{Vertex-wise operations from \textsc{geomAdd}}
     \State $u_\ell \gets u_\ell + \left\{
      \begin{array}{ll}
       Pd_{\ell-1} & \mbox{if c-point} \\
       Pd_{\ell-1} - Pi_{\ell-1} & \mbox{else}
      \end{array}
     \right.$
       \Comment{with modified prolongation}
     \State $\hat u_\ell \gets u_\ell - P_{\ell -1}$
     \State $r_\ell \gets 0$; $\hat r_\ell \gets 0$
     \State $b_\ell \gets 0$ for all $b$ associated to refined vertices
     \SIf{$\ell < \ell_{max}$}
       \Call{geomBPX}{$\ell +1$}
     \EndSIf
     \State $r_\ell \gets -A_\ell u_\ell $; 
            $\hat r_\ell \gets -A_\ell \hat u_\ell $
     \State $r_\ell \gets r_\ell + b_\ell $;
            $\hat r_\ell \gets \hat r_\ell + b_\ell $
     \State $d_\ell \gets \omega  \ diag^{-1}(A_\ell) \  r_\ell $
     \State $u_\ell \gets u_\ell + d_\ell $ for non-c-points
     \SIf{$\ell > \ell _{min}$}
         $b_{\ell-1} \gets R \hat r_{\ell}$;
         $i_{\ell-1} \gets I d_{\ell}$
     \EndSIf
      \Comment{Injection of skipped updates}
     \State $d \gets 0$ for c-points
      \Comment{Skip update}
     \EndFunction
  \end{algorithmic}
 }
\end{algorithm}


\subsection{Spacetree block smoothers}
\label{sections:blocksmoothers}

Strict element-wise matvecs that touch each (fine-grid) cell only once render
smoothers beyond point Jacobi technically challenging if a whole smoothing step
has to be realized within one grid traversal.
We can update an unknown only once a vertex is used for the last time (in
(\ref{equation:03_solver:stencils}) once the fourth cell has been processed),
and information of a vertex update thus can propagate only along the grid
traversal sequence:
A vertex update may only affect vertices that are adjacent to the last cell
processed that is adjacent to this vertex.
We observe that a naive splitting of a stencil into equal parts as realised in
(\ref{equation:03_solver:stencils}) then is not possible anymore and the
splitting has to anticipate the Gauss-Seidel like enumeration.

Gau\ss -Seidel with an unknown enumeration that is not tied to the cell
traversal order or line smoothers can not be realized within one grid
traversal.
Coloured schemes such as red-black Gauss-Seidel require one grid traversal per colour.
Krylov schemes work if we evaluate the matvec as well as all scalar products
in one grid sweep and apply the impact in a second sweep \cite{Bader:08:MemorySierpinski}.
With pipelining, multiple sweeps can be fused and the amortized cost per unknown 
update can be reduced
\cite{Ghysels:13:ModelMG}.
In the present work, we however restrict ourselves to multigrid
ingredients with minimalist memory access and thus stick to Jacobi.

Point Jacobi is a poor choice for many non-trivial parameter combinations in
(\ref{equation:pde}).
To facilitate more powerful smoothers without giving up data locality or single
touch, we augment Jacobi by an additional block smoothing that improves
convergence locally on very small subdomains ($k^d$ patches).
For this, we generalize the tree traversal by a $descend$ event (see
Sect.~\ref{sections:spacetrees}).
In a depth-first traversal code, such an operation makes a recursive step down
within the tree and loads all children of a node before it continues
recursively---a one-level recursion unrolling
\cite{Eckhardt:10:Blocking}.

\begin{figure}[htb]
  \begin{center}
    \includegraphics[width=0.3\textwidth]{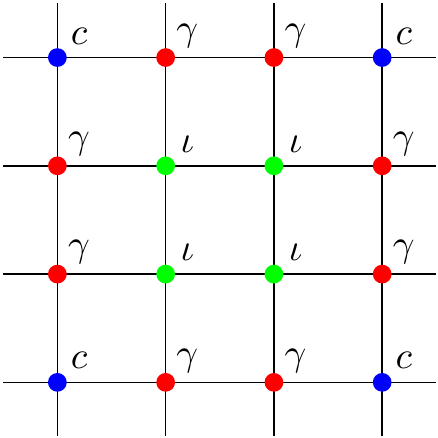}
    \hspace{1.0cm}
    \includegraphics[width=0.3\textwidth]{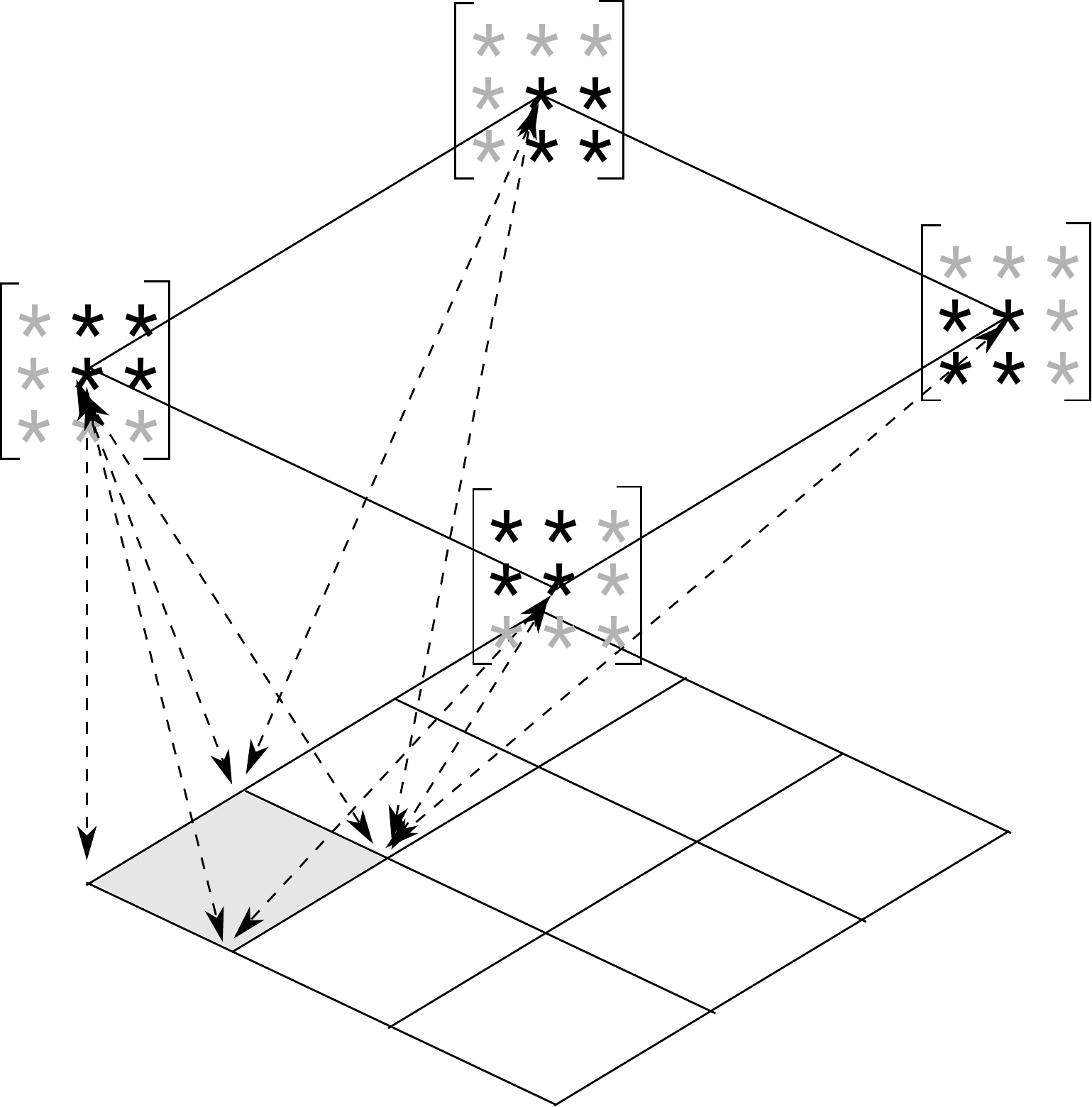}
  \end{center}
  \caption{
    Left: $3^d$ patch available to $descend$ (together with the parent cell). 
    We distinguish three vertex types: c-points coincide with 
    vertices on the next coarser grid level, $\gamma $-points lie on coarse-grid
    lines, and $\iota $-points lie within the patch.
    Right: Illustration of a cell-wise Galerkin coarse grid operator
    computation.
    With $P$ and $R$ known, we can set the interpolation for the vertices of a cell, apply
    the element matrix and add the result back to the coarse cell's stencils
    via $R$. Coarse-grid stencil entries subject to accumulation
    are bold.
    \label{fig:gridpoints}
  }
\end{figure}

A block smoother accepts all $(k+1)^d$ vertices and runs the solver's correction
steps for them if not triggered by adjacent cells already.
We distinguish three different types of \new{vertices
(Fig.~\ref{fig:gridpoints})---the nomenclature follows the BoxMG terminology as
detailed in Sect.~\ref{sections:boxmg}---and} we use the term {\em block Jacobi}
for a smoother that performs Jacobi on the $\gamma$ and $c$ points and processes the $\iota$ vertices differently. 
When block Jacobi loads a patch's $(k-1)^d$ $\iota$-vertices 
for modification (clearance of right-hand sides, e.g.) it first runs Gau\ss -Seidel sweeps on them.
Then it computes the respective operations from
the algorithm, determines the hierarchical surplus and continues in an 
element-wise fashion.

While an exact solve on the interior points of a patch would be possible and
even trivial for $k=2$, we use $k=3$ and find it convenient to make the
patch sweeps run Gau\ss-Seidel iterations.
For the geometric multigrid solver, these iterations use on-the-fly
rediscretization.
For the Galerkin coarse grid variants using $d$-linear or BoxMG
inter-grid transfer operators we make them use stencils that are explicitly
available.
The proposed technique falls into the class of hybrid smoothers
\cite{baker11ultraparallelsmoother}.

\begin{observation}
Through an augmentation of the multiscale element-wise spacetree traversal with
a $descend$ operation, the realization of block Jacobi smoothers is
straightforward. They preserve the data locality of the element-wise traversal.
\end{observation}

\subsection{Galerkin multigrid variants}
\label{sections:galerkin}

Galerkin in our tests denotes all mul\-tigrid variants
where prolongation and restriction are $d$-linear, while coarse grid operators
result from $A_{\ell -1} = R A_{\ell } P$.
\new{The} Galerkin computation rules  
do not imply any assumption about $P$ and $R$.
As a computation of  $A_{\ell -1}$ depends on a cascade of evaluations on finer
levels, an on-the-fly computation of $A_{\ell -1} $ is impossible.
The coarse-grid system has to be held explicitly.
For this, we augment the vertex records with $3^d$ (for multiplicative
multigrid) or $2 \cdot 3^d$ (for additive and BPX) doubles. 
They hold the stencil associated with a vertex that in turn determines the
element-wise matrices.
All $r \gets -A_{\ell} u$ and $\hat r \gets
-A_{\ell} \hat u$ evaluations are modified such that they read the
stencils from the cells' adjacent vertices.
Whenever we create a new vertex, the stencil entries are initialized via
PDE discretization.
To determine the Galerkin operator, we accumulate the coarse-grid operator
element-wisely together with the residual by decomposing $A_{\ell -1}|_v = (R
A_{\ell } P)|_v$ over all $2^d \cdot 3^d$ child cells $\hat c$ of the cells adjacent to $v$
(Fig.~\ref{fig:gridpoints}, right).
This strategy is element-wise w.r.t. level $\ell-1$.

To make the accumulation work, we have to 
clear \new{the} stencils \new{before}.
If a grid is stationary and the PDE is linear, we could skip any re-accumulation
of coarse-grid operators.
If the grid \new{however refines} into a level $\ell +1$, fine grid stencils on
a level $\ell $ that are associated to newly refined vertices transform into stencils that
carry a Galerkin operator and thus have to be recomputed.
This recomputation \new{recursively} triggers stencil updates on coarser levels.
If the PDE is non-linear, the stencils on the finest level depend on the current
approximate solution $u_h$ and thus trigger changes recursively on all coarser
levels as soon as we update $u_h$.
Even if we stick to linear problems such as (\ref{equation:pde}), it is convenient to recompute the coarse-grid
operators in every cycle before we coarsen in the spacetree.
We then do not have to analyze whether the grid changes.
Our vertical integration of multigrid operations suggests that all fine-grid
operators influencing a Galerkin recomputation are held in caches.
The re-accumulation does not increase the pressure on the
memory subsystem.
\new{Finally}, we could omit the storage of the stencils on the finest grids
and rely on on-the-fly rediscretization.
Yet, we introduce a holistic memory compression
applying to all grid levels in Sect.~\ref{sections:stencil-compression} that
realises this storage optimisation automatically.

If we recompute the Galerkin operator in each traversal, we need both a valid
operator and memory to accumulate the coarse grid stencils for the additive solver variants.
Additive multigrid and the BPX variant thus store a copy of the stencil upon
the first read of a vertex in the data structure.
Hereafter, all current stencil entries are cleared and we start the
accumulation.
Matvecs use the backup copy, which can be held temporarily,
i.e., not stored in-between two grid sweeps.
A Galerkin variant of \new{the multiplicative} Algorithm
\ref{algorithm:geometric::multiplicative} does not require duplicated stencils.
Let 
\begin{equation}
recomputeGalerkin(\ell, S) = \left\{ 
  \begin{array}{ll}
   \top & \mbox{if}\ \ell+1=current(S) \wedge \mbox{last smoothing step} \\
   &    
   \mbox{on level } \ell+1 \ \wedge \ vertex \ refined
   \\
   \bot & \mbox{otherwise.}
  \end{array}
 \right.
 \label{equation:recomputeGalerkin}
\end{equation}

\noindent
A vertex's stencil is set to zero when we read the
vertex for the first time if 
$recomputeGalerkin=\top$.
Coarse operator contributions are added to those coarser
vertices where the predicate holds.
They are not altered on finer grids than the active smoothing level or
during the descend process.

\begin{observation}
  If we augment each vertex data structure by $3^d$ (multiplicative) or
  $2 \cdot 3^d$ (additive and BPX) doubles, we can realize Galerkin multigrid
  variants within the element-wise multiscale traversal that work on-the-fly.
\end{observation}

\noindent
A proper choice of $\ell _{max}$ is delicate for multigrid methods:
If $\ell_{max}$ is too fine, multigrid convergence suffers.
If $\ell_{max}$ is too coarse, the Galerkin operators deteriorate, might become
indefinite, and, in the worst case, the coarse grid smoothing contributions
destroy the overall convergence and make the solver diverge \cite{yavnehbendig12nonnsymbb}.
We propose to rely on a dynamic coarsest level which starts with $\ell
_{max}=1$.
On the one hand, the solver increases $\ell _{max}$ in-between two multigrid
cycles/iterations once we observe stagnating or \new{growing residuals}.
On the other hand, the solver immediately increases  $\ell _{max}$ if one
Galerkin operator computation yields a coarse grid stencil that is not diagonal
dominant anymore or carries negative diagonal values.
This approach makes our solver variants start as a multilevel approach and
deteriorate, in the worst case, to (block) Jacobi.
We note that the definiteness check is, in principle, active only at the solver
startup, if the underlying PDE is linear.
The convergence speed criterion however can kick in later.
For more sophisticated applications, it might be reasonable to make $\ell
_{max}$ space-dependent\new{, too}.
We do not follow-up such a sophisticated scheme.

\subsection{BoxMG}
\label{sections:boxmg}

To construct PDE-depen\-dent inter-grid
transfer operators, we rely on BoxMG \cite{dendy82blackbox}
\new{applied to tri-partitioning}
\cite{yavnehbendig12nonnsymbb,mweinzierl13diss,dendy10blackboxCoarseBy3}. 
It has been shown to yield robust and efficient multigrid solvers for a large
class of problems while it can be seen as a special case of classical algebraic
multigrid with a geometric definition of ``strong connections'' \cite{maclachlan12robustAdaptiveMg}.
It thus fits to our geometric multiscale meshing concept.
For studies on the robustness and efficiency of BoxMG, we refer  to 
\cite{dendy83blackboxNonsym,
dendy88boxmgperiodicsingular,moultonetal98boxmghogenization}\new{, and notably
cite} \cite{wienandskoestler10frameworkprolongation} \new{offering} a
framework for the construction of prolongation operators where BoxMG is recovered as a special case of more general multigrid
techniques.

For any refined vertex, this vertex's inter-grid transfer operator affects $5^d$
fine grid vertices.
It carries a $5^d$ stencil for $P$ and $R$.
\new{Instead of discussing how the inter-grid stencils affect the $5^d$ fine
grid vertices, we 
\begin{itemize}
  \item clear $P$ and $R$ upon a coarse grid vertex load,
  \item plug into $descend$ when we descend from an refined cell that is 
  adjacent to the vertex into the finer grids,
  \item collect there the $4^d$ stencils from the fine grid vertices, and
  \item alter the $3^d$ affected stencil entries of the coarse grid,
\end{itemize}
whenever these are to be recomputed.
}
In \new{exchange, $descend$  always} computes $2^d$ partial inter-grid transfer
stencils per cell, \new{i.e.~contributes to the stencils of all the vertices
adjacent to a cell. 
We fit to the concept of element-wise assembly for the multigrid.
}

\new{As we work with cubes only, any vertex/stencil configuration can be
mirrored by a matrix $M$ such that the coarse vertex and its stencil of interest
coincide with the left, bottom vertex.
$M$ reorders the vertices within a $3^d$ patch an all of their stencils, too.
Once we have determined all $P$ entries of interest---$R$ follows due to
$R^T=P$ if not explicitly stated otherwise---we can mirror all entries back
through $M^T$.
We break down BoxMG into multiscale element-wise operations that we
formalise w.r.t.~one vertex configurations from which $2^d-1$ further
operations results through mirroring.
}

Our presentation restrict\new{s} to one \new{partial inter-grid transfer} stencil \new{therefore}.
\new{Following the literature, we} distinguishe c-, $\gamma $- and $\iota
$-points (Fig.~\ref{fig:gridpoints}) \new{within a $3^d$ patch}
and re-order the equation system accordingly.
All formalism and techniques are here described for $d=2$ and extend naturally
to the three-dimensional case.
For the complete \new{operators}, we refer to Appendix
\ref{appendix:boxmgtensors}.
Let the stencil $P$ be equivalent to the vector $P = (P_c\ P_\gamma \ P_\iota
)^T$.
In $2d$, the entries read
\[
  P = 
  \left[
   \begin{array}{ccccc}
     p_{0,4} & p_{1,4} & p_{2,4} & p_{3,4} & p_{4,4} \\
     p_{0,3} & p_{1,3} & p_{2,3} & p_{3,3} & p_{4,3} \\
     p_{0,2} & p_{1,2} & p_{2,2} & p_{3,2} & p_{4,2} \\
     p_{0,1} & p_{1,1} & p_{2,1} & p_{3,1} & p_{4,1} \\
     p_{0,0} & p_{1,0} & p_{2,0} & p_{3,0} & p_{4,0} 
   \end{array}
  \right]
  \
  \mapsto
  \
  (
    \underbrace{p_{2,2}}_{P_c},
    \underbrace{p_{3,2},p_{4,2},p_{2,3},p_{2,4}}_{P_\gamma},
    \underbrace{p_{3,3},p_{4,3},p_{3,4},p_{4,4}}_{P_\iota}
  )^T
\]
\noindent
with lexicographic stencil entry enumeration for the prolongation stencil
associated to the bottom left coarse grid vertex of a $3\times 3$ patch
(Fig.~\ref{fig:gridpoints}).
BoxMG makes the impact of a coarse-grid correction $u_{\ell -1}$ onto 
$u_{\ell }$ fall into the PDE's nullspace.
Its operator-dependent inter-grid transfer implies that an interpolation of
coarse data fits to the PDE.

 \begin{equation}
 A_\ell P u_{\ell -1} = 
 \left(
  \begin{array}{ccc}
   A_{cc} & A_{c\gamma } & A_{c\iota } \\
   A_{\gamma c} & A_{\gamma \gamma } & A_{\gamma \iota } \\
   A_{\iota c} & A_{\iota  \gamma } & A_{\iota \iota } 
  \end{array}
 \right) 
 \left(
  \begin{array}{c}
   P_c \\
   P_\gamma \\
   P_\iota
  \end{array}
 \right)
 u_{\ell -1} = 
 \left(
  \begin{array}{c}
   b_c \\
   b_\gamma \\
   b_\iota
  \end{array}
 \right),
 \label{equation:boxmg}
 \end{equation}
 \noindent
\new{$p_{2,2}$=1} induces $5^d-1$ interpolated fine grid
 values that disappear under the PDE operator.
 To achieve this, $P = (P_c\ P_\gamma \ P_\iota )^T$ is constructed in five
 steps:

 \begin{enumerate}
 \item c-points are assigned the value of their coarse counterpart coinciding
 spatially. $P_c = I^T$ with $I$ from FAS.
 \item We ignore the impact of $\gamma $- and $\iota $-points on
 c-points and from $\iota $-points on $\gamma $-points. $A_{c
 \gamma } = A_{c \iota } = A_{\gamma \iota } = 0$ and, therefore, $A_{cc} =
 id$. 
 We bring $A_{\gamma c}$ and $A_{\iota c}$ to the right-hand side and obtain 
\[
 \left(
  \begin{array}{cc}
   A_{\gamma \gamma } & 0 \\
   A_{\iota  \gamma } & A_{\iota \iota } 
  \end{array}
 \right) 
 \left(
  \begin{array}{c}
   P_\gamma u_{\ell -1}\\
   P_\iota  u_{\ell -1}
  \end{array}
 \right)
 =
 \left(
  \begin{array}{c}
   b_\gamma - A_{\gamma c } P_c u_{\ell -1} \\
   b_\iota  - A_{\iota c } P_c u_{\ell -1}
  \end{array}
 \right).
 \]
 \item This system remains hard to solve as the matrices are large.
 BoxMG therefore decomposes the level $\ell $ into patches (Fig.~\ref{fig:gridpoints}).
 To reduce inter-patch dependencies, the two-dimensional 
 stencils belonging to $\gamma$ points are collapsed to one-dimensional 
 stencils by summing up all stencil entries in the dimension perpendicular to
 the corresponding coarse grid line. 
 In $d=2$ and for coarsening by a factor of three, each two $\gamma$ points on a
 coarse-grid line can be computed from the two neighbouring $c$ points and
 themselves by solving two equations in two unknowns:
 \[
 P_\gamma u_{\ell -1} = \tilde A_{\gamma \gamma}^{-1} ( b_\gamma - \tilde
 A_{\gamma c } P_c u_{\ell -1} ).
 \]
 \item  As multigrid is defined over residual equations, it is
 reasonable to assume $b_\gamma = b_\iota =0$. This yields a linear equation for $P_\gamma$. 
 More efficient BoxMG variants apply a postsmoothing step similar to 
 smoothed aggregation to $P$ and do not neglect the right-hand side
 in (\ref{equation:boxmg}) \cite{maclachlan12robustAdaptiveMg}.
 We do not follow-up this technique though our software base in principle
 allows for nonhomogeneous right-hand sides.
 \item Finally, the four $\iota$ points are
 computed by solving four equations in four unknowns.
 \end{enumerate}

\begin{observation}
BoxMG yields operator-dependent inter-grid transfer operators that can be represented by
$5^d$-stencils per vertex.
We thus can realize an algebraic-geometric multigrid solver without any external
global matrix if we store these stencils within the vertices.
\end{observation}

\noindent 
The extension of the scheme to three dimensions is straightforward
\cite{Behie1983boxmg3d, Scott1985boxmg3d,dendy86boxmg3d,yavnehbendig12nonnsymbb}:
The stencils are collapsed into $1d$ stencils along patch cube edges and 
into $2d$ stencils on the patch faces. An $8 \times 8$ equation system
is to be solved in the patch interior.
\new{We summarise the key property of BoxMG's construction from a traversal's
point of view: Whenever the tree traversal descends within a refined cell, it
alters exactly $3^d$ entries of any $P$ stencil of any adjacent coarse
stencil.
For this, it has to know the $3^d$ stencils of the fine grid vertices affected.
Within a patch however no $P$ entry depends on any stencil of any vertex that is
not contained within the same patch.
}


\begin{observation}
 \new{
 We can implement BoxMG in a strict multiscale element-wise sense through the
 introduction of $descend$.
 All inter-grid transfer operators of one level become available within one grid
 sweep. We propose a single-sweep, single-touch inter-grid transfer operator
 computation.
 }
\end{observation}

\noindent
\new{
One advantage of additive multigrid variants is the possibility to merge coarse
grid operator computations with the smoothing process
(\cite{AlOnazi:17:TaskBasedAlgebraicMG} and references therein).
Our approach offers this propery also for the multiplicative variant:
Coarse-grid operator computation, level $\ell -1$ smoothing and restriction
from level $\ell +1$ all are interwoven.
}

Plugging a symmetrized system operator $A_{sym} = \frac{1}{2}(A_{\ell} + A_{\ell}^T)$
into the computation of either restriction or prolongation computation in
the \new{BoxMG scheme} improves the convergence and robustness for non-symmetric setups \cite{dendy83blackboxNonsym, yavnehbendig12nonnsymbb}.
Such a Petrov-Galerkin multigrid scheme however changes the memory access
pattern for $P$:
As ``half of $A_{sym}$'' stems from $A_{\ell}^T$, a patch's computation of
$P$ entries needs the stencils of vertices surrounding a patch.
It effectively uses a patch with $5^d$ cells, i.e.~$6^d$ vertices per linear
equation system solve.
It currently is unclear whether we can avoid a spreading of the influence area and stick to
a strictly patch-wise, localized data evaluation pattern and,
at the same time, use a symmetrized operator.
This has to be subject of future studies.
Yet, we can use simpler symmetric restriction operators $R$ in combination with
the BoxMG prolongation $P$. 
Notably, simple injection $R=I$ \cite{Griebel:90:HTMM} or aggregation of $5^d$
fine grid points into one coarse grid point are trivial to implement.
The resulting inter-grid operator combination then lacks the accuracy required by
multigrid efficiency proofs, but we know that the resulting Galerkin
coarse grid operators tend to remain more stable---they do not degenerate
that fast into central differences for convection-dominated flows when we
coarse \cite{yavneh98singularPertubation}.
We thus can expect that our adaptive coarse grid selection does not increase
$\ell _{max}$ as fast as for traditional BoxMG.

The computation of $P$ and $R$ fits to $descend$, and
we can store the results as stencils in the vertices.
Again, we need (temporary) backups of the inter-grid operators in the additive
variants.
All statements and details on (re-)accumulation of stencils for the coarse-grid
operators apply to the BoxMG-operators, too.

The patch-based locality makes BoxMG well-suited
for spacetree-based non-uniform grids.
The interpolant on hanging vertices is no longer determined geometrically.
It instead results from $P$ and $R$ according to the BoxMG formalism.
At the hanging vertices, well-suited stencils that can be
collapsed are required.
We use $d$-linear interpolation of the parent stencils to obtain them.
The elimination of dependencies through stencil collapsing along patches
preserves the high memory access locality of an augmented element-wise multiscale
traversal.

\section{Stencil compression}
\label{sections:stencil-compression}

Our Galerkin and BoxMG multigrid variants are not literally matrix-free.
They do not hold a dedicated matrix data structure, but they store the
stencils within the grid.
Since sparse matrix storage formats exist
that introduce a small administrative overhead \cite{King:16:DynamicCSR},
the savings through this in-situ storage are limited.
Yet, significant savings can be made if we omit storage on the finest grid
levels and recompute the stencils there on-the-fly.
This does not introduce any savings on coarser grids.
\new{If} the discretization
is costly---through material parameters $\epsilon $ that have to be integrated
or challenging boundary conditions, for example---it might be better for
performance\new{-wisely} to store the stencils on the finest grid level, too.

Galerkin coarse-grid operators resemble
rediscretizations for smooth $\epsilon$, small $v=0$ or fine mesh sizes.
We thus introduce hierarchical operators
\begin{eqnarray*}
  \hat A & = & A - A_{rediscretized}, \quad  
  \hat P = P - P_{d-linear} \quad \mbox{and} \quad
  \hat R = R - R_{d-linear}.
\end{eqnarray*}

\noindent
$\hat A, \hat P$ and $\hat R$ can be computed whenever we use a vertex
for the last time throughout a grid traversal.
Either $A$ or $\hat A$ have to be held in-between two
iterations, i.e., we can either store the original operator or reconstruct it
from the hierarchical representation upon the subsequent vertex
load.
The argument holds for $\hat P$ and $\hat R$ analogously.

\begin{observation}
  For setups where $\epsilon $ is smooth and $v$ is
  small in most of the domain, the entries of the hierarchical operators $\hat
  A$, $\hat P$ and $\hat R$ are small.
  They hold fewer valid digits than made available through the IEEE standard.
\end{observation}

\noindent
For an almost matrix-free multigrid realization, we thus propose \new{to}
rewrite all three operators held within the vertex into their hierarchical representation.
If the operators are zero, i.e., the stencil equals rediscretization and the 
inter-grid transfer operators are $d$-linear we mark the vertex
and discard the operator's stencil.
Otherwise, we convert all entries $x$ of the hierarchical representation into a
format $f^{-1}_{bpa}(x) = m \cdot 2^e$, where the exponent $e$ is stored in one
byte (C data type \texttt{char}) and $m \in \mathbf{N}_0$, with
$e$ chosen such that $m$ fits exactly into $bpa-1$ bytes as a natural
number.
Here, $bpa \in \{0,2,\ldots,8\}$ (bytes per attribute) is the number of
bytes that we use to store the exponent $e$ plus the integer value $m$.
Upon a vertex store, we determine per operator the smallest $bpa$ such that
$|f_{bpa}(\hat x) - \hat x| \leq \epsilon _{mf}$ with $\epsilon _{mf} \ll 1$.

Within the vertex, solely $bpa$ per stencil is held in-between
two iterations. 
All three $bpa$ values \new{for $A$, $P$ and $R$} fit into 9 bits.
The values $e$ and $m$ per stencil entry are piped into a separate byte stream.
When we read the vertex for the first time, we take
$bpa$ from the vertex record, apply $f_{bpa}$, add $A_{discretized}$ or
$P_{d-linear}$ respectively, and from hereon continue to work with the standard
IEEE precision.

Such a conversion \new{computationally} is not for free but reduces the memory
footprint.
First, few vertices with an uncompressed operator representation are
required simultaneously at \new{the same} time.
Those vertices which have not been used yet or where all adjacent cells have
been processed already can be held in compressed form.
Second, all stencils on fine-grid vertices are removed completely from the
persistent data structures as $\epsilon$ and $v$ here are simple.
Third, all Galerkin and BoxMG inter-grid transfer
operators are held persistently.
Yet, their hierarchical surplus often is very small, yields small $bpa$ and thus
is compressed aggressively.
The coarser the grid, the more bytes have to be invested in operator storage. 
This is not problematic as the number of coarser grid vertices is
small.

\begin{observation}
Our implementation is almost matrix-free in terms of storage.
\end{observation}

\noindent
\new{
Our fine grid statement is invalid if $\epsilon $ and $v$ on the fine grid are
homogenized from subgrid sampling. 
Homogenized stencils differ from rediscretization though we may expect
again that the difference is small.
We also note that our approach is orthogonal to \cite{Eckhardt:16:SPH} and
\cite{Bungartz:10:Precompiler} where we use a hierarchical transform to reduce
the memory footprint of unknowns or hierarchical surpluses, respectively.
In the present approach, the unknown per vertex is only one double compared to
$2\cdot 5^d+3^d$ (multiplicative) or $4\cdot 5^d+2\cdot 3^d$ (additive/BPX)
doubles required to store the stencils. 
A $bpa$-based compression of $u$ (and probably $b$) on top of the stencil
compression thus would only have a minor impact.
}

\section{Parallelization}
\label{sections:parallelisation}

Our parallelization considerations focus on shared memory and distributed memory
via tasking or MPI, respectively.
We rely on static
load balancing and \new{task} stealing only. 
That is, we concentrate on an academic discussion of potential concurrency in the 
linear algebra and postpone performance engineering.

\subsection{Shared memory}
\label{section:parallelisation:sharedmem}

The dynamic adaptivity plus vertical integration render standard loop-based
shared memory parallelization problematic, as we we do not assume that there are
larger regular grid region \new{\cite{Eckhardt:10:Blocking}.
There are no major loops well-suited for a \texttt{parallel
for} and tiling \cite{AlOnazi:17:TaskBasedAlgebraicMG}.}
We therefore derive a task-based parallelism formalized via
operation dependencies on the element-wise traversal.
These dependencies \new{can be resolved by a directed acyclic graph and then directly} can be
\new{passed to} any task-based library.

For all solver variants, first accesses to vertices may run in
parallel as long as the traversal preserves a top-down ordering for
$touchVertexFirstTime$. 
As soon as a set of vertices on level $\ell $ is loaded, all vertices on level
$\ell +1$ that share only the level $\ell $
vertices as parents can be handled in parallel.
On one level, $touchVertexFirstTime$ (and thus prolongation) is
embarrassingly parallel, with read-only access to the coarser level.

\begin{figure}
  \begin{center}
    \includegraphics[width=0.4\textwidth]{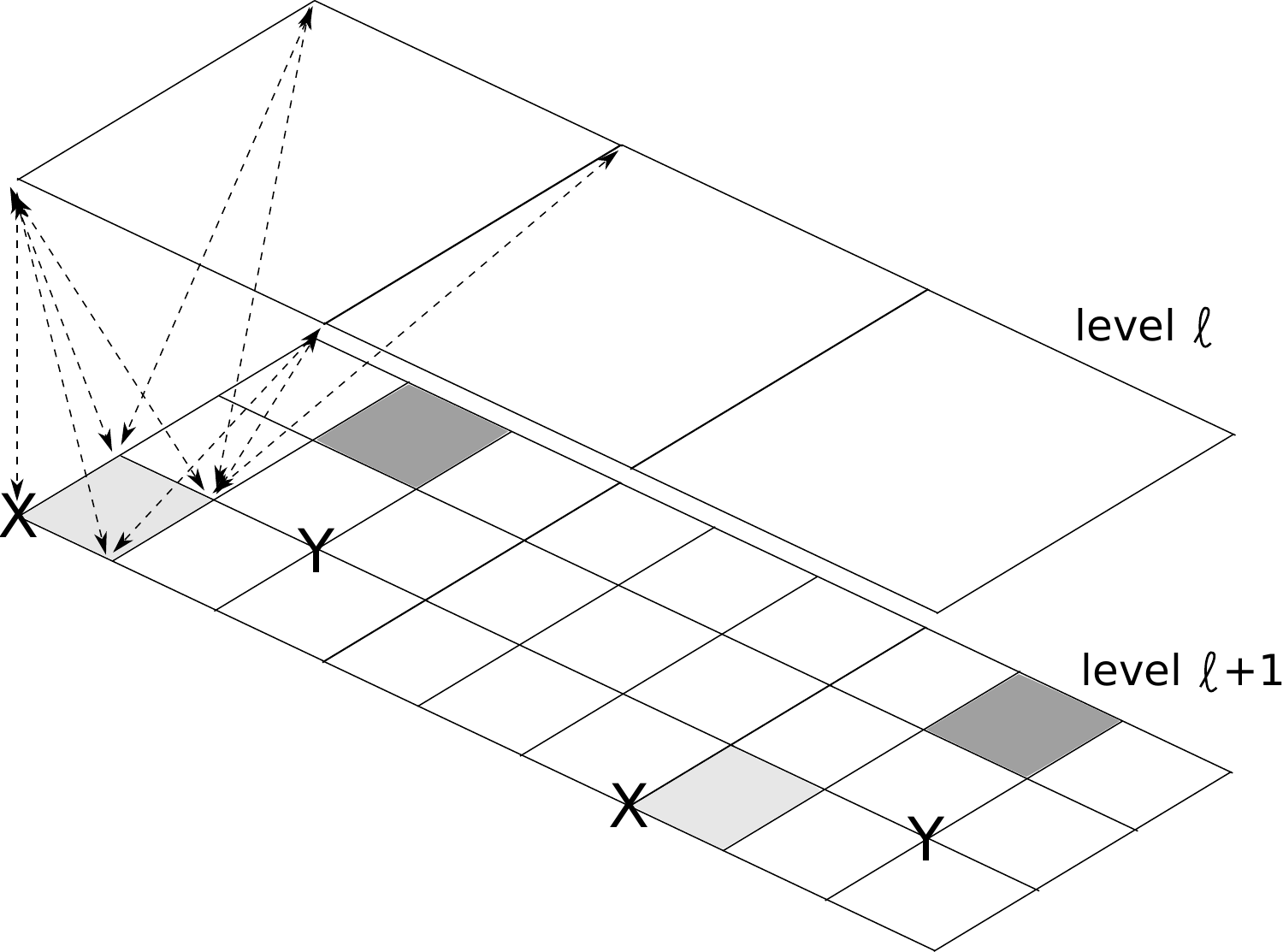}
    \includegraphics[width=0.4\textwidth]{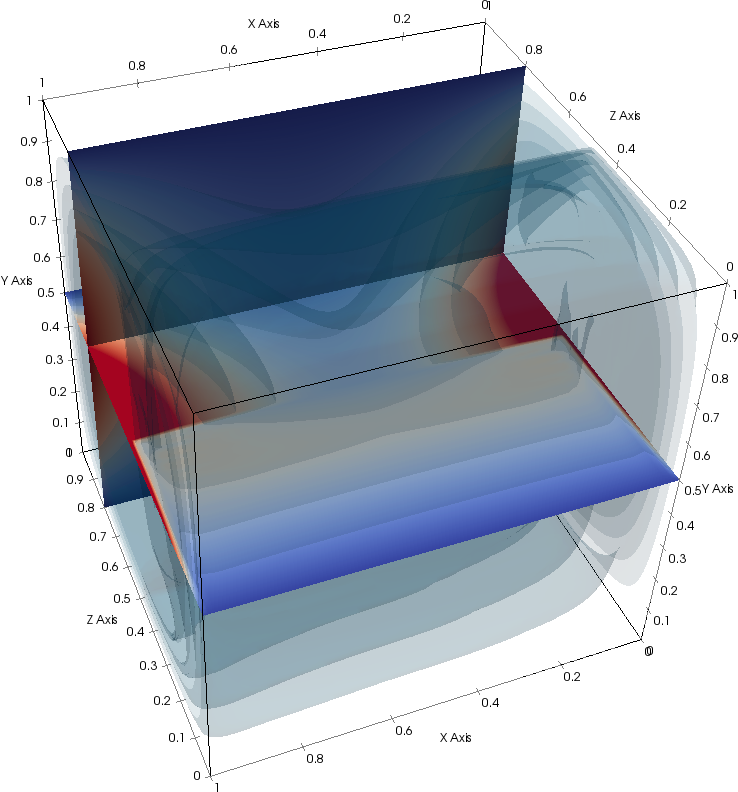}
  \end{center}
  \caption{
    Left: BoxMG concurrency within the spacetree. Cells of same grey shade
    can be processed in parallel.
    All vertices can be updated in parallel upon their very first usage.
    Those vertices carrying an \texttt{X} or \texttt{Y} marker, e.g., can be
    processed in parallel when the last operations per vertex per traversal are
    executed.
    Right: Solution to \new{the $\texttt{circle}$ benchmark} for $d=3$ with some contour
    faces. Diffusion dominates in the back half of the setup while $\epsilon
    = 10^{-4}$ in the front half allows the convection to yield an asymmetric
    solution.
    \label{figure:shared-memory:concurrency}
  }
\end{figure}

The element-wise residual computation exhibits a lower level of concurrency. 
For purely geometric multigrid solvers, no two cells may be updated concurrently
that share a vertex.
This induces a red-black type colouring \new{of the cells} with $2^d$ colours.
Galerkin and BoxMG solvers are more restrictive as they compute the residual
and modify the coarse-grid stencil and the inter-grid
transfer operators (Fig.~\ref{figure:shared-memory:concurrency}).
Here, two cells on level $\ell $ may be updated in parallel if their parent
cells on level $\ell -1$ do not share any common adjacent vertex.
If we read $\ell $ as a regular grid, this is a $(2k)^d$ colouring of the cells.
Such a multiscale dependency reduces the algorithm's concurrency severely.
However, we can work with $k^d$ colouring where possible and only use $(2k)^d$ colours 
while $recomputeGalerkin$ \new{in (\ref{equation:recomputeGalerkin})} holds.

$touchVertexLastTime$ updates the unknowns, restricts right-hand
sides and injects data. 
Updates are embarrassingly parallel.
As each vertex on level $\ell $ coincides spatially with at most one vertex on
level $\ell +1$, the injection of vertices on level $\ell +1$ is embarrassingly
parallel too.
Again, the multiscale traversal synchronizes the individual levels and ensures
that all vertices on level $\ell $ receive a $touchVertexLastTime$ before any of
their shared parents is handed over to this event.
The restriction imposes additional constraints.
Where the interpolation reads coarse-grid data only, the restriction reads
fine-grid data (modified by the update) and writes to the coarse grid.
No two vertices on level $\ell +1$ \new{that
share a parent} may thus be updated concurrently.
If we read level $\ell +1$ as a regular grid, this implies a $(2k+1)^d$
colouring of the vertices. 
We were not able to obtain reasonable speedups if we always \new{sticked} to
$(2k+1)^d$ colouring.
\new{Therefore,} we apply this colouring only for levels and
grid regions where the right-hand side needs to be re-determined.
Otherwise, we \new{process} vertices embarrassingly parallel.

The $descend$ events require a $2^d$ colouring on the coarser level $\ell $ if
$P$ and $R$ are re-computed.
Theoretically, all BoxMG patches can be computed in parallel.
\new{Yet, along} the interface
of two patches all adjacent patches compute the same entries due to the stencil
collapsing.
\new{Though this is} 
a race condition where multiple threads
determine the same data and write the same entries, a concurrent, redundant computation of some stencil entries
without a synchronization led 
invalid \new{stencil} entries \new{in our implementation}.
We thus fall back to $2^d$ (coarse \new{cell}) colouring where no two adjacent
$3^d$ patches determine BoxMG operators. 
For the block smoothers, we run \new{the patches}
parallel without locks.
\new{They} never modify any data on level $\ell +1$ (Fig.~\ref{fig:gridpoints}).

\subsection{Distributed memory}
\label{section:parallelisation:distributedmem}

The MPI parallelization of multigrid is an active area of research.
There are dozens of different strategies for any combination of solver variant,
machine and problem.
In line with the shared memory parallelization, we conduct a basic
data flow analysis for the multiplicative algorithms on
non-overlapping multiscale decompositions here.

 \begin{figure}
   \begin{center}
     \includegraphics[width=0.35\textwidth]{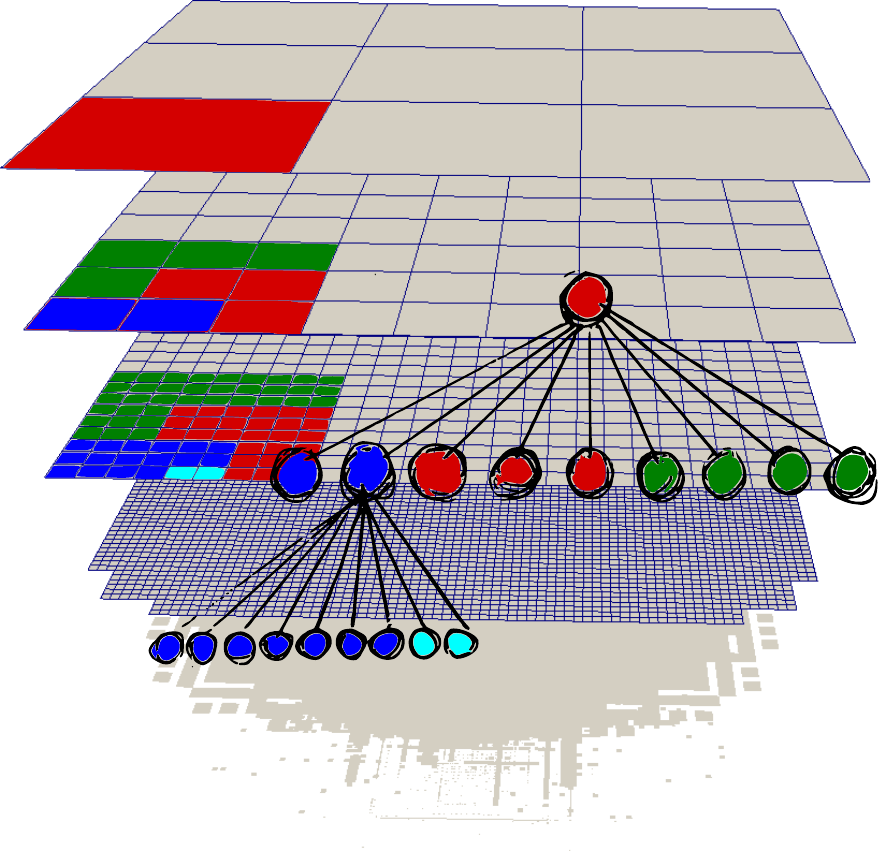}
     \hspace{0.8cm}
     \includegraphics[width=0.4\textwidth]{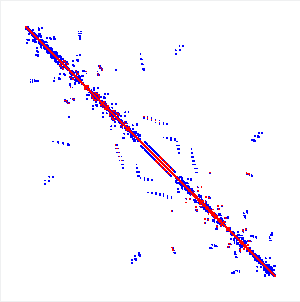}
   \end{center}
   \caption{
     Left: Domain decomposition along an SFC with four
     ranks highlighted.
     The decomposition induces a logical tree topology with masters and workers
     among the MPI ranks (from \cite{Weinzierl:17:PeanoSoftware}).
     \new{Right: Sparsity pattern of an explicitly assembled (PETSc) matrix for
     the $d=2$ setup on regular grids.
     As we use a space-filling curve to traverse the finest grid, the matrix 
     pattern is not dominated by diagonals.}
     \label{figure:mpi:spacetree-decomposition}
   }
 \end{figure}

%
A decomposition scheme fitting to our notion of element-wise spacetree
traversals is a classic non-overlapping domain decomposition.
Each cell in the fine grid is assigned to a unique rank, while vertices along
the domain boundary are replicated on each adjacent rank.
Two options exist how to handle coarser levels \cite{Weinzierl:17:PeanoSoftware}:
We can hold a refined spacetree cell on any rank that also holds one of
its children.
Such a bottom up construction of ownership implies that refined spacetree
nodes are replicated on multiple ranks---we are not non-overlapping
\new{in a multiscale sense}---and that vertices can be replicated on more than
$2^d$ ranks.
Notably, it implies that each rank holds the spacetree's root.
Alternatively, we can assign a refined spacetree cell uniquely to one of the set
of ranks that hold one its children
(Fig.~\ref{figure:mpi:spacetree-decomposition}).
Such a bottom-up construction yields a non-overlapping domain decomposition on
each and every level.
It implies that each cell has a unique owner, that vertices are adjacent to at
most $2^d$ ranks, that there is a logical tree topology induced on the MPI
ranks, and that most ranks holds only a fragment of the overall spacetree.
In our implementation, we follow the latter approach though all data flow
insights holds for both approaches.

%
%
We assume that all replicated data, i.e.~all values and stencils stored within
the vertices, are consistent after the initial grid construction.
All element-wise stencil evaluations can be
done without communication on all ranks holding a cell.
This yields a rank-local partial residual.
If cells are held redundantly, also the partial results are determined
redundantly.
As soon as the local residual on a domain boundary vertex is accumulated, we
send out this residual to all other ranks that hold a copy of this vertex and
continue the rank-local spacetree traversal.
A traversal realization that sticks to the single-touch policy for the
unknowns and hides data exchange behind computation postpones the update
of a vertex's unknown. i.e.~no update of unknowns is conducted right away.
Instead, we retain the residual in-between two grid traversals and add a prelude
to the vertex's load process of the subsequent traversal.
Data is sent out 
when a vertex has been processed for the last time, but it is not merged into
data structures prior to the subsequent grid traversal.
The prelude receives all residual contributions from all other adjacent ranks,
adds them to the local residual and performs the unknown update.
If we had redundant cells, the accumulation of the residual has to anticipate
that some residual fragments might be determined multiple times.

\new{The} postponing of unknown updates decreases the
speed of the smoother by one grid sweep in total.
Amortized over all sweeps, our parallelization does not alter the
convergence speed.
In return, the data exchange can be realized in a non-blocking manner in the
background.
However, we experience a slight reduction of the speed along grid resolution
changes.
Our FAS-based handling of coarser grid regions requires the injected fine-grid solution.
If the smoother's fine grid updates are postponed to the subsequent iteration,
no valid injected data from the current iteration is available there yet.
In the parallel code, updates are injected into coarser levels one sweep later
than in the serial case.
This has an impact on smoothers working on adaptive grids.
At hands of the domain decomposition description of such smoothers, we see
that the fine-to-coarse domain coupling is delayed while the
coarse-to-fine coupling through interpolation along the hanging vertices remains
tight.

While a relaxed inter-level coupling is acceptable, we face a more severe data
consistency challenge throughout the reduction \new{of the multiplicative
algorithm}.
The hierarchical residuals---partially computed along the domain
boundaries---are restricted per rank.
The restricted values then are exchanged to determine the right-hand side for
the coarse grid correction.
This postponed data flow scheme mirrors the exchange of residuals.
Our serial/shared memory algorithms propose to fuse the coarse grid smoothing 
with the restriction.
This relies on the facts that (i) the right-hand side is
correctly restricted when a vertex is used for the last time throughout a
traversal, (ii) the matvec accumulation of the correction is independent of
the right-hand side, and (iii) the fine grid computation of $\hat r$ has no
influence on the injected coarse grid representation $u_{\ell -1} = Iu_{\ell
}$.
Obviously, constraint (i) and (ii) are shifted to the begin of the subsequent
grid sweep in our scheme.
\new{Through this,}
constraint (iii) is harmed.

We thus run $\mu _{pre}$ pre-smoothing steps per level, and then add an
additional grid traversal to complete the smoothing, inject the solution to the
next coarser grid and restrict the right-hand side locally.
The partially restricted right-hand sides are sent out at the end of the grid
sweep and thus are available on remote ranks once the first smoothing step on
the coarser level starts.
This break-up of smoothing and restriction into two separate grid
traversals reduces the solver efficiency---additional $\ell _{max}-\ell _{min}$ grid sweeps
are required per V-cycle---but it allows us to keep the right-hand sides
consistent.
No data consistency problems arise during the steps down within the V-cycle as
we realize our multigrid solver within depth-first tree traversals.

\begin{observation}
\label{observation:mpi}
We have to invest one additional grid sweep per multigrid restriction \new{in
the multiplicative case}.
Otherwise, the data consistency can not be preserved with data exchange in the
background of the computation.
\end{observation}

%
%
%

Our algorithm exchanges two doubles per refined vertex ($r$
and $b$) and only the residual for fine-grid vertices.
For the dynamic adaptivity criterion, additional
quantities are exchanged.
While the exchanged attribute cardinality is low, it is important to recognize
that the two residuals on boundary vertices now have to be stored persistently.
Without MPI, we are able to discard them after each traversal.
This increases the memory footprint along domain boundaries.

The Galerkin coarse grid stencils are computed additively over cells. 
Consequently, we may send out the partial stencils in the synchronization
traversal and receive and accumulate them in the preamble of the follow-up
smoothing sweep.
BoxMG determines $P$ and $R$ through local patch computations. 
While we use a non-overlapping multiscale domain decomposition, we weaken this
concept and replicate those cells required to compute all patch operations on
one rank. 
Consequently, the $descend$ event requires no special attention.
Along the boundary, all rank-local $descends$ yield solely partial
inter-grid transfer operators. 
As BoxMG computes entries along a patch boundary
redundantly---this statement also holds for injection and trivial
aggregation---all entries affecting local vertices are always available.

Additive multigrid is sketched in \cite{Mehl:06:MG} and can be realized
following the present data flow ideas as long as no FAS is employed.
If we use FAS, Observation \ref{observation:mpi} immediately implies that the
scheme runs into inconsistent data.
A solution to this approach is the pipelining approach from
\cite{Reps:15:Helmholtz} that also fixes the weakened domain coupling of meshes
of different resolutions in the adaptive case.
BPX follows these lines.
It is an open question whether the pipelining concept applied to the
multiplicative setup yields a realization that is superior to the present
approach with an additional grid sweep per restriction step.

We conclude our distributed memory discussion with the observation that our code
family keeps redundantly held vertex data consistent---if required by additional
spacetree sweeps.
As a consequence, stencil compression can be applied in parallel on boundary
vertices, too, where it yields the same compression factors for redundantly
held grid entities.
It is a straightforward decision thus to exchange the compressed byte streams
instead of the real stencils.

\section{Results}
\label{sections:results}


\begin{table}
  \tbl{Overview of the benchmark parameter sets for
    (\ref{equation:pde}).
    \label{table:parameters}
  }{
  {\footnotesize 
  \begin{tabular}{l|ccc}
    Identifier
     & $\texttt{sin}$
     & $\texttt{jump}$
     & $\texttt{checkerboard}$
     \\
   \hline
   \\[0.002cm]
   $f$
     & $2\pi ^2 \Pi _i sin (\pi x_i)$
     & 1
     & 1 \\[0.1cm]
   $u|_{\partial \Omega}$ 
     & 0 & 0 & 0 \\[0.4cm]
   $v$   
     & 0 & 0 & 0 \\[0.6cm]
  $\epsilon$ 
    & 1 
    & $\left\{\begin{array}{rl}
        1    & \mbox{if}\ x_1<0.5 \\
        0.1  & \mbox{otherwise} \\
        \end{array}\right. $
    & $\epsilon _i = \left\{   
       \begin{array}{rl}
       1 \ \mbox{if}\ x_i<0.5, \\
     0.1 \ \mbox{otherwise}
      \end{array} \right. $        
      \\[0.4cm]
   \hline
   \hline
   Identifier
   & && $\texttt{circle}$ \\
   \hline
   \\[0.002cm]
   $f$ && 
     & 0 \\[0.1cm]
   $u|_{\partial \Omega}$ && 
     & $\left\{\begin{array}{ll}
       1-4(x_2-0.5)^2 & \mbox{if}\ x_1=0 \vee x_1=1\\
       0 & \mbox{otherwise}
       \end{array}\right.$ \\[0.4cm]
   $v$ &&  
    & $\left( 
      \begin{array}{l}
       sin(\pi x_2 -0.5)cos(\pi x_1 -0.5) \\
       -cos(\pi x_2 -0.5)sin(\pi x_1 - 0.5) \\
       0 \\
      \end{array}
    \right)$   \\[0.6cm]
  $\epsilon$ &&  
    & $ \left\{ \begin{array}{lcl}
     \{10^{-1}, 10^{-2}, 10^{-3}, 10^{-4}\}
    & \mbox{if} & x_3 \leq 0.5 \ \mbox{or $d=2$}  \\
    1 & & \mbox{otherwise}
    \end{array} \right\} $
  \end{tabular}
  }} 
\end{table}


\new{Our benchmarks study
 (\ref{equation:pde}) with the parameter sets from Table
 \ref{table:parameters}.
}
\new{To realize}
dynamic adaptivity, we evaluate the mean value of
the $3^d-1$ surrounding vertices per non-hanging vertex and compute the absolute value of the difference of this mean value to the actual vertex
value.
The feature-based criterion assumes that refinement pays
off where the problem changes rapidly, i.e.~\new{where} this difference is
significant.
A region around a vertex is a refinement candidate if the vertex is unrefined
and if the residual in the particular vertex falls below $10^{-2}$.
Per grid sweep, we refine the 10 percent of the candidates with the biggest
absolute mean value differences \new{as long as} the minimal mesh width permits.
To avoid a global sorting step per iteration, our code uses binning with ten
bins, with each bin representing a certain range of the refinement criterion's
differences.
This range is adapted after each traversal as we keep statistics per bin on how
many vertices fit into it.
All vertices fitting into the bin representing the mean value difference are
refined\new{---the 10 percent goal is approximated}.

%
%
%

\new{
All experiments were conducted  on a cluster with Intel E5-2650V4 (Broadwell)
nodes with 24 cores per node.
They run at 2.4 GHz.
Furthermore, we reran the multicore experiments on an Intel Knights Landing
chip (Xeon Phi 7250) at 1.4 GHz.
For the shared memory parallelization, we rely on Intel's Threading Building
Blocks (TBB).
For the distributed memory parallelization, we use Intel MPI.
Intel's C++ compiler 2017 version 2 translates all codes.
The realization is based on the spacetree PDE
framework Peano \cite{Software:Peano}.
As we stick to a low order discretization on dynamically adaptive meshes and do
not exploit any mesh regularity, we only exploit the compiler's vectorization of
the BoxMG matrix-vector products.
Further vectorization along the lines of \cite{Reps:15:Helmholtz} for
multi-parameter runs or \cite{Eckhardt:10:Blocking} for grids with some (patch)
regularity is beyond scope.
}

\begin{figure}
 \begin{center}
  \includegraphics[width=0.9\textwidth]{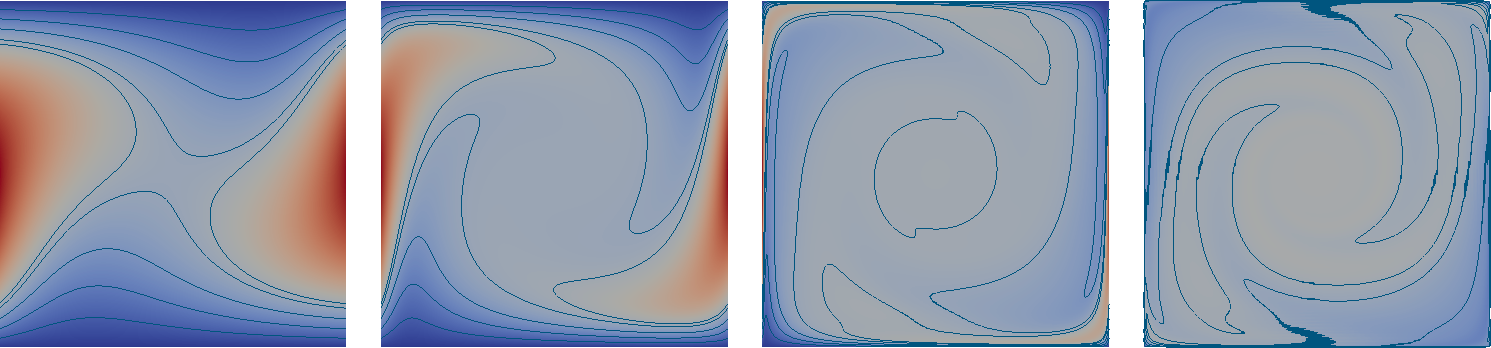}
 \end{center}
 \caption{
   Solution of \new{the $\texttt{circle}$ problem (Table
   \ref{table:parameters}) } for $d=2, \epsilon \in
   \{10^{-1},10^{-2},10^{-4},10^{-8}\}$ (left to right) with isolines
   at $u \in \{0.1,0.2,0.3,0.4,0.44\}$.
   \label{figure:results:circle-2d}
}
\end{figure}

\subsection{Diffusion with constant coefficients}
\label{section:results:sin} 

%
%
An analytical solution is known for the \new{$\texttt{sin} $ benchmark}.
Geometric rediscretization 
yields the Galerkin coarse-grid operator if $d$-linear prolongation and
restriction are chosen \new{that} 
BoxMG yields.
The setup thus acts as validation scenario.

\begin{table}
  \tbl{
    Number of grid sweeps for additive multigrid that are required to reduce the
    residual for the \new{$\texttt{sin} $ benchmark (Table \ref{table:parameters})} by a factor of $10^{-8}$. 
    The sections show $d=2$ (top) and $d=3$ (bottom) with $\omega = 0.8$. 
    All experiments are computed on regular Cartesian grids spanned by the
    spacetree, i.e.~we build up the whole grid and then start the solve.
    The tuples denote cycle counts with exponential damping (left) and undamped
    coarse grid relaxation (right).
    Jac denotes a Jacobi smoother, BJ is a block Jacobi with the number of
    Gau\ss -Seidel sweeps per block in brackets. 
    The /e postfix implies that we use an exact coarse grid solve.
%
%
    \label{table:experiments:sin:additive}
  }
  {
  \footnotesize
  \begin{tabular}{l|ccccc|ccccc}
   h & Jac & BJ(1) & BJ(2) & BJ(4) & BJ(8) & Jac/e & BJ(1)/e & BJ(2)/e & BJ(4)/e & BJ(8)/e \\
\hline
$3^{-2}$ & 34/26 & 19/19 & 16/16 & 13/13 & 12/12 & 34/26 & 16/16 & 13/13 & 12/12 & 12/12 \\
$3^{-3}$ & 48/41 & 24/22 & 20/18 & 17/18 & 17/18 & 48/41 & 22/21 & 19/18 & 17/18 & 17/18 \\
$3^{-4}$ & 63/44 & 29/24 & 24/21 & 22/21 & 22/21 & 63/44 & 27/23 & 24/21 & 22/21 & 22/21 \\
$3^{-5}$ & 82/47 & 34/24 & 30/23 & 28/23 & 28/23 & 82/47 & 33/24 & 30/23 & 28/23 & 28/23 \\
$3^{-6}$ & 98/45 & 41/25 & 37/26 & 35/26 & 35/26 & 98/45 & 40/25 & 37/25 & 35/26 & 35/26 \\
  
   \hline 
$3^{-2}$ & 21/19 & 19/19 & 17/17 & 17/17 & 17/17 & 21/19 & 17/17 & 17/17 & 17/17 & 17/17 \\
$3^{-3}$ & 42/39 & 26/25 & 23/22 & 22/22 & 22/22 & 42/39 & 25/24 & 23/22 & 22/22 & 22/22 \\
$3^{-4}$ & 51/39 & 32/27 & 29/25 & 29/25 & 28/25 & 51/39 & 31/27 & 29/25 & 29/25 & 28/25 \\
  
  \end{tabular}
  }
\end{table}

All required records per vertex are enlisted in Table \ref{table:unknowns}.
For the pure geometric solver, we only hold three or four, respectively,
doubles per grid vertex. 
Exhaustive search over potential relaxation parameters for the additive
solver results in $\omega \approx 0.8$, yielding reasonable
convergence rates though no relaxation or slight overrelaxation is even faster
(Table \ref{table:experiments:sin:additive}).
As soon as we switch to a dynamically adaptive solver, we find that
$\omega \geq 1.0$ becomes unstable. 
Thus, we stick to $\omega = 0.8$ from hereon.
Better convergence rates might be obtained for alternating
relaxation parameters, where we use a different parameter for each sweep.
We refer to \cite{Huckle:08:FourierToeplitzAnalysis} for some symbol analysis
or \cite{Reps:15:Helmholtz} for a Helmholtz \new{example}.
If we damp $\omega$ exponentially---\new{use $\omega $ on the finest grid,
$\omega ^2 $ on the first coarse grid, and so forth}---we harm the
speed.
However, we obtain a stable scheme, while otherwise the additive solver tends
to overshoot 
\cite{Bastian:98:AdditiveVsMultiplicativeMG,Reps:15:Helmholtz}.
The speed gap between exponential $\omega $ damping and 
smoothing with uniform relaxation factor narrows if we use a block
smoother, but it does not close completely.
%
%
Block smoothers can double the convergence speed, but more than four
Gau\ss -Seidel block sweeps rarely pay off.
An exact coarse grid solve does not pay off for the additive solvers.

\begin{figure}
 \begin{center}
  \includegraphics[width=0.4\textwidth]{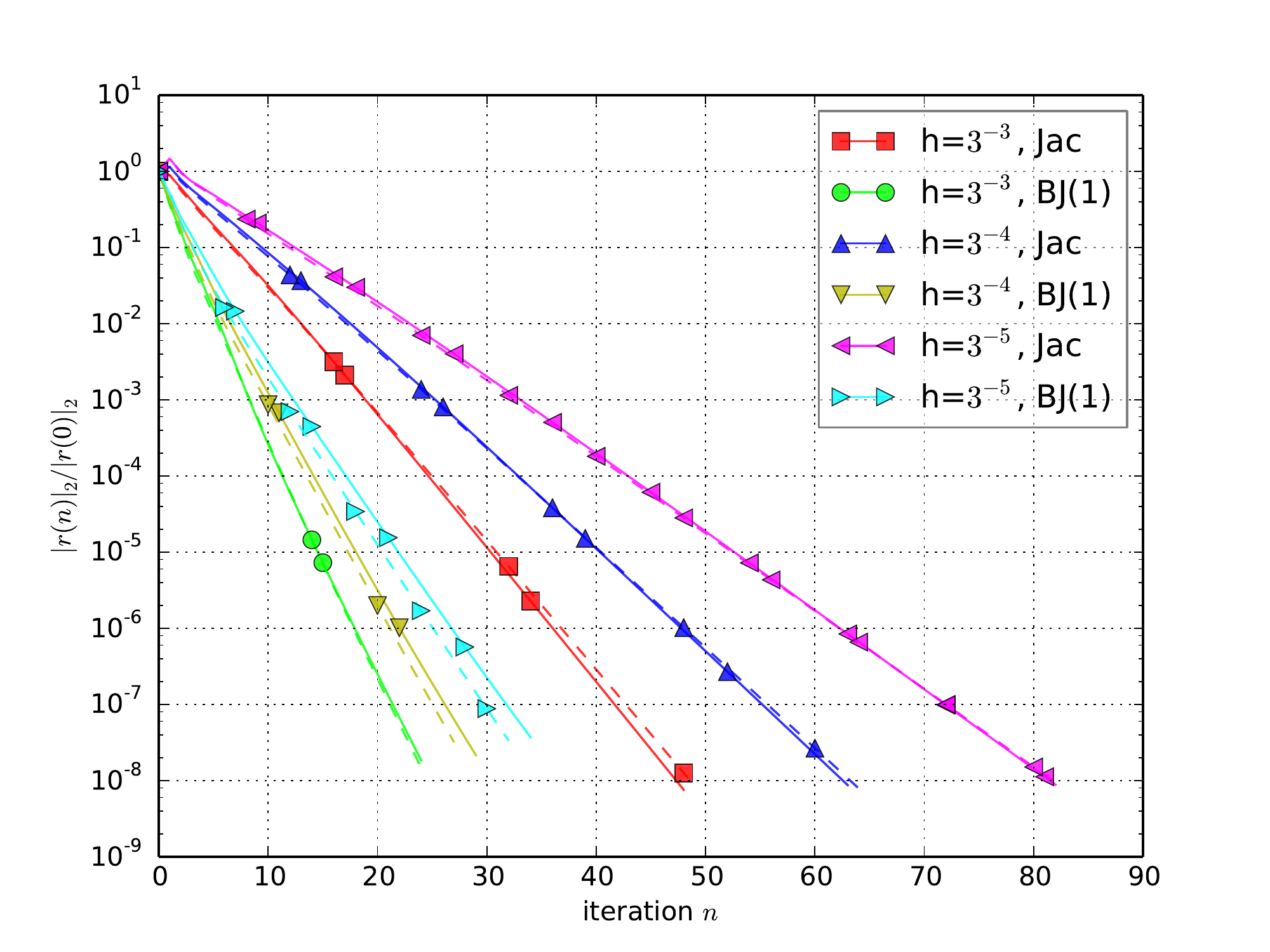}
  \includegraphics[width=0.4\textwidth]{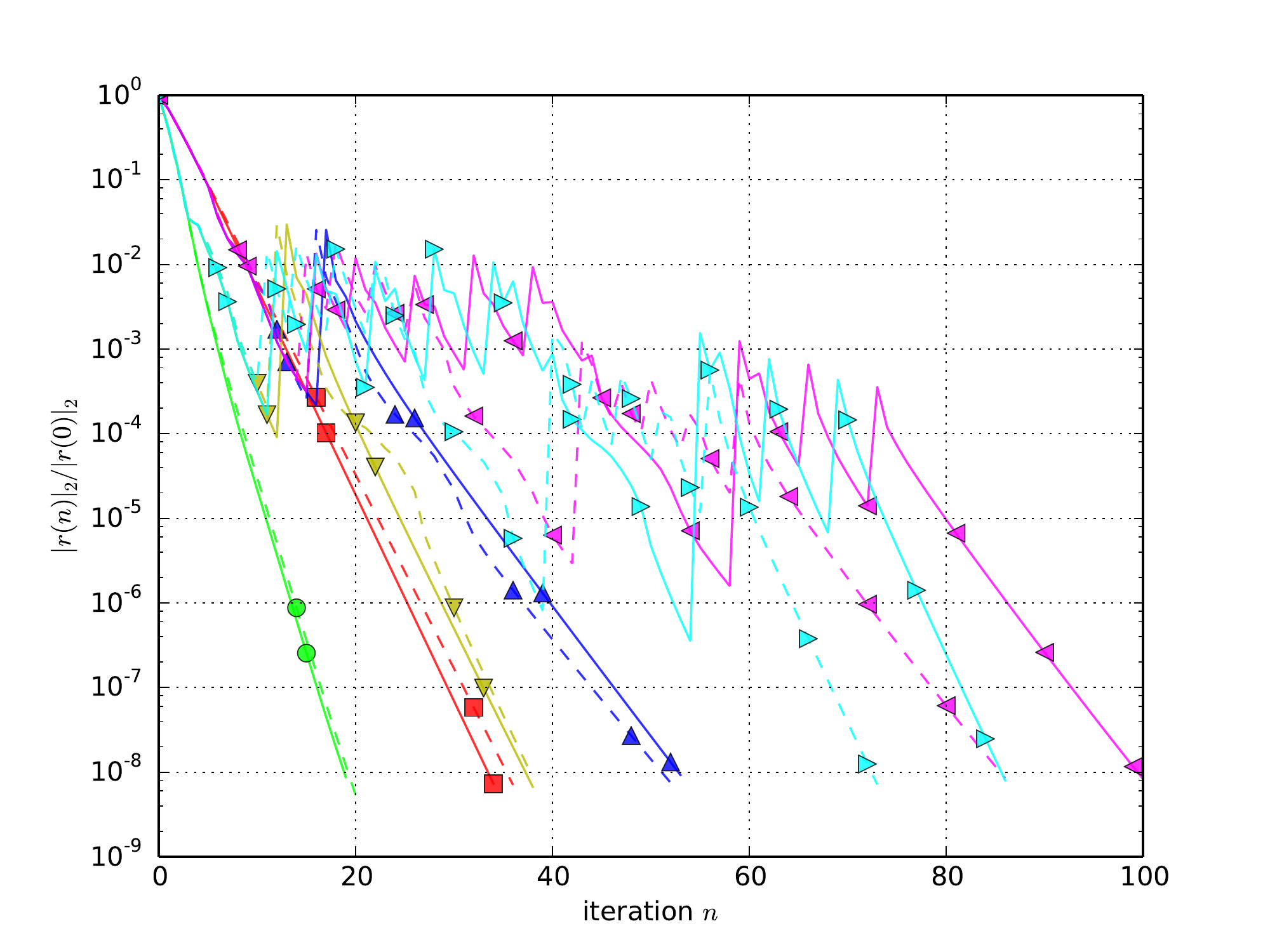}
 \end{center}
 \caption{
   Convergence behaviour of additive multigrid with exponential damping for
   the \new{$\texttt{sin} $ benchmark} on regular grids (left) and dynamically adaptive grids
   (right) with $\omega = 0.8$ and $d=2$.
   Dynamical implies that we start F-cycle like with a coarse grid and leave it
   to the refinement criterion to yield the final mesh.
   Solid lines show results
   for Algorithm \ref{algorithm:geometric::additive} compared to a synchronized variant
   from \cite{Reps:15:Helmholtz} (dotted lines).
   \label{figure:experiments:sin:synchronisation}
 }
\end{figure}

In Algorithm \ref{algorithm:geometric::additive}, fine-grid updates are
immediately injected to the coarser grids. 
In turn, coarse-grid computations might work with outdated coarse solutions,
which change during the element-wise assembly.
Our experiments show that this inconsistency does not make a difference for
regular grids (Fig.~\ref{figure:experiments:sin:synchronisation}).
It, however, slightly deteriorates the convergence for adaptive grids.
Here, we start from $h=3^{-1}$ and make the adaptivity criterion add further
vertices.
The effect is studied \new{and a single-touch solution is proposed in
\cite{Reps:15:Helmholtz}, and we conclude for our experiments that their
multilevel synchronization with pipelining should be used}.

\begin{figure}
 \begin{center}
  \includegraphics[width=0.4\textwidth]{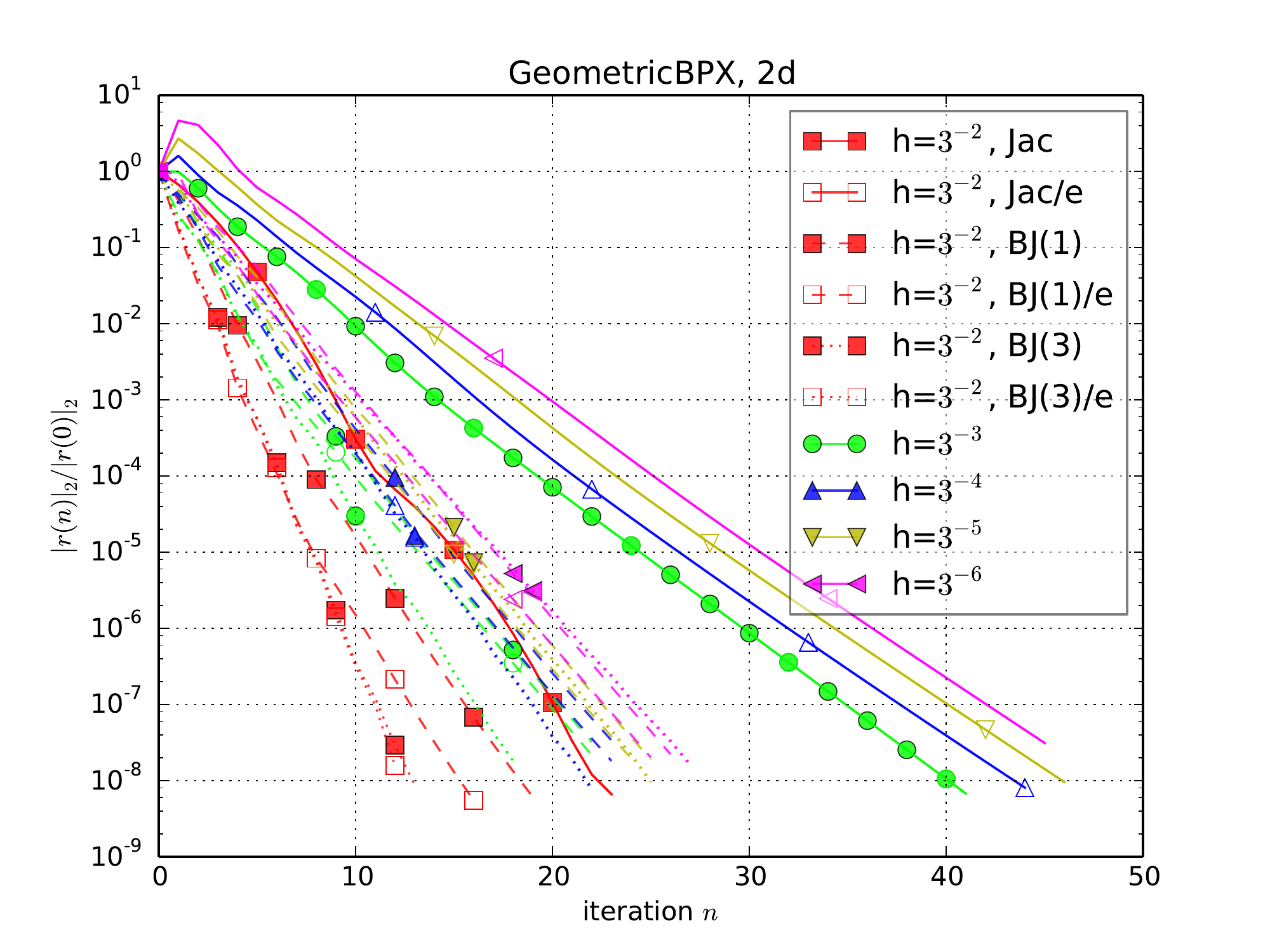}
  \includegraphics[width=0.4\textwidth]{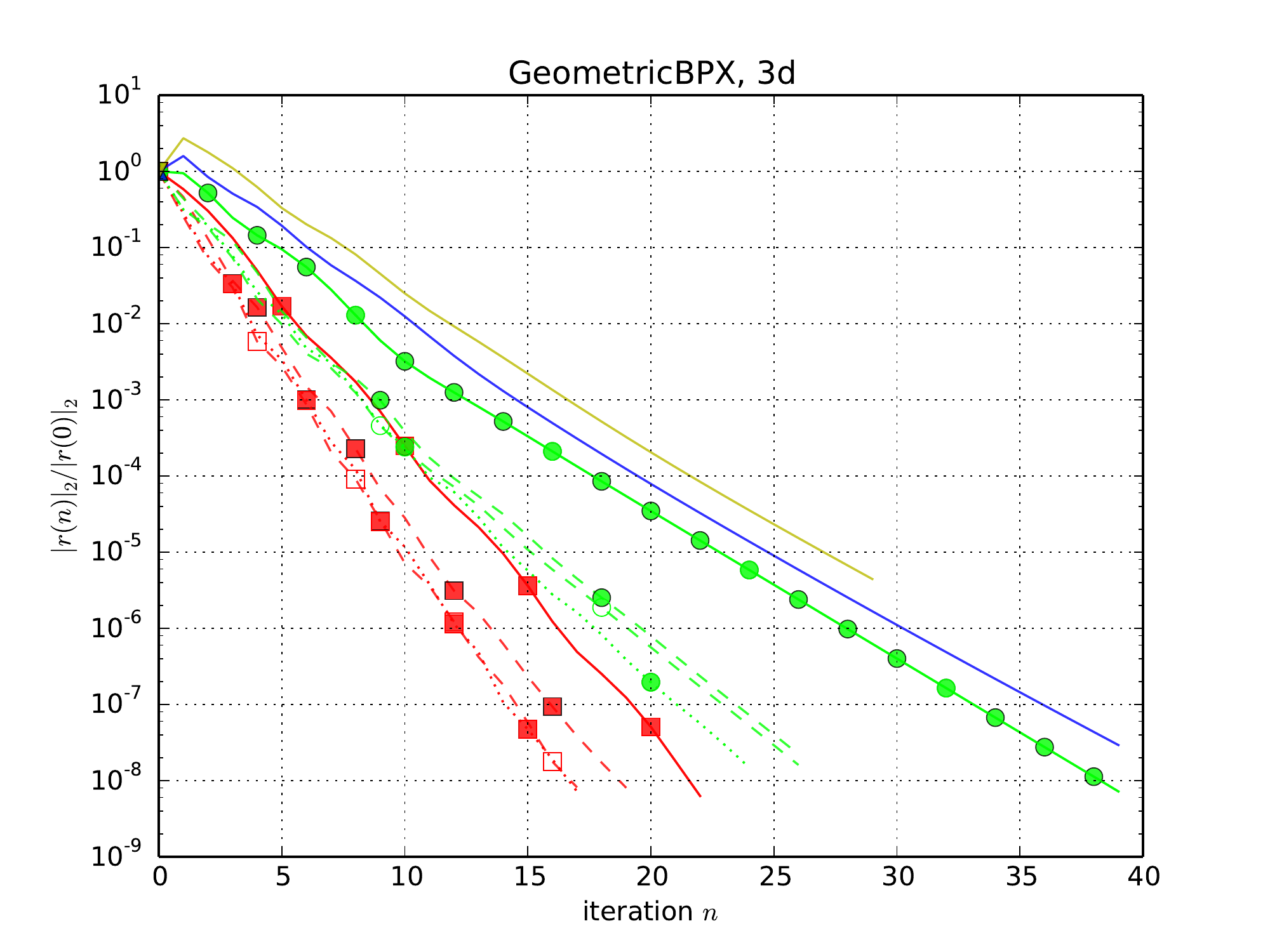}
  \includegraphics[width=0.4\textwidth]{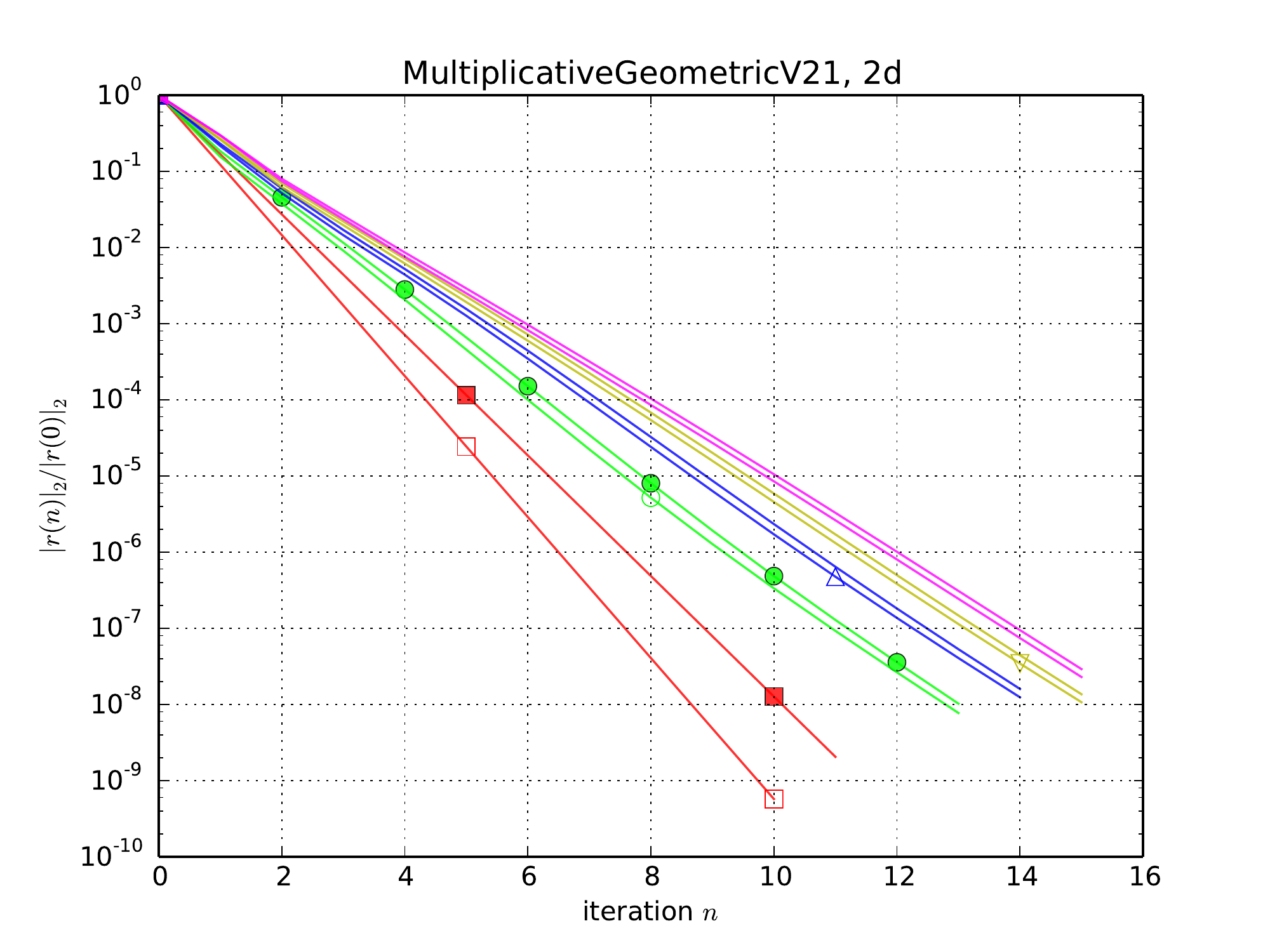}
  \includegraphics[width=0.4\textwidth]{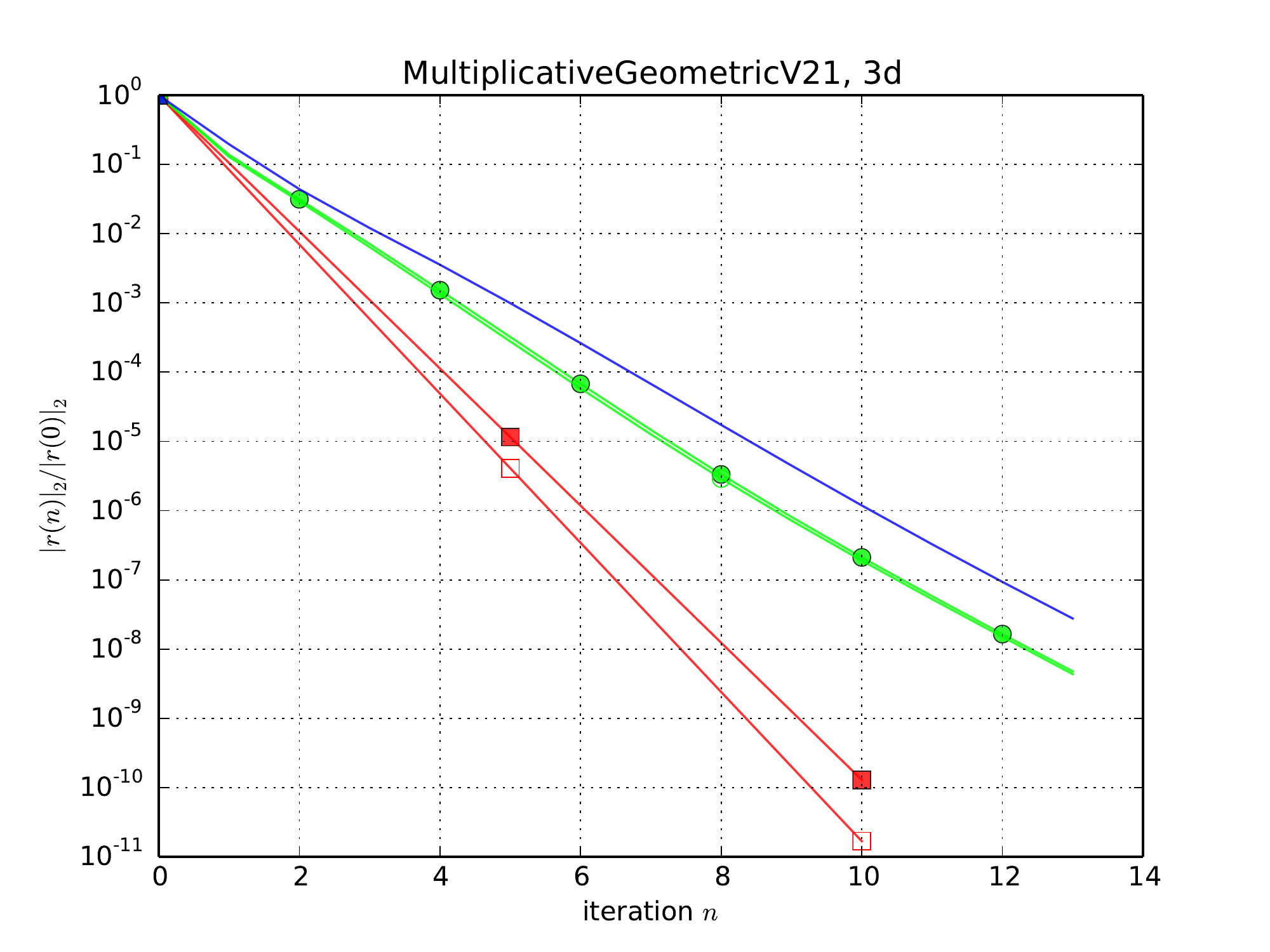}
  \includegraphics[width=0.4\textwidth]{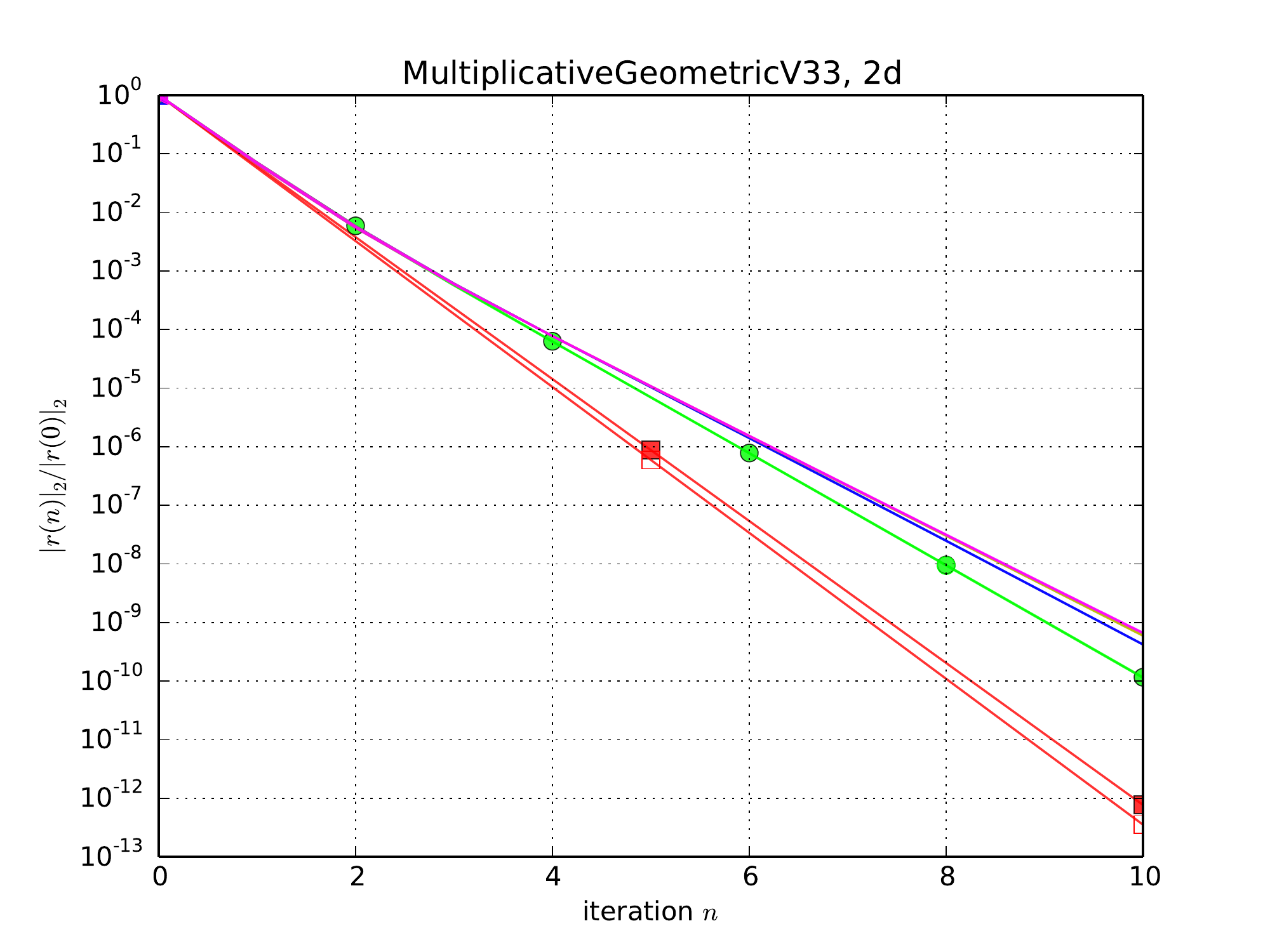}
  \includegraphics[width=0.4\textwidth]{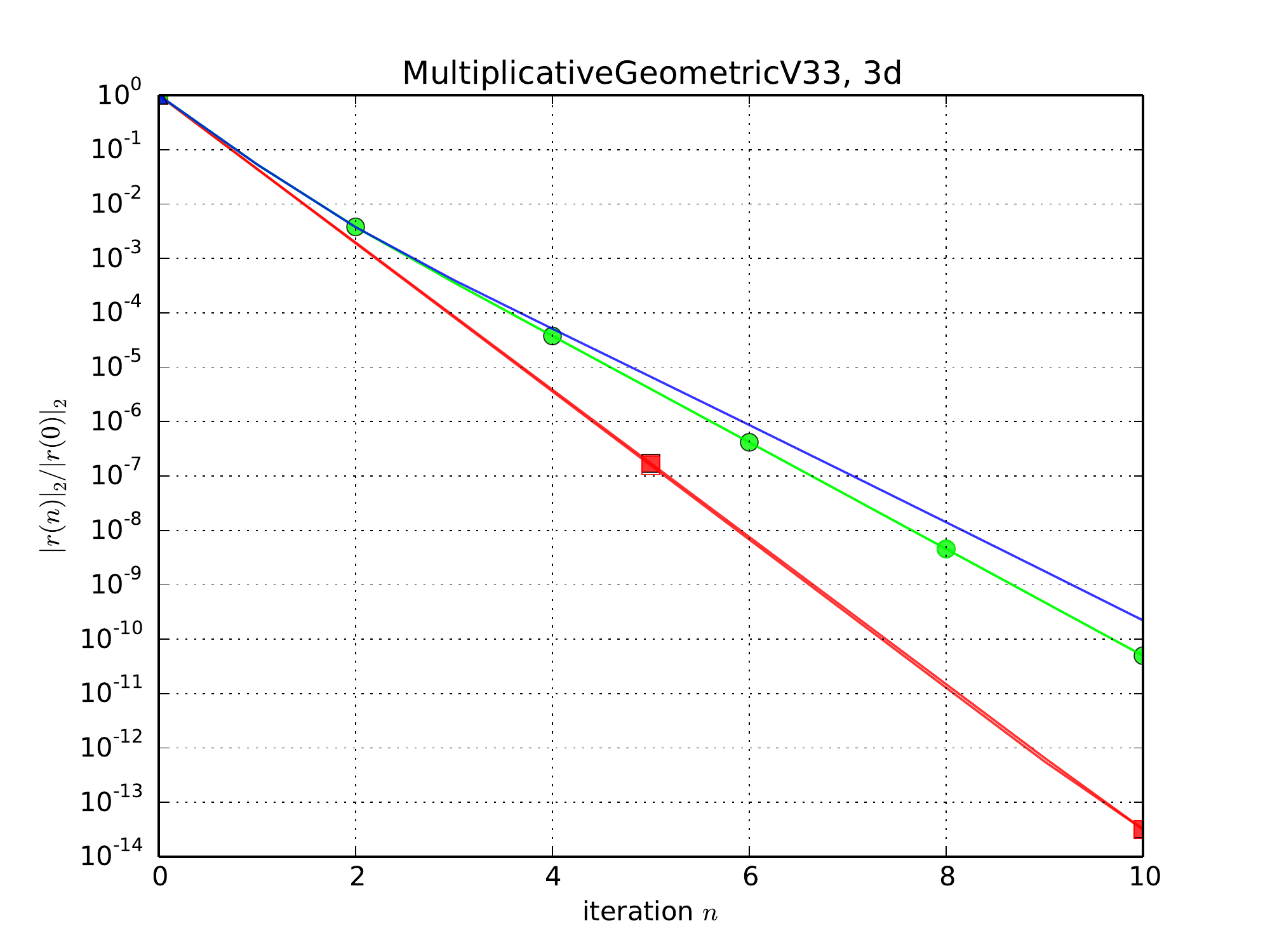}
 \end{center}
 \caption{
   Convergence of various solvers for the 
   \new{$\texttt{sin} $ benchmark} on regular grids with $\omega = 0.8$.
   Solid lines with filled symbols use a Jacobi smoother, dashed
   lines use block Jacobi with one Gau\ss -Seidel sweep, and
   dotted lines apply three sweeps.
   Empty symbols furthermore are for solving the coarse-grid problem exactly
   while solid symbols are for applying only $\mu _{pre}+\mu _{post} $ sweeps on
   the coarsest grid.
   \label{figure:experiments:sin:solvers}
 }
\end{figure}

All statements on additive solvers also hold for BPX, besides the fact that an
exact coarse grid solve here has no major positive impact
(Fig.~\ref{figure:experiments:sin:solvers}).
In general, BPX from Algorithm \ref{algorithm:geometric::bpx}
outperforms the additive solver.
Block smoothing pays off.
The impact of block smoothing on multiplicative multigrid is even more
significant.
For the latter, exact coarse-grid solves are advantageous and we need 
\new{10--15} iterations in total.

\begin{figure}
 \begin{center}
  %
  %
  %
  \includegraphics[width=0.4\textwidth]{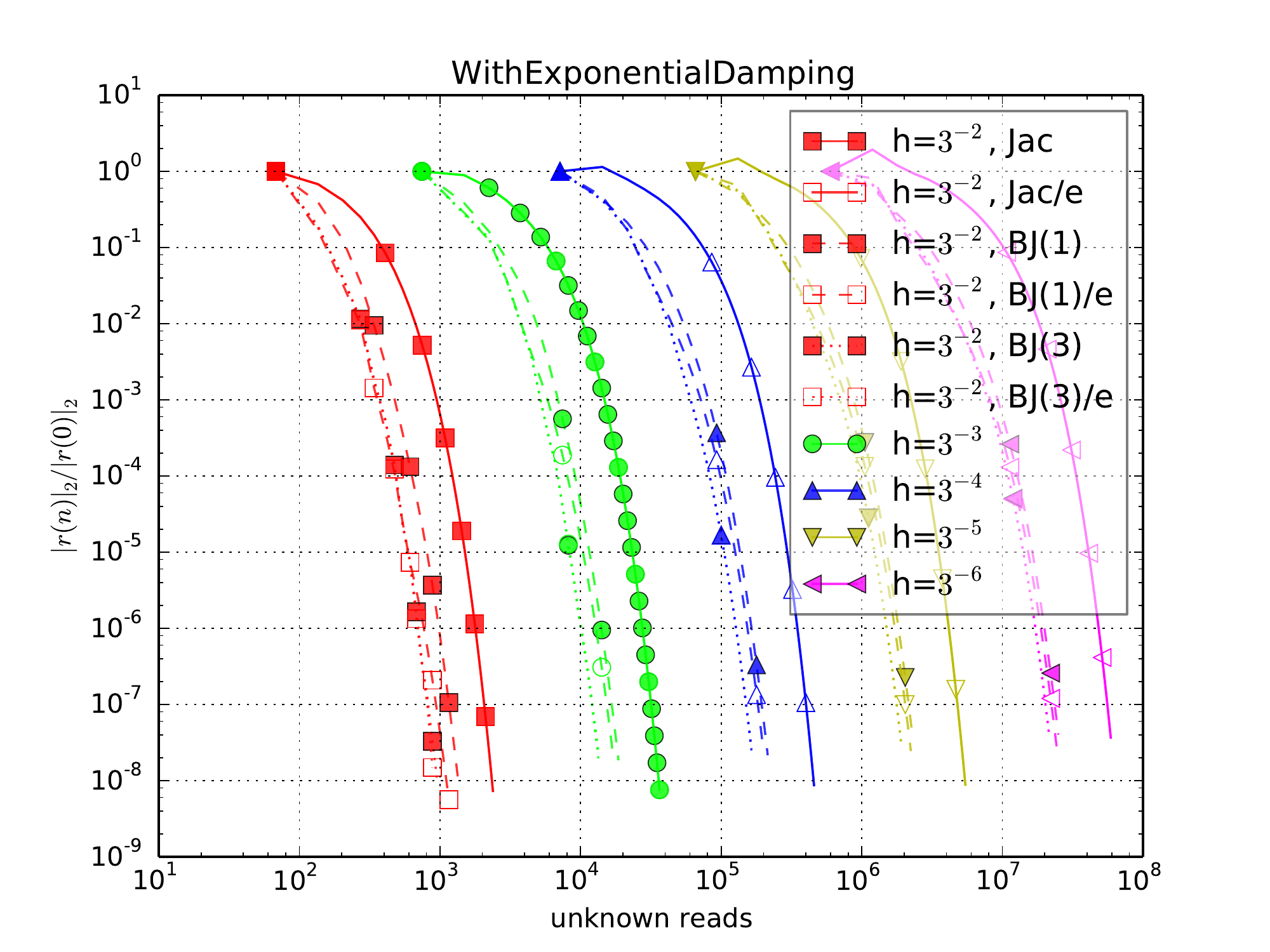}
  \includegraphics[width=0.4\textwidth]{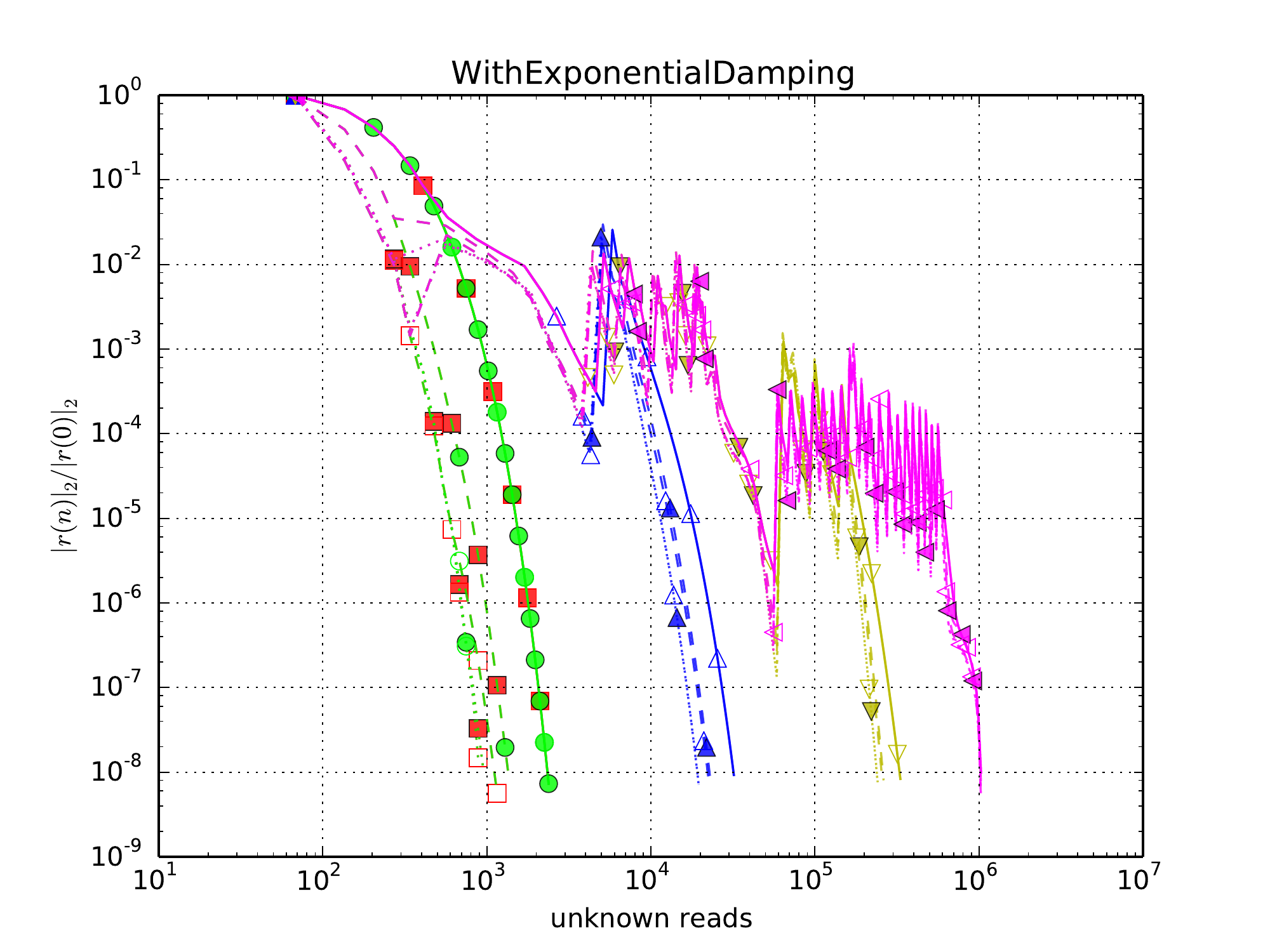}
  \includegraphics[width=0.4\textwidth]{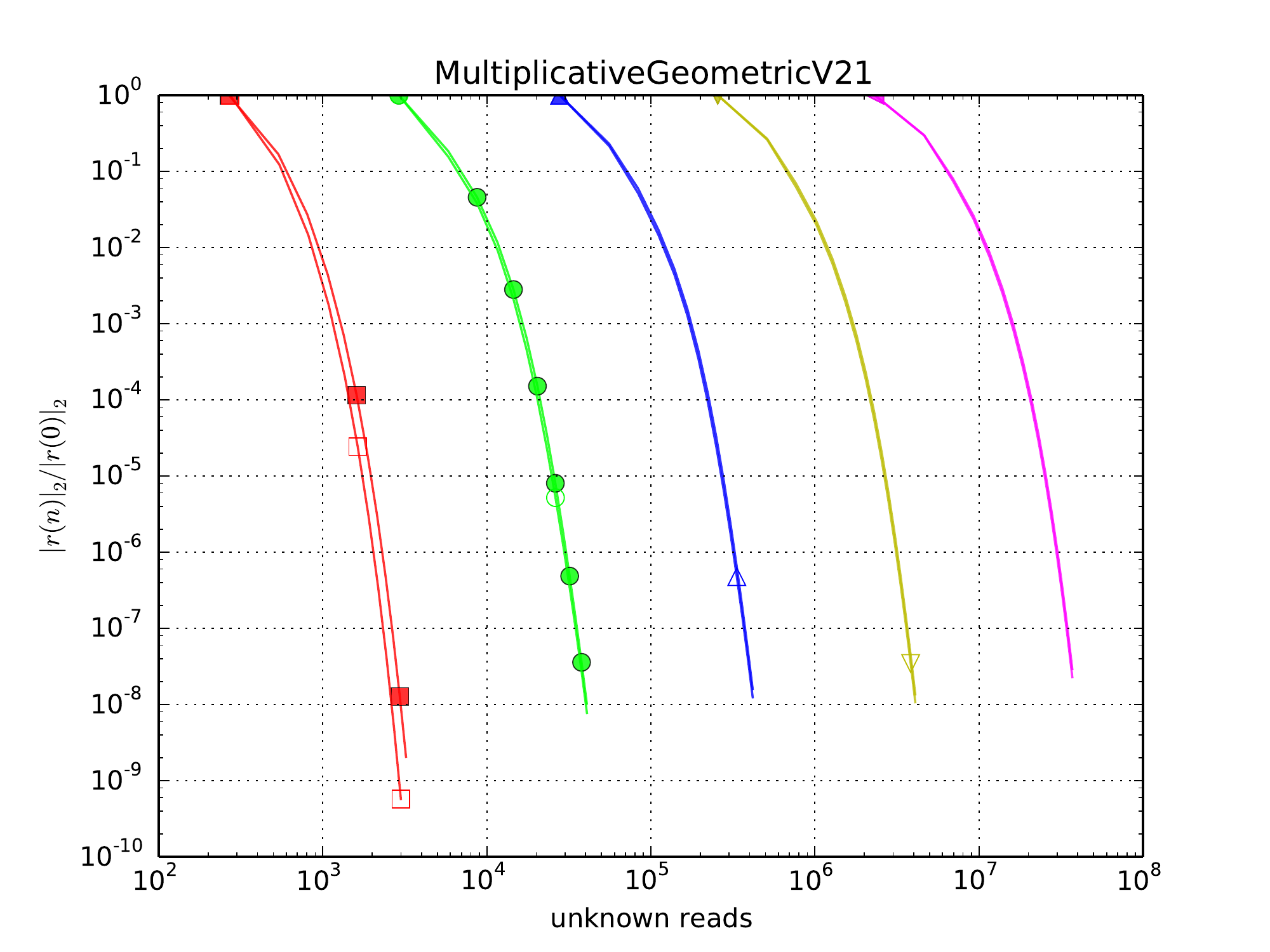}
  \includegraphics[width=0.4\textwidth]{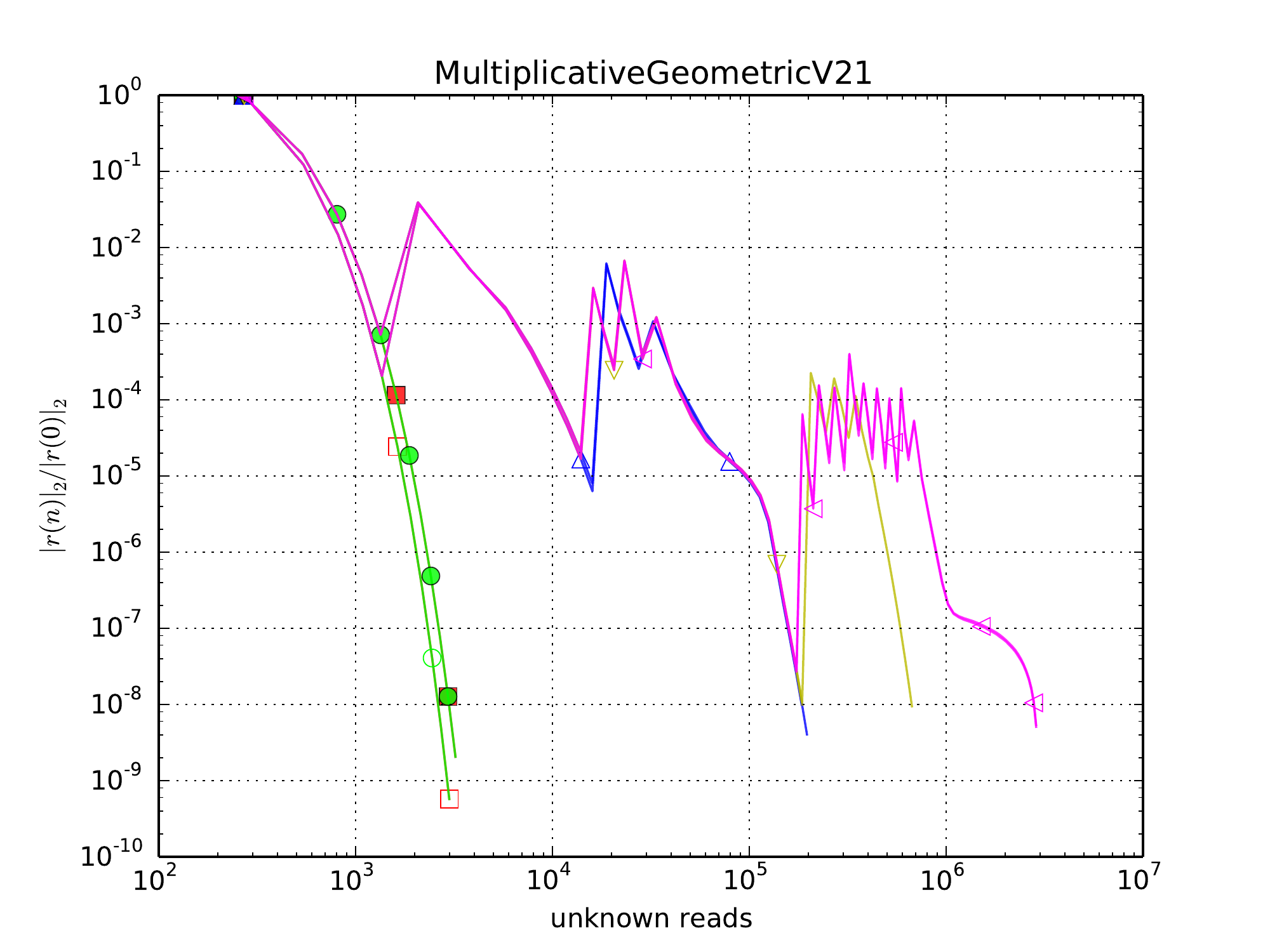}
 \end{center}
 \caption{
   Cost in terms of unknown reads for the \new{$\texttt{sin} $ benchmark} with $\omega = 0.8$ and $d=2$.
   Regular grids (left column) are compared to grids that unfold
   dynamically through the adaptivity criterion (right).
   The additive solver (top) is slower by a factor of one up to two compared to
   BPX (not shown) which is in the same order of cost as the multiplicative
   variant (bottom).
   \label{figure:experiments:sin:cost}
 }
\end{figure}

\begin{observation}
  Our experiments validate at hands of the \new{$\texttt{sin} $ benchmark} that all three
  algorithm variants yield multigrid behaviour.
\end{observation}

\noindent
Next we compare the residual reduction to unknown reads from memory.
The latter \new{scale with} $touchVertexFirstTime$ counts.
Block operations and accumulations all are expected to happen in the cache due
to the vertical integration of the algorithm's phases within one tree sweep
\cite{Weinzierl:11:Peano}.
We validate this theoretic statement experimentally in
Sect.~\ref{section:results:parallel}.
Multiplicative solvers are superior to the other solver variants. 
However, the vast difference in iteration counts does not translate directly
into speed/data reads---the difference in actual runtime is smaller
(Fig.~\ref{figure:experiments:sin:cost}).
All variants yield cost in the same order of magnitude.

\begin{observation}
  For the smooth \new{$\texttt{sin} $} setup, multiplicative multigrid is not
  significantly superior to additive variants in terms of memory access cost.
\end{observation}

\noindent
We observe an opposite effect regarding cost vs.~iteration count in Fig.~\ref{figure:experiments:sin:cost}:
The dynamic adaptivity unfolds the grid in a
full multigrid (FMG) way for BPX and the additive solver and naturally yields
FMG \new{character lacking} higher-order operators.
This decreases the time-to-solution \new{yet} increases the iteration count.


\begin{observation}
  For a very smooth setup such as \new{in the $\texttt{sin}$ benchmark}, a dynamic refinement
  criterion unfolding the grid from a coarse start solution yields almost F-cycle-like
  convergence even in the absence of higher order interpolation. 
\end{observation}

\subsection{Jumping and anisotropic material parameters $\epsilon$}
\label{section:results:convergence:jump}

We continue with \new{the $\texttt{jump}$ setup, where} the material
parameter \new{changes} in the middle of the domain.
The \new{change/}jump is not aligned with the grid.
For the geometric multigrid variants, 
any coarse-grid update's $d$-linear interpolation $Pd_{\ell -1}$ that
overlaps the parameter jump does not anticipate the lack of smoothness in the
solution and, thus, introduces a localized fine-grid error
around the material parameter transition.
If we project update from left and right of the jump, the update's linear
interpolation lacks the discontinuity in the derivative, it pollutes the solve around
the $\epsilon $ jump.
For reasonable big $\omega $, the solvers start to oscillate locally.

\begin{figure}
 \begin{center}
  %
  %
  %
  \includegraphics[width=0.4\textwidth]{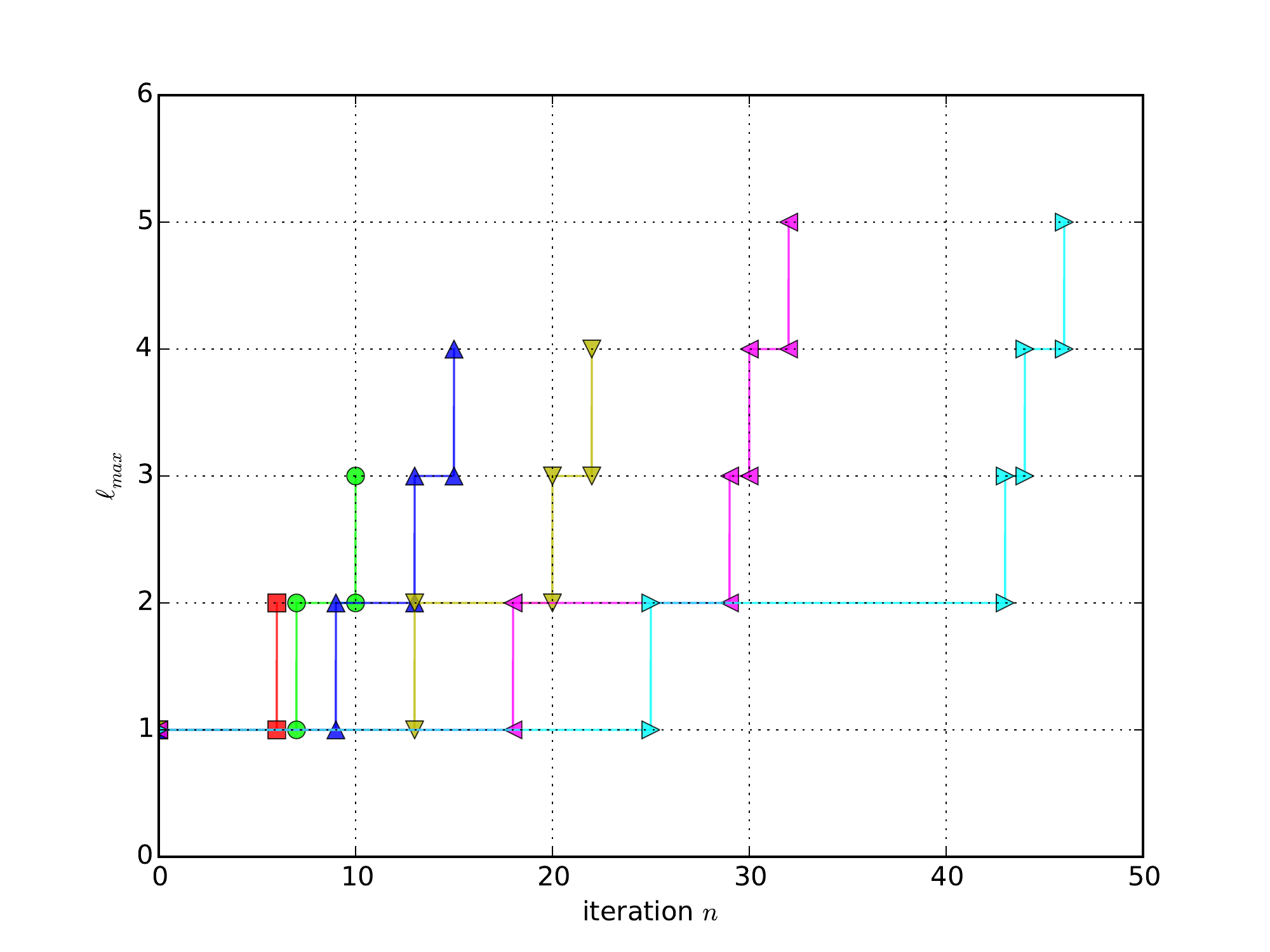}
  \includegraphics[width=0.4\textwidth]{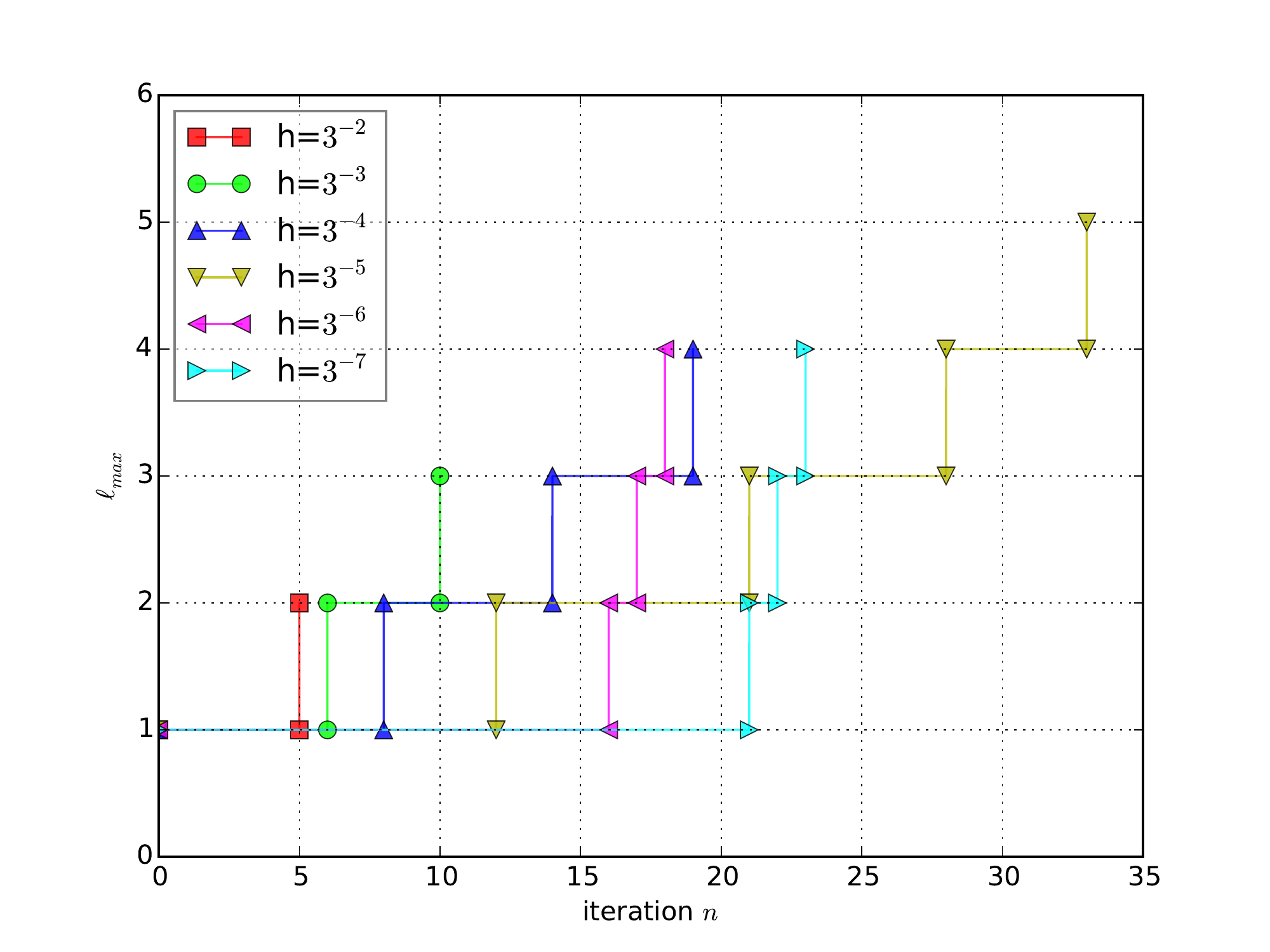}
 \end{center}
 \caption{
   Development of $\ell _{max}$ over the number of iterations if the algorithm
   is allowed to increase the coarsest level whenever the residual for the $2d$
   version of \new{the $\texttt{jump}$ benchmark} starts to grow. 
   Data for geometric BPX with the \texttt{BlockJacobi(4)} smoother (left) and 
   multiplicative V11 with Jacobi smoothers (right).
   \label{figure:experiments:jump:lmax}
 }
\end{figure}

The problem can be tackled by stronger smoothers or higher $\mu$
counts applied to the finest grids, or we can use smaller $\omega $.
Both approaches harm multigrid performance. 
They furthermore suffer from the fact that it is often not clear which
$\ell _{max}$ still yields a robust solver for a particular setup.
Our code family identifies non-diagonal dominant operators, stagnation or
amplifying oscillations, increases the coarsest mesh level $\ell _{max}$ autonomously and 
thus avoids some instabilities (Fig.~\ref{figure:experiments:jump:lmax}).
Yet,  we are not able to solve any setup from \new{the $\texttt{jump}$ benchmark} with
less than 300 iterations with Jacobi once $h<3^{-2}$.
For $h=3^{-2}$, the additive multigrid with a 4-sweep block smoother
requires already 133 iterations---it \new{deteriorates}.
In general, the $\ell _{max}$ modifications remove more coarse grid resolutions
from the additive and BPX schemes than for the multiplicative multigrid, 
while the increase of $\ell _{max}$ does not kick in for any Galerkin solver
variant:

\begin{observation}
  An adaptive coarse grid choice mitigates the effect of the absence of Galerkin coarse grid
  operators w.r.t.~stability and, at the same time, identifies the
  coarsest resolution level where geometric multigrid remains robust.
  \new{With respect to the performance}, it remains a workaround in the absence
  of proper \new{unstructured coarse grids or algebraic/direct coarse grid
  solves}.
\end{observation}

\begin{table}
  \tbl{
    Iteration counts for the additive multigrid, additive multigrid with exact
    coarse-grid solve, and BPX (left to right) with $\omega = 0.8$ solving
    \new{the $\texttt{jump}$ benchmark} for $d=2$.
    $\bot $ is used to denote that a solver was not able
    to reduce the residual by a factor of $10^{8}$ within 300 iterations.
    \label{results:jump:additive}
  }
{
\footnotesize
 \setlength{\tabcolsep}{2.4pt}
 \begin{tabular}{l|ccccc|ccccc|ccccc}
   %
   %
   %
   %
   h & Jac & BJ(1) & BJ(2) & BJ(4) & BJ(8) & Jac/e & BJ(1)/e & BJ(2)/e & BJ(4)/e & BJ(8)/e & Jac & BJ(1) & BJ(2) & BJ(4) & BJ(8) \\
 \hline  
$3^{-2}$ & 62      & 41  & 36  & 35  & 34  & $\bot $ & 41  & 36  & 35  & 34  & 62      & 41 & 36 & 35 & 34 \\
$3^{-3}$ & 90      & 58  & 52  & 49  & 49  & $\bot $ & 57  & 52  & 49  & 49  & 76      & 52 & 47 & 44 & 44 \\
$3^{-4}$ & 131     & 81  & 72  & 68  & 67  & $\bot $ & 79  & 71  & 68  & 67  & 96      & 67 & 60 & 57 & 56 \\
$3^{-5}$ & 179     & 106 & 93  & 88  & 87  & $\bot $ & 103 & 93  & 88  & 87  & $\bot $ & 81 & 72 & 68 & 67 \\
$3^{-6}$ & 201     & 120 & 108 & 103 & 102 & $\bot $ & 119 & 108 & 103 & 102 & $\bot $ & 85 & 77 & 73 & 72 \\
$3^{-7}$ & $\bot $ & 123 & 112 & 107 & 106 & $\bot $ & 123 & 112 & 107 & 106 & $\bot $ & 80 & 73 & 70 & 69 \\

 \end{tabular}
}
\end{table}

%
%
%

\begin{table}
  \tbl{
    Iteration counts for \new{the $\texttt{jump}$} and $d=2$ with multiplicative V21-cycle.
    $\bot $ denotes that a solver is not able
    to reduce $|r|_2$ by a factor of $10^{8}$ within 300 iterations.
    \label{results:jump:MGV21}
  }
{
\footnotesize
  \setlength{\tabcolsep}{2.4pt}
  \begin{tabular}{l|ccccc|ccccc}
   %
   %
   %
   %
   h & Jac & BJ(1) & BJ(2) & BJ(4) & BJ(8) & Jac/e & BJ(1)/e & BJ(2)/e & BJ(4)/e & BJ(8)/e \\
 \hline  
$3^{-2}$ & 21 & 14 & 12 & 12 & 12 & $\bot $ & 14 & 12 & 12 & 12 \\
$3^{-3}$ & 27 & 20 & 18 & 18 & 18 & $\bot $ & 20 & 18 & 18 & 18 \\
$3^{-4}$ & 35 & 27 & 25 & 24 & 24 & $\bot $ & 27 & 25 & 24 & 24 \\
$3^{-5}$ & 42 & 34 & 32 & 31 & 30 & $\bot $ & 34 & 32 & 31 & 30 \\
$3^{-6}$ & 42 & 35 & 33 & 32 & 32 & $\bot $ & 35 & 33 & 32 & 32 \\
$3^{-7}$ & $\bot $ & 33 & 31 & 30 & 30 & $\bot $ & 33 & 31 & 30 & 30 \\

  \end{tabular}
}
\end{table}

\begin{figure}
  \begin{center}
  \includegraphics[width=0.3\textwidth]{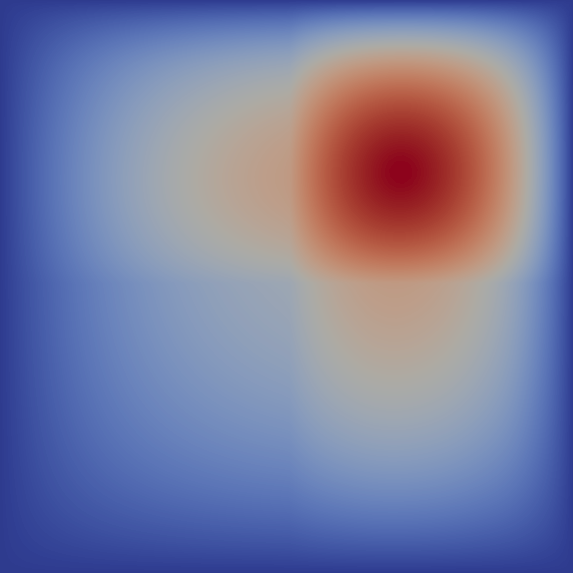}
  \includegraphics[width=0.3\textwidth]{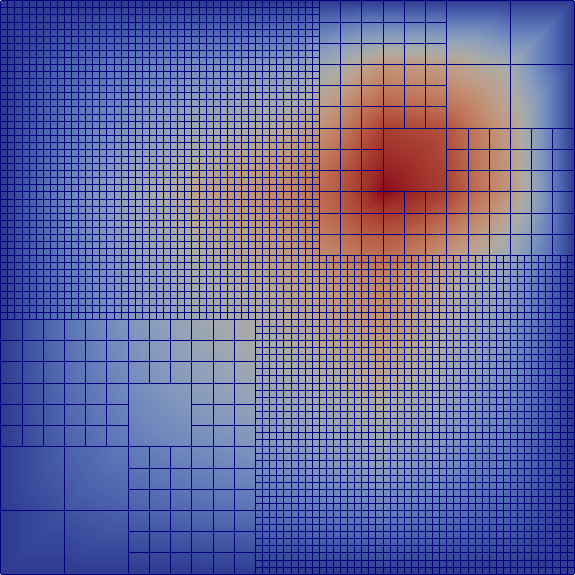}
  \includegraphics[width=0.3\textwidth]{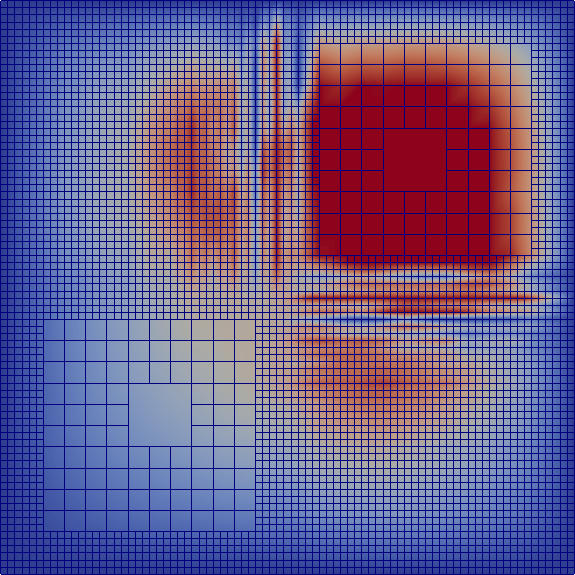}
  \end{center}
  \caption{
    From left to right: Solution of \new{the $\texttt{checkerboard}$ setup} for $d=2$. 
    Typical adaptive grid.
    Instability pattern arising from additive geometric multigrid with
    rediscretization after 28 iterations.
    \label{figure:jump:solutions}
  }
\end{figure}

\noindent
Experiments reveal that the problem can be solved in around 107 iterations
($h=3^{-7}$) if we use Galerkin coarse-grid operators and block smoothers.
BPX converges in 35 ($h=3^{-2}$) to 70 ($h=3^{-7}$) iterations
(Table \ref{results:jump:additive}).
We recognize that a grid spacing reduction from $3^{-6}$ to
$3^{-7}$ reduces the total iteration count here:
in regions of interest the grid is refined aggressively, but if the maximum level
is too constrained, these regions spread out and increase the vertex count
unnecessarily.
This anomaly carries over to multiplicative multigrid (Table
\ref{results:jump:MGV21}), which now clearly
outperforms the other two solvers.
An exact coarse-grid solve \new{now} does not pay off anymore.
Though the solver remains stable with $\ell _{max}=1$,  the coarsest
level can not contribute with any useful correction as it is too coarse.
As an improvement, one could choose a larger $\ell
_{max}$ right from the start and solve exactly there.

\begin{table}
 \tbl{
  Cost (number of unknown reads) of experiments of Table
  \ref{results:jump:MGV21} with regular grid (top) and the dynamically
  adaptive grid (bottom).
  \label{results:jump:MGV21-cost}
 }
 {
 \footnotesize
   \begin{tabular}{l|ccccc}
    %
    %
    %
    %
    h & Jac & BJ(1) & BJ(2) & BJ(4) & BJ(8) \\
 \hline  
$3^{-2}$ & 5.90e+03 & 4.02e+03 & 3.48e+03 & 3.48e+03 & 3.48e+03 \\
$3^{-3}$ & 8.14e+04 & 6.11e+04 & 5.53e+04 & 5.53e+04 & 5.53e+04 \\
$3^{-4}$ & 1.00e+06 & 7.79e+05 & 7.24e+05 & 6.96e+05 & 6.96e+05 \\
$3^{-5}$ & 1.10e+07 & 8.95e+06 & 8.44e+06 & 8.18e+06 & 7.93e+06 \\
$3^{-6}$ & 9.96e+07 & 8.34e+07 & 7.88e+07 & 7.65e+07 & 7.65e+07 \\

    \hline  
    $3^{-2}$ & $\bot $ & 4.02e+03 & 3.48e+03 & 3.48e+03 & 3.48e+03 \\
$3^{-3}$ & $\bot $ & 4.02e+03 & 3.48e+03 & 3.48e+03 & 3.48e+03 \\
$3^{-4}$ & $\bot $ & 3.58e+04 & 3.37e+04 & 3.17e+04 & 3.17e+04 \\
$3^{-5}$ & $\bot $ & 6.35e+06 & 6.35e+06 & 6.31e+06 & 6.31e+06 \\
$3^{-6}$ & $\bot $ & $\bot $ & 5.12e+06 & 5.10e+06 & 5.10e+06 \\

   \end{tabular}
 }
 \end{table}

\begin{table}
  \tbl{
    Number of iterations to solve \new{the $\texttt{checkerboard}$ benchmark} with
    additive multigrid (left,middle) and BPX (right) for $d=2$.
    All inter-grid transfer operators are bilinear, coarse-grid
    operators realize the Galerkin idea.
    \label{results:checkerboard:bpx}
  }
{
\footnotesize
   %
   %
   %
   %
 \setlength{\tabcolsep}{2.4pt}
 \begin{tabular}{l|ccccc|ccccc|ccccc}
   h & Jac & BJ(1) & BJ(2) & BJ(4) & BJ(8) & Jac/e & BJ(1)/e & BJ(2)/e & BJ(4)/e & BJ(8)/e & Jac & BJ(1) & BJ(2) & BJ(4) & BJ(8)\\
 \hline  
$3^{-2}$ & 120 & 67 & 57 & 52 & 51 & 120 & 68 & 57 & 52 & 51 & 115 & 67 & 57 & 52 & 51 \\
$3^{-3}$ & 270 & 117 & 91 & 80 & 78 & 270 & 118 & 91 & 79 & 78 & $\bot $ & 116 & 91 & 79 & 78\\
$3^{-4}$ & $\bot $ & 140 & 112 & 100 & 98 & $\bot $ & 137 & 111 & 100 & 98 & $\bot $ & 136 & 107 & 94 & 92 \\
$3^{-5}$ & $\bot $ & 163 & 132 & 118 & 116 & $\bot $ & 158 & 129 & 118 & 116& $\bot $ & 146 & 116 & 102 & 100 \\
$3^{-6}$ & $\bot $ & 167 & 140 & 130 & 129 & $\bot $ & 166 & 140 & 130 & 129 & $\bot $ & 137 & 113 & 102 & 101\\
$3^{-7}$ & $\bot $ & 156 & 135 & 127 & 125 & $\bot $ & $\bot $ & $\bot $ & $\bot $ & $\bot $ & $\bot $ & 119 & 101 & 92 & 91\\

 \end{tabular}
}
\end{table}

Dynamic adaptivity again pays off (Table \ref{results:jump:MGV21-cost})
and reduces the number of unknown reads by an order of magnitude:
the grid quickly refines around the material transition and, thus, injects the
critical behaviour into the coarse-grid corrections.
At the same time it acts as FMG facilitator.
Empty entries in Table \ref{results:jump:MGV21-cost} result from the fact that
we always stop after 300 iterations; an iteration count that is
quickly met if grid ``setup'' iterations are counted as part of an FMG cycle.


\begin{table}
  \tbl{
    Experiments from Table \ref{results:checkerboard:bpx} with
    multiplicative V21 cycle.
    \label{results:checkerboard:mg}
  }
{
\footnotesize
   %
   %
   %
   %
  \begin{tabular}{l|ccccc|ccccc}
    h & Jac & BJ(1) & BJ(2) & BJ(4) & BJ(8) & Jac/e & BJ(1)/e & BJ(2)/e & BJ(4)/e & BJ(8)/e \\
 \hline  
$3^{-2}$ & 44      & 23 & 18 & 16 & 16 & 44      & 23 & 18 & 16 & 16  \\
$3^{-3}$ & $\bot $ & 40 & 32 & 29 & 29 & $\bot $ & 40 & 32 & 29 & 29 \\
$3^{-4}$ & $\bot $ & 49 & 42 & 39 & 39 & $\bot $ & 49 & 42 & 39 & 39 \\
$3^{-5}$ & $\bot $ & 58 & 52 & 50 & 49 & $\bot $ & 58 & 52 & 50 & 49 \\
$3^{-6}$ & $\bot $ & 56 & 51 & 50 & 49 & $\bot $ & 56 & 51 & 50 & 49 \\
$3^{-7}$ & $\bot $ & 50 & 47 & 46 & 45 & $\bot $ & 50 & 47 & 46 & 45 \\

  \end{tabular}
}
\end{table}

Our follow-up experiments \new{with the $\texttt{checkerboard}$ setup} continue to have
$\epsilon $ jumps, but furthermore introduce anisotropic regions.
Anisotropic behaviour poses even harder \new{challenges} to $d$-linear
inter-grid transfer operators:
All geometric solvers are ill-suited for this problem, they do not
converge (Fig.~\ref{figure:jump:solutions}).
While the regular Galerkin solvers converge, mesh-independent convergence is
lost completely (Tables
\ref{results:checkerboard:bpx},\ref{results:checkerboard:mg}).
Again, an exact coarse-grid solve in the multiplicative setting is not required: 
no matter how exact the coarsest problem is solved, the prolongation of the
correction always yields wrong fine-grid modes for varying $\epsilon $ or
anisotropic $\epsilon _i \not= \epsilon _j$.

 \begin{table}
   \tbl{
       Cost for the experiments of Table
      \ref{results:checkerboard:mg} with regular (left) or dynamically
      adaptive grid (right).
     \label{results:checkerboard:MGV21-cost}
   }
 {
 \footnotesize
  \setlength{\tabcolsep}{2.4pt}
  \begin{tabular}{l|ccccc|ccccc}
   %
   %
   %
   %
   h & Jac & BJ(1) & BJ(2) & BJ(4) & BJ(8) & Jac & BJ(1) & BJ(2) & BJ(4) & BJ(8) \\
 \hline  
$3^{-2}$ & 1.21e+04 & 6.43e+03 & 5.09e+03 & 4.56e+03 & 4.56e+03 
         & 1.21e+04 & 6.43e+03 & 5.09e+03 & 4.56e+03 & 4.56e+03 \\
$3^{-3}$ & $\bot $ & 1.19e+05 & 9.60e+04 & 8.72e+04 & 8.72e+04 
         & 1.21e+04 & 6.43e+03 & 5.09e+03 & 4.56e+03 & 4.56e+03 \\
$3^{-4}$ & $\bot $ & 1.39e+06 & 1.20e+06 & 1.11e+06 & 1.11e+06 
         & $\bot $ & 7.63e+04 & 5.80e+04 & 5.05e+04 & 5.05e+04 \\
$3^{-5}$ & $\bot $ & 1.51e+07 & 1.36e+07 & 1.30e+07 & 1.28e+07 
         & $\bot $ & 1.13e+06 & 9.98e+05 & 9.79e+05 & 9.63e+05 \\
$3^{-6}$ & $\bot $ & 1.32e+08 & 1.20e+08 & 1.18e+08 & 1.16e+08 
         & $\bot $ & 4.63e+06 & 4.59e+06 & 4.55e+06 & 4.55e+06 

  \end{tabular}
 }
 \end{table}


The dynamic adaptivity criterion continues to make a dynamic approach outperform
its regular-grid counterpart in terms of cost
(Table \ref{results:checkerboard:MGV21-cost}), while the actual number of grid
sweeps is higher by a factor of four.
Yet, the sweeps are cheap as long as the grid has not unfolded substantially.
\new{Due to} the irregularity of \new{$\texttt{checkerboard}$}, the
feature-based adaptivity refines
anisotropic regions, regions around $\epsilon $ changes and along the
boundary (Fig.~\ref{figure:jump:solutions}).

\begin{observation}
 \new{
  Our global, uniform choice of a smoothed grid level where fine grid regions
  coarser than the current smoothing level are updated per cycle, too,  
  plus the dynamic
  refinement criterion being active in each cycle imply that rough solution
  regions are subject to immediate refinement.
  At the same time, smooth regions resolved by a rather coarse grid are subject
  to more smoothing steps.
 }
\end{observation}


\subsection{BoxMG}

\begin{figure}[htb]
  \begin{center}
    \includegraphics[width=0.49\textwidth]{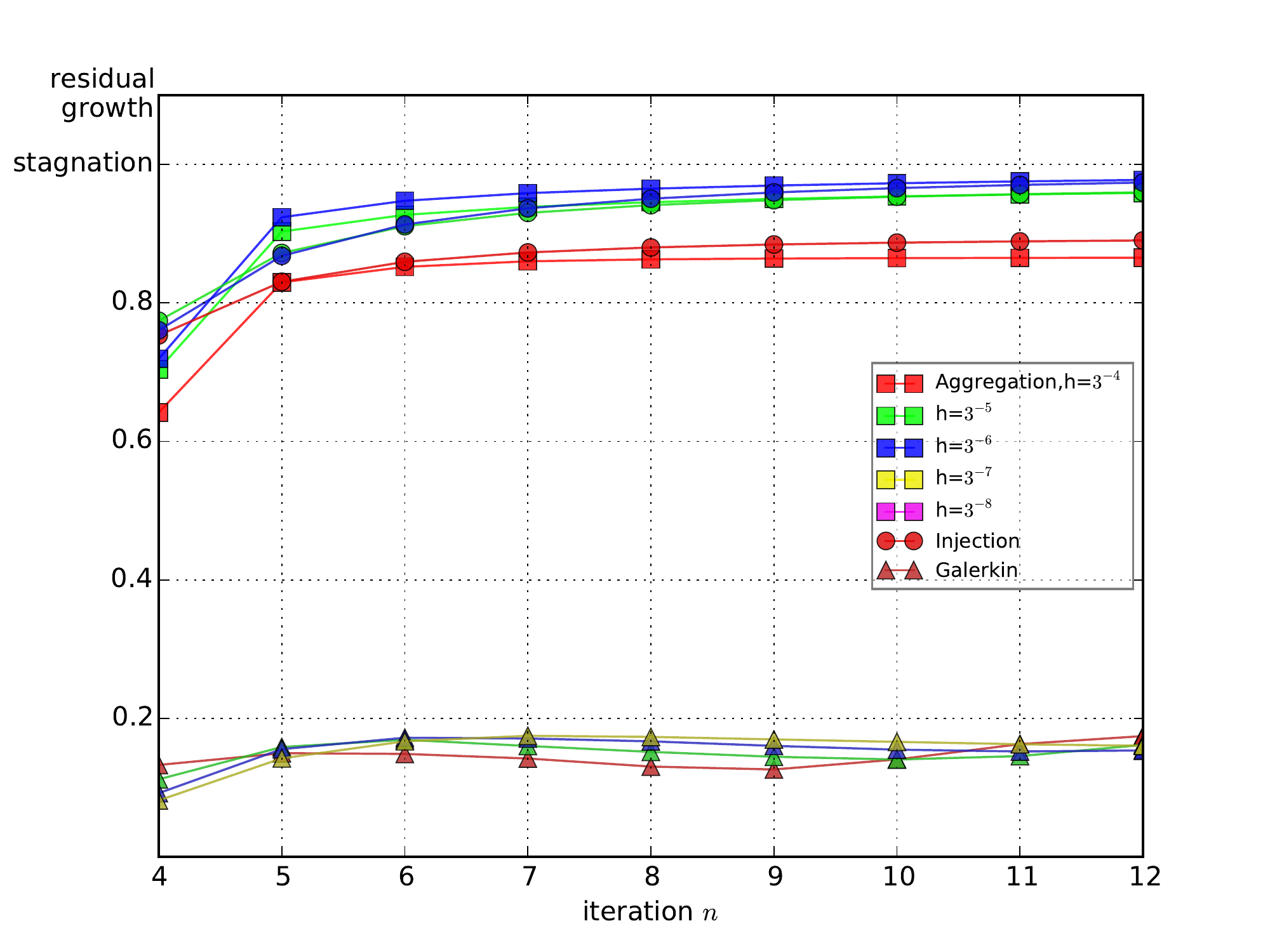}
    \includegraphics[width=0.49\textwidth]{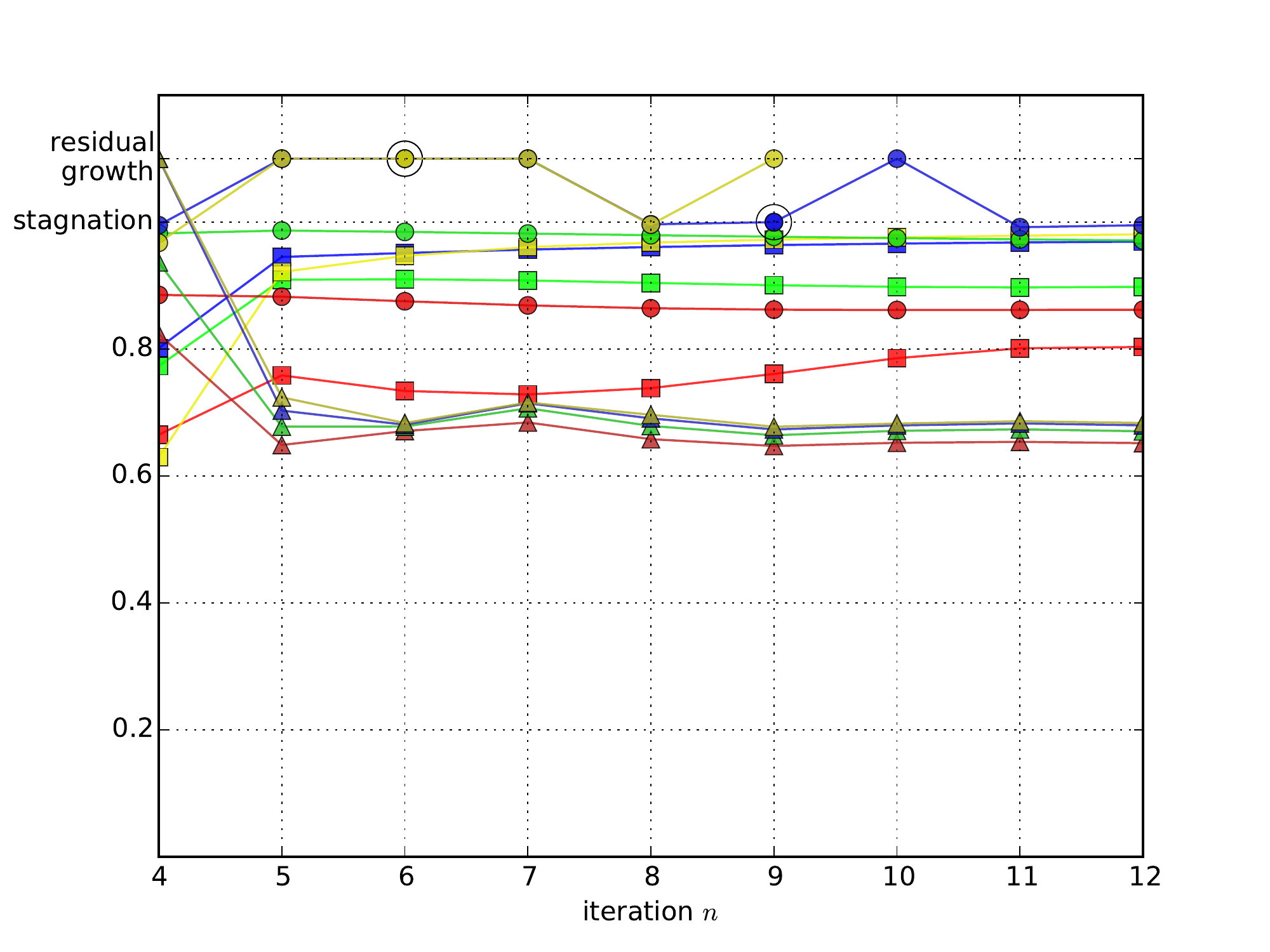}
    \\
    \includegraphics[width=0.49\textwidth]{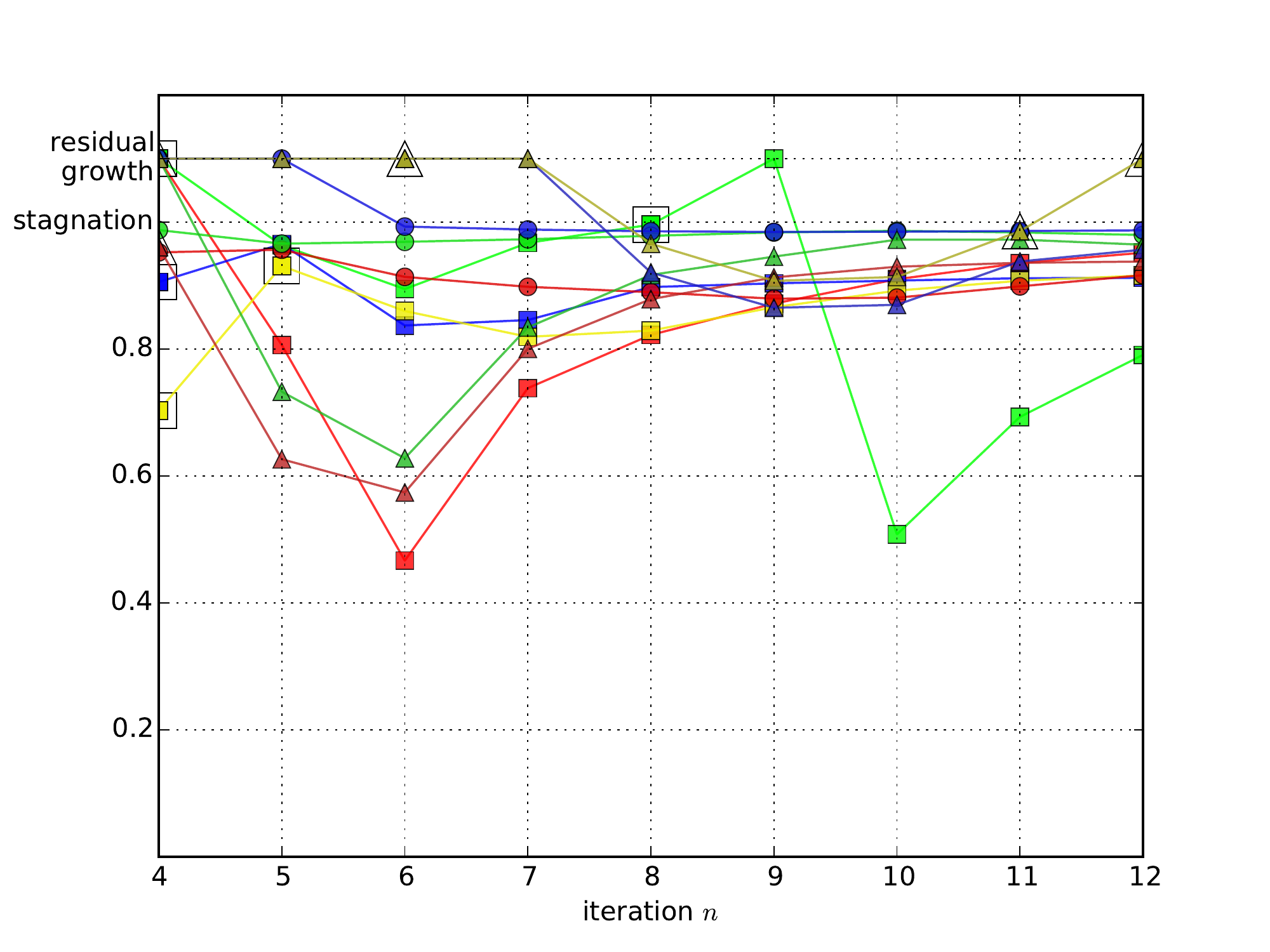}
    \includegraphics[width=0.49\textwidth]{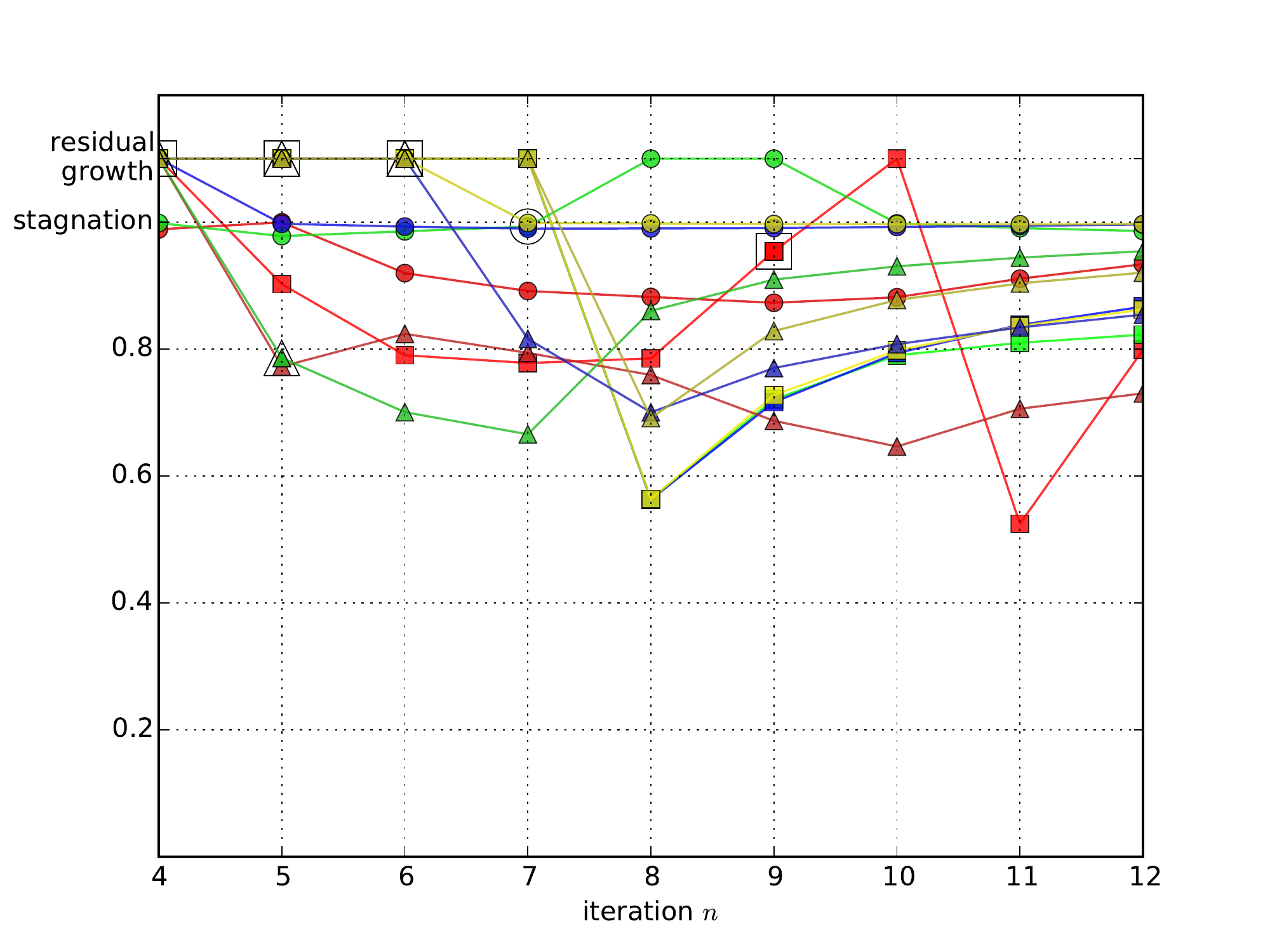}
  \end{center}
  \caption{
    Convergence speed $\rho = |res(n)|_2/|res(n-1)|_2$ as ratio between current
    residual to residual of previous iteration.
    The experiments pass $\epsilon = \{ 10^{-1},10^{-2},10^{-4},10^{-8} \}$
    (left to right, top-down) into a V21 multigrid with a block Jacobi smoother
    running four iterations per patch. 
    Large empty symbols indicate that the solver
    increases $\ell _{max}$ in that iteration.
    \label{results:convection-2d:rho}
  }
\end{figure}

The recirculating flow problem \new{$\texttt{circle}$ (Table
\ref{table:parameters}}) with inflow
boundary conditions from the left and right side (Fig.~\ref{figure:shared-memory:concurrency}) introduces non-zero
convection which destroys the symmetry of $A$.
The smaller the (symmetric) diffusion \new{weight} $\epsilon $ of the
operator relative to the convection the more
challenging the setup becomes for multigrid solvers.
They struggle with convection-dominated systems
(Fig.~\ref{figure:results:circle-2d}).
In our setup, the convection coefficients furthermore vary in space and there is
a singularity in the middle of the domain. 
Finally, the closed characteristics prevent any solver
error from being ``pushed out'' of the domain by relaxation. 
Geometric multigrid and Galerkin multigrid with $d$-linear inter-grid
transfer operators solve this problem if we use dynamic coarse-grid adaptation
and proper upwinding.
However, $\ell _{max}$ increases fast and, hence, the solvers
\new{transition} quickly into pure (block) Jacobi.

Once we supplement the Galerkin coarse grid operator computation with BoxMG, we
obtain a more robust solver.
To show this, we track the convergence speed $\rho (n)$ as ratio of the
residual norm in an iteration $n$ relative to the residual norm in the iteration
$n-1$ (Fig.~\ref{results:convection-2d:rho}).
Furthermore, we compare plain BoxMG with a Petrov-Galerkin scheme where $R=I$
(injection) or $R$ is an aggregation operator.
We observe the convergence speed deteriorates with decreasing $\epsilon $,
i.e.~with an increasing impact of the convection, for all solvers.
We start \new{loosing} multigrid behaviour.
\new{Convergence} speed depends
on the mesh size, and measurements for one type of solver (same symbol) fan out.
For $\epsilon = 10^{-1}$ the automatic coarse grid increase is not triggered.
The solvers work with $\ell _{max} = 1$ all the time.
The smaller $\epsilon $ the more often $\ell $ is increased. 
Injection yields consistently the worst convergence speed and is thus not
studied further.
We can not identify a significant reduction of the number of $\ell _{max}$ increases if we
switch from pure BoxMG to the symmetric aggregation-based restriction.
While literature would expect us to obtain more stable coarse grid operators
(cmp.~the discussion in \cite{yavnehbendig12nonnsymbb} and references therein),
we do not observe this effect here---\new{probably due to a lack of a higher order symmetric operator}.
However, we do observe that aggregation yields slightly better convergence rates
than a pure Galerkin approach with $P=R^T$ if the grid is sufficiently fine and
$\epsilon $ is small enough.

\begin{table}
  \tbl{
    \new{Overview of the solver variants and their performance for the different
    problem setups. A checkmark indicates that the problem can be solved by the
    respective solver with multigrid performance, a checkmark in brackets
    indicates that the solver converges, but mesh-independent (multigrid)
    performance is lost. A cross denotes that the solver diverges for the
    problem. For \texttt{circle}, the behaviour depends strongly on the weight
    of the convection term relative to the diffusion. }
    \label{table:solvers}
  }{
  {\footnotesize 
   \begin{tabular}{l|cccc}
       & \texttt{sin} & \texttt{jump} & \texttt{checkerboard} 
       & \texttt{circle} \\
     \hline 
      Geo.~Add.& $\checkmark$ & $(\checkmark)$ & $\times$ & $\times$ \\
      Geo.~BPX& $\checkmark$ & $(\checkmark)$ & $\times$ & $\times$ \\
      Geo.~Mult.& $\checkmark$ & $(\checkmark)$ & $\times$ & $\times$ \\
      Galerkin.~Add.& $\checkmark$ & $\checkmark$ & $(\checkmark)$ & $\times$ \\
      Galerkin.~BPX& $\checkmark$ & $\checkmark$ & $(\checkmark)$ & $\times$ \\
      Galerkin.~Mult.& $\checkmark$ & $\checkmark$ & $(\checkmark)$ &  $\times$
      \\
      BoxMG Add.& $\checkmark$ & $\checkmark$ & $\checkmark$ &
      $\checkmark\ \mbox{to}\ (\checkmark)$ \\
      BoxMG BPX& $\checkmark$ & $\checkmark$ & $\checkmark$ & $\checkmark\ \mbox{to}\ (\checkmark)$ \\
      BoxMG Mult.& $\checkmark$ & $\checkmark$ & $\checkmark$ & $\checkmark\ \mbox{to}\ (\checkmark)$
      \\
    \end{tabular}
    }
  }  
\end{table}

A study of the exact behaviour of the various multigrid compositions \new{(Table
\ref{table:solvers})} is beyond scope here though we state that the BoxMG
solvers are reasonably robust to tackle non-trivial problems.
Yet, the experiments also illustrate that stronger
smoothers such as ILU, Kaczmarz or alternating line Gauss-Seidel are highly
desirable for dominating convection.
Keeping this in mind, our ideas can act as a reasonable code building block for robust solvers
on spacetrees that work strictly element-wise and matrix-free.
\new{
All solvers implement a single-touch policy, i.e., one tree traversal realizes
one multigrid cycle (additive and BPX) or smoothing step (multiplicative).
In practice, this means that a multiplicative cycle is by a factor of $(\mu
_{pre}+\mu _{post})$ times slower than an additive cycle.
Both the problem character and the multigrid's role---is it used as solver or
as preconditioner---determine which approach is superior in terms of time to
solution.
A generic ``better than'' statement is impossible to make.
}

\subsection{Memory consumption}

\begin{table}
  \tbl{
    Multiplicative BoxMG with dynamic
    adaptivity criterion for $d=2$ and \new{$\texttt{sin}$} (top), 
    \new{$\texttt{jump}$} and \new{$\texttt{circle}$} (bottom).
    We \new{present one tuple per run}.
    \new{The left entry denotes} how many vertex updates are 
    required to reduce the residual to $10^{-8}$\new{, while the right
    entry} gives the
    \new{memory footprint with compressed algebraic operators} relative
    to the uncompressed run.
    \label{results:memory:sin}
  }
  {
  \footnotesize
  \begin{tabular}{l|cccc} 
   h & 0 & 1e-2 & 1e-4 & 1e-8 \\
\hline 
$3^{-3}$ & 2.24e+03 \new{ / 1.0} & 2.24e+03 / 1.33e-02 & 2.24e+03 / 1.99e-02 & 2.24e+03 / 3.31e-02 \\
$3^{-4}$ & 1.36e+05 \new{ / 1.0} & 1.27e+05 / 5.68e-03 & 2.06e+05 / 8.35e-03 & 1.36e+05 / 1.38e-02 \\
$3^{-5}$ & 1.37e+05 \new{ / 1.0} & 2.97e+05 / 4.28e-03 & 4.57e+05 / 5.12e-03 & 1.37e+05 / 1.39e-02 \\
$3^{-6}$ & 2.09e+06 \new{ / 1.0} & 2.92e+06 / 2.83e-03 & 2.09e+06 / 4.39e-03 & 2.09e+06 / 7.29e-03 \\
$3^{-7}$ & 2.09e+06 \new{ / 1.0} & 4.80e+06 / 2.82e-03 & 2.09e+06 / 4.39e-03 & 2.09e+06 / 7.29e-03 \\
  
   \hline 
   $3^{-3}$ & 2.65e+03 \new{/ 1.0} & 1.02e+03 / 1.63e-02 & 1.84e+03 / 2.48e-02 & 2.65e+03 / 4.10e-02 \\
$3^{-4}$ & 2.81e+04 \new{/ 1.0} & 4.41e+05 / 9.57e-03 & 2.51e+04 / 1.47e-02 & 2.81e+04 / 2.39e-02 \\
$3^{-5}$ & 4.54e+06 \new{/ 1.0} & 2.72e+05 / 7.90e-03 & 4.54e+06 / 8.41e-03 & 4.54e+06 / 1.35e-02 \\
$3^{-6}$ & 4.13e+06 \new{/ 1.0} & 2.72e+05 / 7.90e-03 & 4.13e+06 / 6.65e-03 & 4.13e+06 / 1.07e-02 \\
$3^{-7}$ & 4.17e+06 \new{/ 1.0} & 2.72e+05 / 7.90e-03 & 4.18e+06 / 6.72e-03 & 4.17e+06 / 1.08e-02 \\
  
   \hline 
   $3^{-3}$ & 1.08e+04 \new{/ 1.0} & 5.10e+03 / 1.68e-02 & 7.75e+03 / 2.53e-02 & 2.33e+04 / 4.10e-02 \\
$3^{-4}$ & 2.43e+04 \new{/ 1.0} & 2.41e+04 / 9.57e-03 & 2.43e+04 / 1.42e-02 & 2.43e+04 / 2.25e-02 \\
$3^{-5}$ & 3.43e+05 \new{/ 1.0} & 2.05e+05 / 1.20e-02 & 3.43e+05 / 1.64e-02 & 3.43e+05 / 2.67e-02 \\
$3^{-6}$ & 1.84e+07 \new{/ 1.0} & 1.78e+07 / 1.14e-02 & 1.60e+06 / 1.62e-02 & 1.84e+07 / 2.51e-02 \\
$3^{-7}$ & 1.45e+07 \new{/ 1.0} & 1.86e+07 / 1.14e-02 & 1.60e+06 / 1.62e-02 & 1.45e+07 / 2.56e-02 \\
  
  \end{tabular}
  }
\end{table}

To quantify the memory \new{demands}, we study dynamically adaptive
grids with a multiplicative V11-BoxMG solver.
Our measurements compare an uncompressed version with $\epsilon _{mf}\in \{
10^{-2}, 10^{-4}, 10^{-8} \} $.
The memory savings are enormous for the \new{$\texttt{sin}$} benchmark where
all operators are well-known to reduce to their rediscretization or $d$-linear
counterpart.
%
%
%
%
%
With decreasing $h$, the code becomes matrix-free.

\begin{table}
  \tbl{
    \new {Characteristic runtimes per V22 cycle for a regular (left) and
    an adaptive (right) mesh on a single core
    for the \texttt{checkerboard} setup.
    }
    \label{results:memory:cost}
  }
  {
  \footnotesize
  \begin{tabular}{l|rr|rr} 
$h_{min}$ & no compr. & compr. & no compr. & compr. \\
\hline
$3^{-4}$ & $1.54 \cdot 10^{-2}$ & $2.03 \cdot 10^{-3}$ 
   & $1.08 \cdot 10^{-2}$ &  $1.55 \cdot 10^{-2}$	 \\ 
$3^{-5}$ & $6.28 \cdot 10^{-2}$ & $8.92 \cdot 10^{-2}$  
   & $2.31 \cdot 10^{-2}$ &  $3.49 \cdot 10^{-2}$	 \\ 
$3^{-6}$ & $3.89 \cdot 10^{-1}$ & $5.26 \cdot 10^{-1}$  
   & $0.56 \cdot 10^{-1}$ &  $1.01 \cdot 10^{-1}$	 \\ 
$3^{-7}$ & $2.74 \cdot 10^{-0}$ & $3.53 \cdot 10^{-0}$  
   & $2.42 \cdot 10^{-1}$ &  $1.53 \cdot 10^{-1}$	 \\ 
  \end{tabular}
  }
\end{table}

For \new{the $\texttt{jump}$ and $\texttt{circle}$ setups}, we preserve
significant memory savings, though
they are by a factor of two to four smaller than for the \texttt{sin}
setup.
We may choose rather small
$\epsilon _{mf} = 10^{-8}$ and nevertheless obtain both high compression rates
and preserve the solvers' semantics. 
If the compression is too aggressive, the adaptivity criterion yields slightly
different adaptivity patterns, presumably through inexact arithmetics.

Stencils plus the unknowns from Table \ref{table:unknowns} yield a
memory footprint of $3+3^d+2\cdot 5^d$ doubles plus another byte for grid management \cite{Weinzierl:11:Peano}.
This already is a small memory footprint for dynamically adaptive grids.
\new{Compression reduces} the average
footprint to close to $3$ doubles plus a byte per vertex without loosing 
algebraic multigrid operators.
\new{Yet, its cost, on a single core, is not negligible
(Table \ref{results:memory:cost}) unless the grid is extremely adaptive.}

\begin{observation}
 \new{
 The reduction of the solver's memory footprint to almost the pure footprint of
 a purely geometric approach can double the computational cost on one
 core.
 }
\end{observation}

\subsection{Comparisons to PETSc's GAMG}
\label{section:results:petsc}

\begin{figure}
 \begin{center}
  \includegraphics[width=0.45\textwidth]{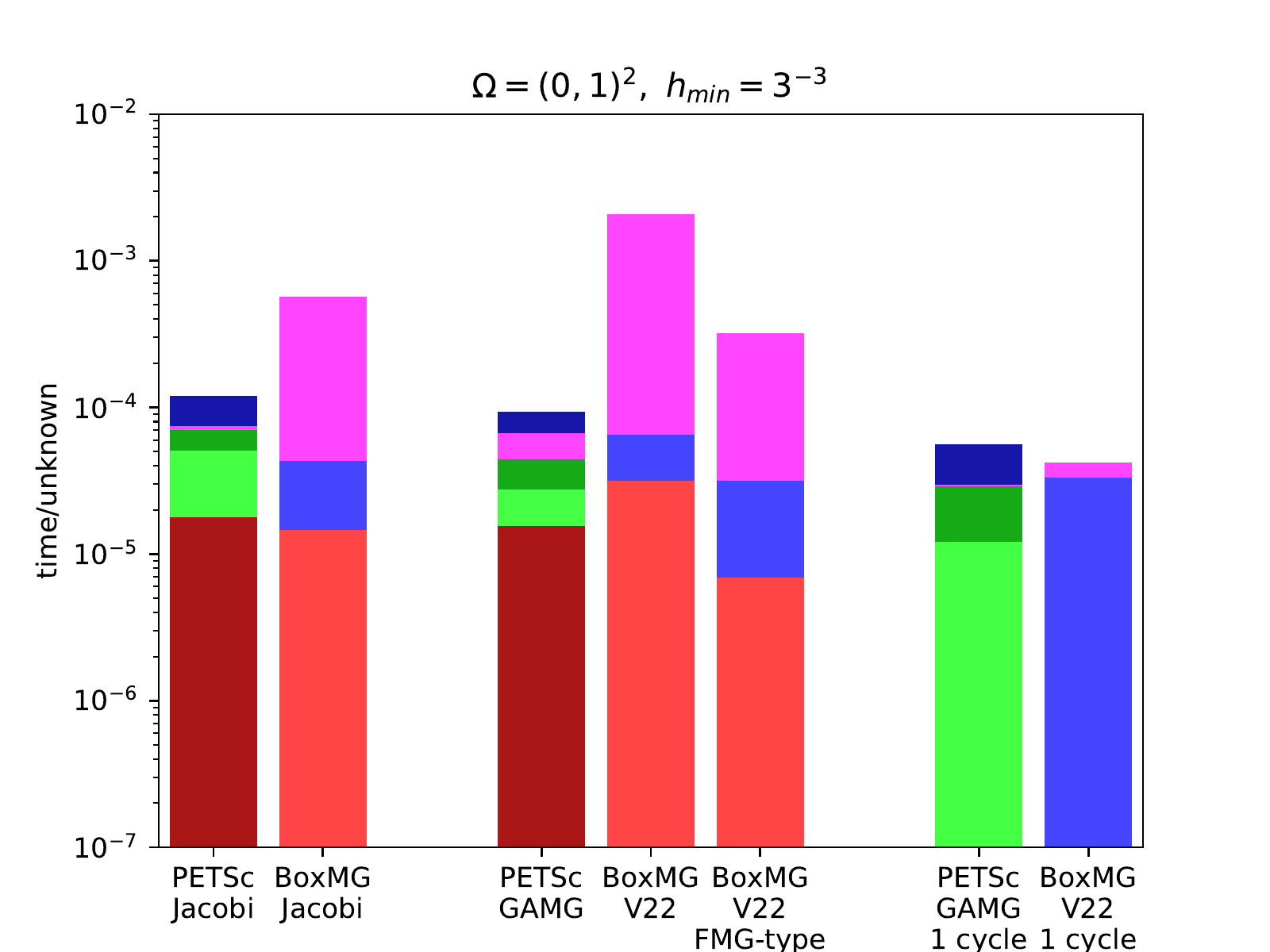}
  \includegraphics[width=0.45\textwidth]{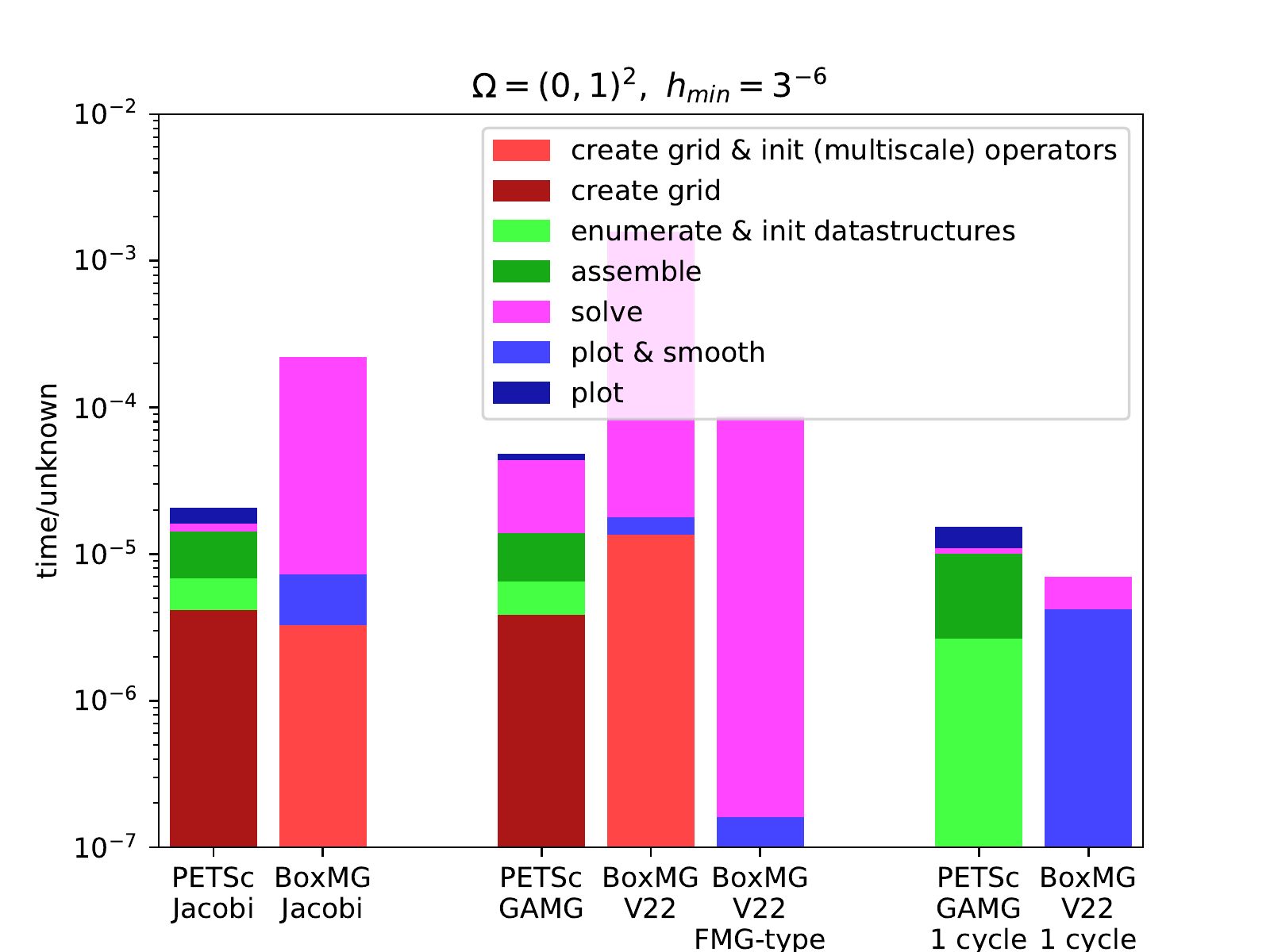}
 \end{center}
 \caption{
  \new{
   Comparison of the runtime of our BoxMG solver to PETSc's GAMG solver for the
   \texttt{checkerboard} setup in $2d$ for two different mesh sizes. We compare
   the cost of a Jacobi solver (single grid code), a full multigrid code, and
   one step of an FMG cycle.
  }
   \label{figure:results:comparison-to-petsc}
}
\end{figure}

\new{
We start our runtime studies with a brief comparison of our code's runtime to
PETSc \cite{Software:PETSc} with the aggregation-based GAMG.
Our naive realisation runs through the grid twice: In a first grid sweep, we
enumerate the unknowns and determine the matrix sparsity pattern.
In a second sweep, we assemble the matrix.
A third sweep is necessary once we have solved the equation system if we want to
plot the result that ties the PETSc solution to the grid.
If the grid changes, a complete re-assembly with two grid sweeps becomes
necessary.
As we use a space-filling curve as enumeration scheme, we obtain a sparsity
pattern alike Fig.~\ref{figure:mpi:spacetree-decomposition}.
}

\new{
We make both BoxMG and PETSc work with the same settings: 
PETSc's GAMG is configured to stop once the relative residual is reduced to
$10^{-8}$, uses a Jacobi smoother (Richardson update with Jacobi diagonal
preconditioner) with $\omega =0.8$, and an AMG connectivity threshold of $0.19$
for its coarse grid identification is applied. 
This empirically chosen value yields, for the present V22-cycle,
comparable coarse grids to BoxMG in terms of unknown counts.
}

\new{
We observe that PETSc's explicit assembly consisting of grid construction, grid
enumeration and sparsity identification plus matrix entry assembly is slightly
slower than the monolithical approach of our BoxMG implementation where
everything is done in one place (Fig.~\ref{figure:results:comparison-to-petsc}).
However, we might be able to save the PETSc enumeration phase if the grid
construction determined the sparsity pattern on-the-fly.
PETSc is, even though the coarse grids have to be determined algebraically and
it has to maintain the coarse matrices, significantly faster than our code
if we kick off from the finest grid.
}

\new{
The picture changes if we run a simulation where we start from a coarse grid,
add (applying a refinement criterion) one level after another and thus make each
solve act as prediction for the next finer solver.
The picture changes if we use an FMG-type cycle.
Our results show data for only two steps of such a cycle.
}

\begin{observation}
\new{
 Our approach is able to outperform a black-box solver such
 as PETSc's GAMG if and only if
 \begin{enumerate}
   \item the grid changes after each (or very few) solver steps, 
   \item the problem can be solved robustly with our hybrid geometric-algebraic
   ansatz, and
   \item block Jacobi/Gau\ss-Seidel smoothers are sufficient.
 \end{enumerate}
}
\end{observation}

\noindent
\new{If the problem is ill-suited for our code due to a lack of 
robustness, our approach however may act as building block applied on the finer
mesh levels while coarser problems are solved algebraically
\cite{Gmeiner:15:HHG,Lu:14:HybridMG,Sundar:12:ParallelMultigrid,Rudi:15:GordonBell}.
We reiterate that such a level $\ell_{max}$ can be determined automatically.
If frequent visualisation of the solve is required, our approach also is
promising as we can merge plotting and solution updates.
If frequent remeshing is required and, thus, we face non-negligible
assembly overhead, our approach becomes competitive. 
This is in line with other papers.
\cite{Cleary:2000:RobustAMG}, e.g., report AMG's setup time to be equivalent to
six multigrid cycles for similar algorithmic components,
while in the more sophisticated setting of
\cite{Treister:2010:MGStochastic}, the coarse grid and operator construction even seems to dominate the runtime. 
}

\new{
We conclude with some memory observations made through PETSc's
\texttt{PetscMallocGetCurrentUsage} function.
Our BoxMG with memory compression can reduce the memory footprint to close
to 25 bytes per degree of freedom.
With PETSc, the grid requires slightly more than one byte per fine grid vertex
to store the linearized spacetree---the finest grid resolution level holding degrees of
freedom dominates the memory footprint---plus one integer used as unknown index.
These five bytes per unkown are supplemented by a total of 38 kByte PETSc 
administration overhead.
The lion share of memory is allocated once we trigger the sparsity pattern 
analysis and enumeration.
We end up with PETSc alone requiring between $26.95 \cdot 3^d$ (strongly adaptive or very small 
grids) down to $15.13 \cdot 3^d$ (more regular and/or large grids) bytes per
degree of freedom.
}


\subsection{Runtime and scalability studies}
\label{section:results:parallel}

\begin{figure} 
  \begin{center}
    \includegraphics[width=0.4\textwidth]{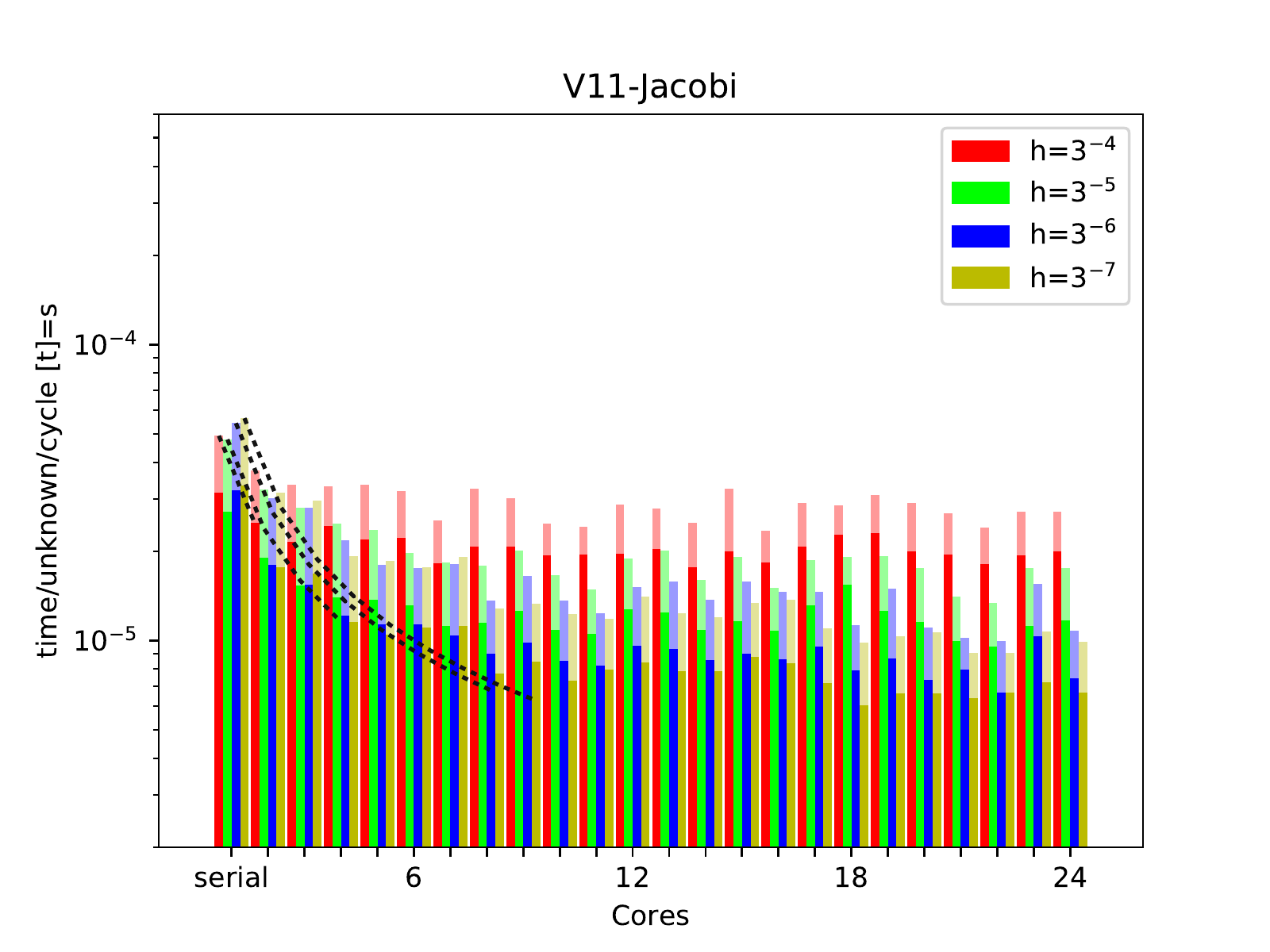}
    \includegraphics[width=0.4\textwidth]{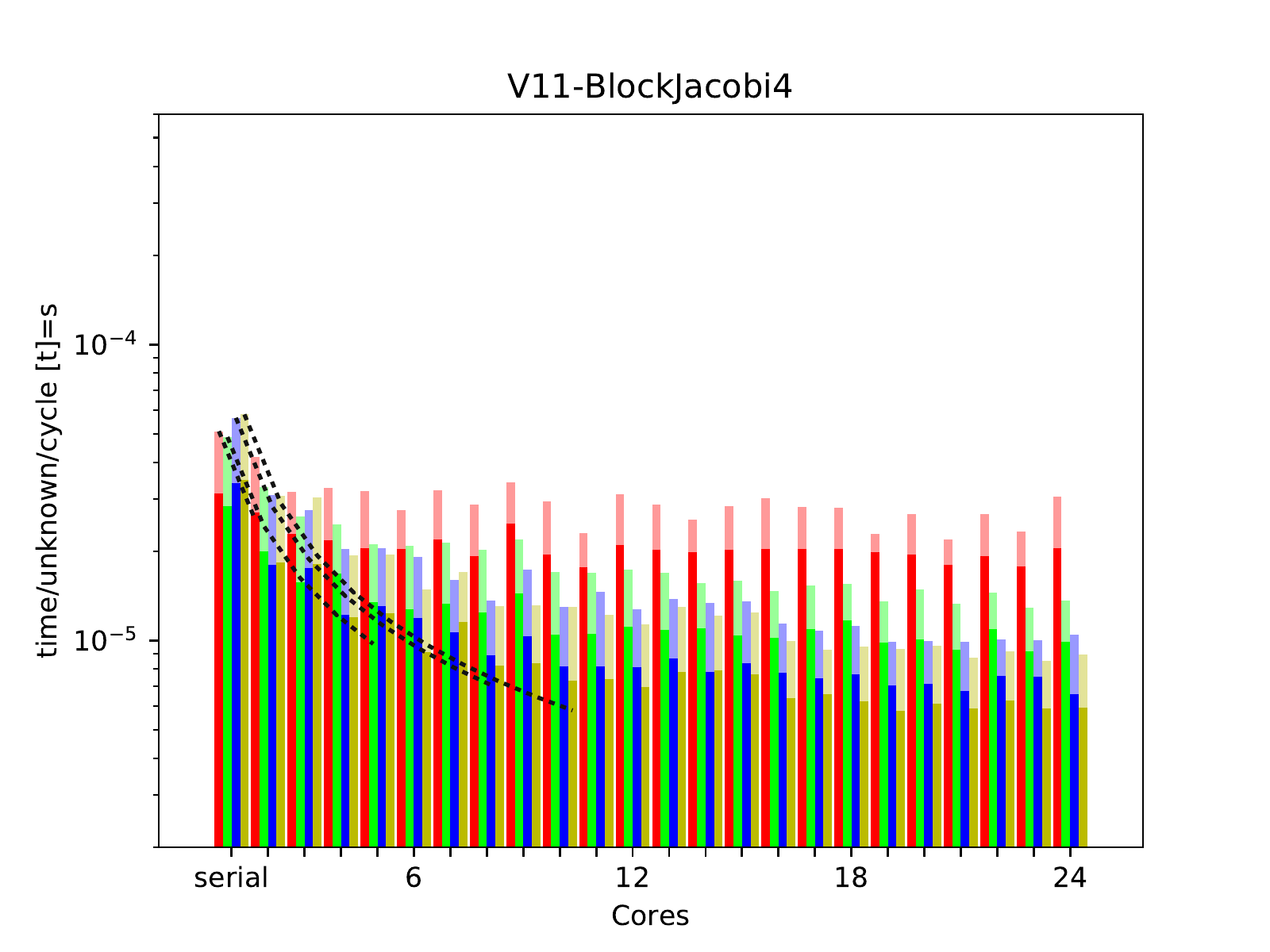}
    \\
    \includegraphics[width=0.4\textwidth]{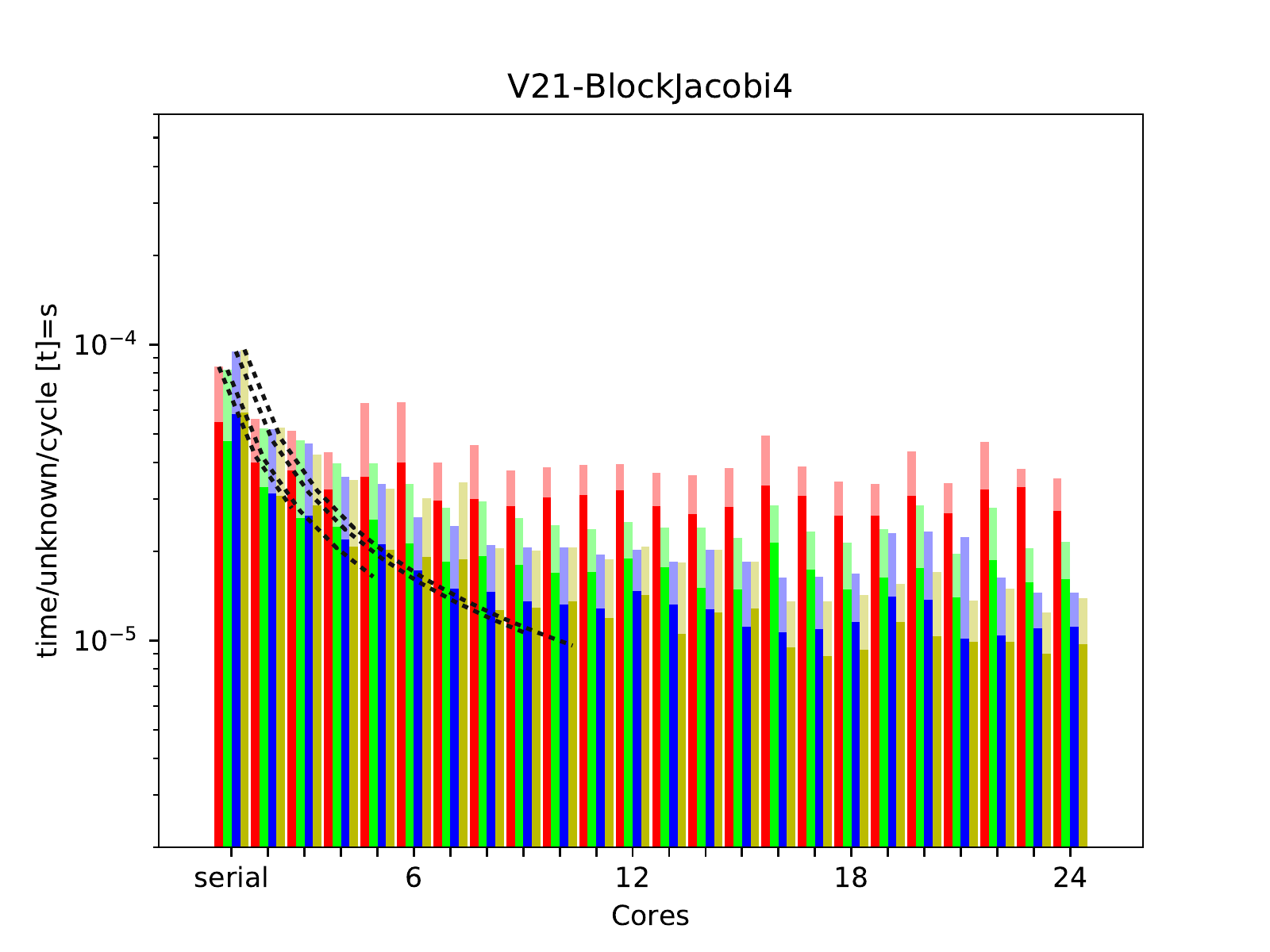}
    \includegraphics[width=0.4\textwidth]{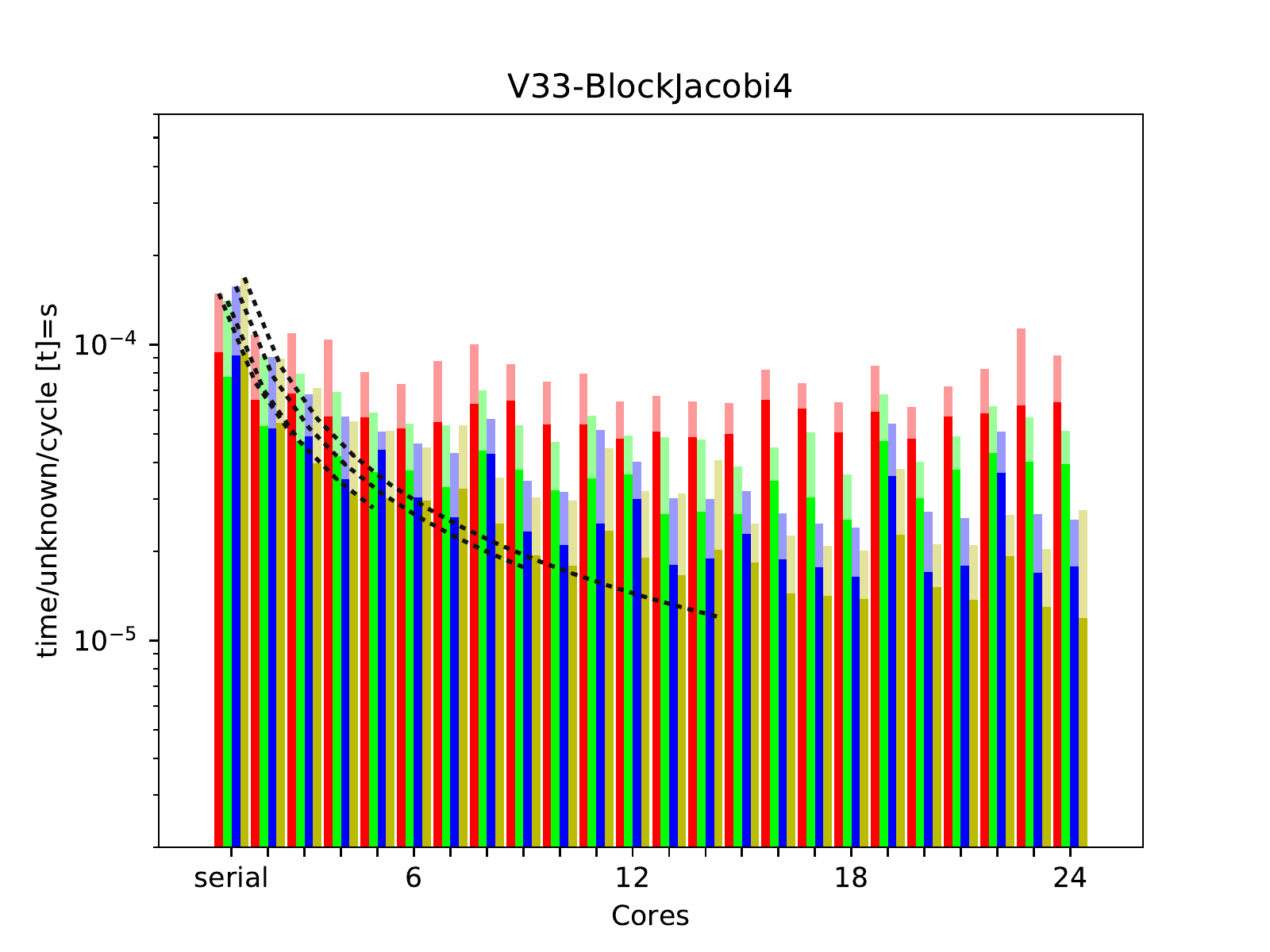}
  \end{center}
  \caption{
    Performance of one cycle of \new{the} multiplicative BoxMG solver on the
    \new{Broadwell} for \new{the  $\texttt{sin} $ benchmark} and $d=2$. 
    The dotted line illustrates linear speedup (100\% efficiency) 
    \new{truncated by the minimum runtime cost obtained}.
    Each measurement consists of two bars. The bar in the background
    (lighter, higher runtime) uses
    $\epsilon _{mf}=10^{-8}$.
    No compression is used for the measurement in the foreground.
    \label{figure:results:shared-memory}
  }
\end{figure}

We \new{wrap up} our experiments with feasibility
studies validating the parallel well-suitedness of the proposed techniques.
Our shared memory study maps the dependencies from Sect.~\ref{section:parallelisation:sharedmem} directly onto TBB tasks.
Such a \new{na\"ive tactic} is well-known to yield non-optimal performance as
the tasks exhibit small arithmetic intensity and a significant
tasking overhead.
Nevertheless, we observe speedup and
we are able to derive qualitative properties of the proposed
scheme.

%
%

We start with $d=2$ runs on the \new{Broadwell}
(Fig.~\ref{figure:results:shared-memory}) and observe that the cost of
the block smoothing is negligible.
All block data is cached and thus the flops for the small blocks are
almost for free.
\new{Despite the fact that}
the grid management overhead amortizes,
\new{we} observe that the cost per degree of freedom \new{per cycle} increases
when we increase the number of mesh levels \new{as}
additional coarse grid problems \new{are} introduced.
\new{Finally, the} runtime \new{penalty} of the on-the-fly compression on a
single core \new{is pessimistically bounded} a factor of two
\new{(cmp.~Table \ref{results:memory:cost})}.

%
%
\new{Once} we use more than one core, \new{the compression yields better
speedups} than a plain implementation.
The arithmetic intensity per grid entity is higher due to \new{the}
computation of $\hat R, \hat P$ and $\hat A$
\new{and the pressure on the memory subsystem is lower}.
Yet, this speedup improvement cannot close the gap between the code with
compression and without compression \new{completely}.
As the theoretical concurrency of the scheme increases with additional grid
levels, the speedup \new{increases with an increase of levels}.
\new{
We observe a significant strong scaling behaviour manifested by the fact that
bigger problems for large thread counts outperform smaller problems.
We further observe that the more smoothing sweeps the better the scalability.
Again, this is due to the fine grid which parallelises best.} 
\new{Block smoothers in general have comparably high arithmetic intensity and
thus improve the scalability}.
\new{At the same time, block data accesses} and inter-grid transfer operators
\new{however} reduce the concurrency level. 

\begin{observation}
  For $d=2$, block smoothers are for free in terms of computational cost. 
  Cost-per-vertex models that are linear in the number of smoothing steps and
  agnostic of the mesh size are inappropriate here.
  The cost for data compression has to be evaluated carefully for any
  application \new{though the technique seems to be promising for
  manycores}.
  Overall, the scalability is very limited; an effect due to the low order of
  the discretization inducing a low arithmetic intensity and the rigorous task
  formalism that introduces a higher administrative overhead than classic
  \texttt{for loop}-based parallelism.
\end{observation}

\noindent
\new{Overall, scalability and performance are limited. 
On all cores, Stream TRIAD \cite{McCalpin:95:Stream} yields 3,346.99 MFlops/s with a total
used bandwidth of 57,838.15 MBytes/s for the machine.
The present V33 simulation however uses only 7,061.49 MBytes/s bandwidth and
yields 317.9540 MFlops/s.
Its cache miss rate is 0.73\%.
Switching on the compression increases the compute load to 650.60 MFlops/s
but reduces the bandwidth
demands to less than 6,000 MBytes/s.
However, the cache miss rate increases to around 30\%.
}

\new{
All characteristic data highlight
the feasibility character of the study---and thus put the comparisons to PETSc
into perspective---where no
performance engineering is done, everything is
modelled with (tiny) tasks, memory access to the individual stencils held in
unoptimised hash maps are scattered, and where we study low order
discretisations tackled by multiplicative multigrid cycles coarsening up to the trivial level.
}


\begin{figure}
  \begin{center}
    \includegraphics[width=0.4\textwidth]{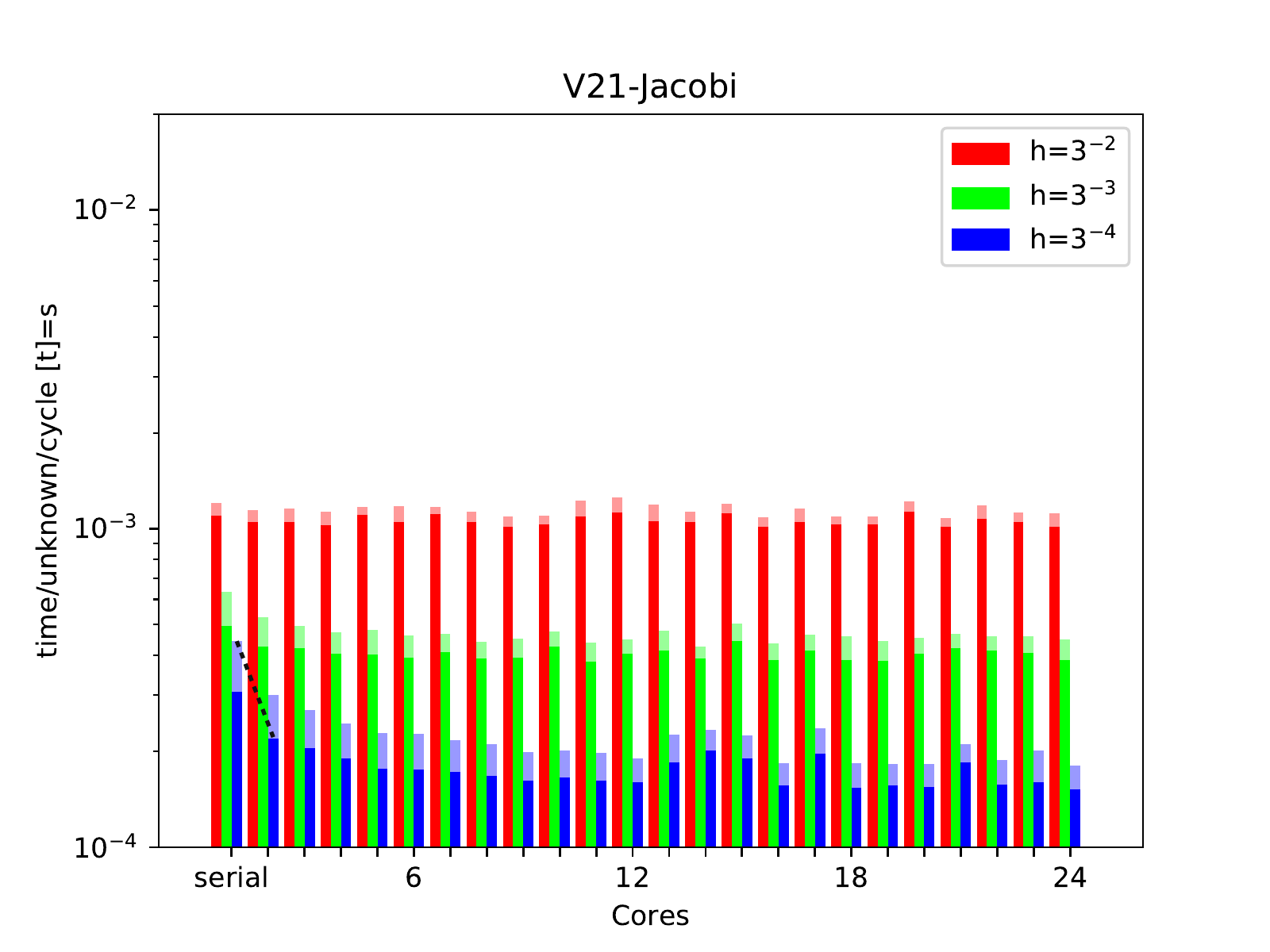}
    \includegraphics[width=0.4\textwidth]{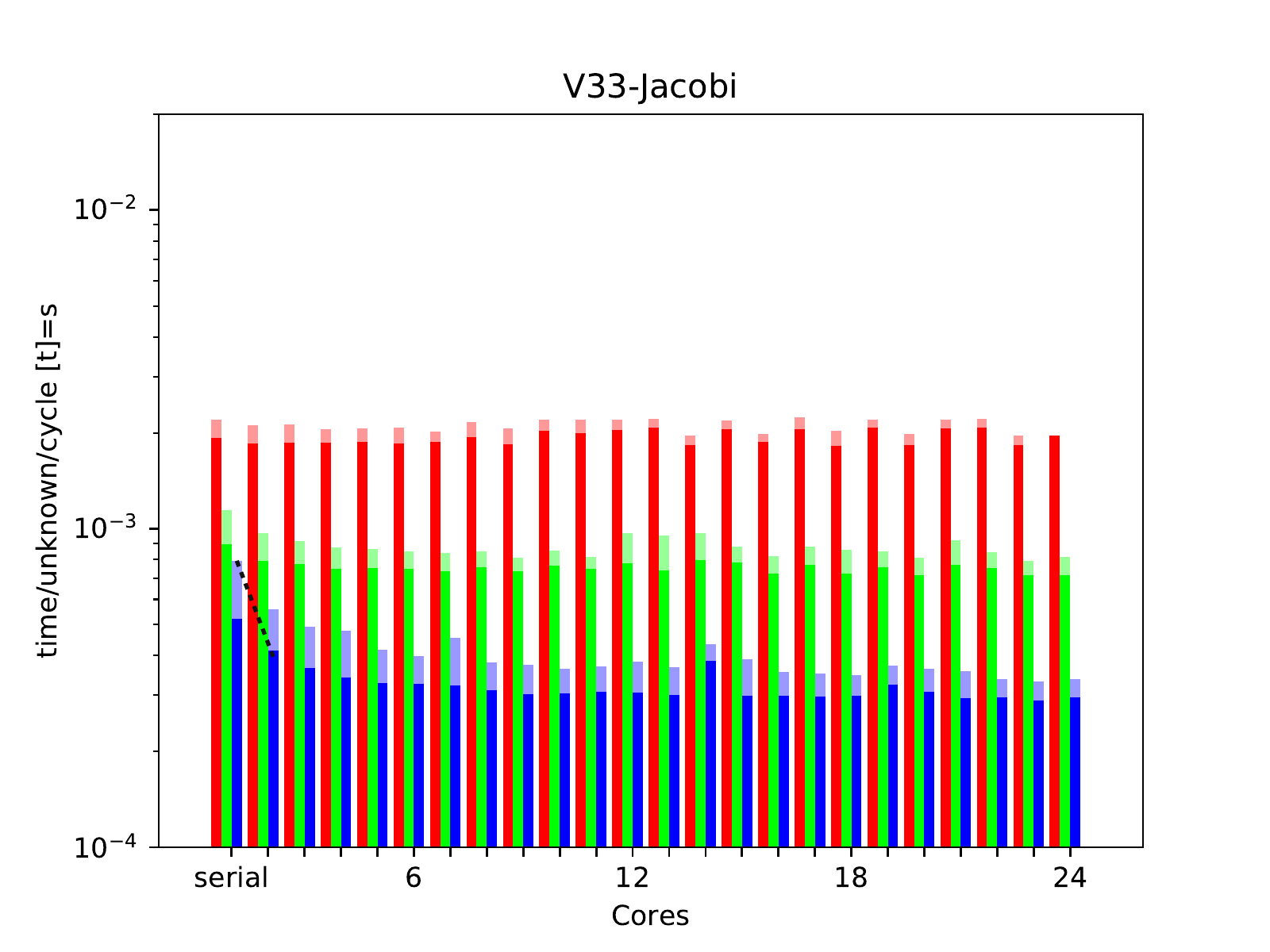}
     \\
    \includegraphics[width=0.4\textwidth]{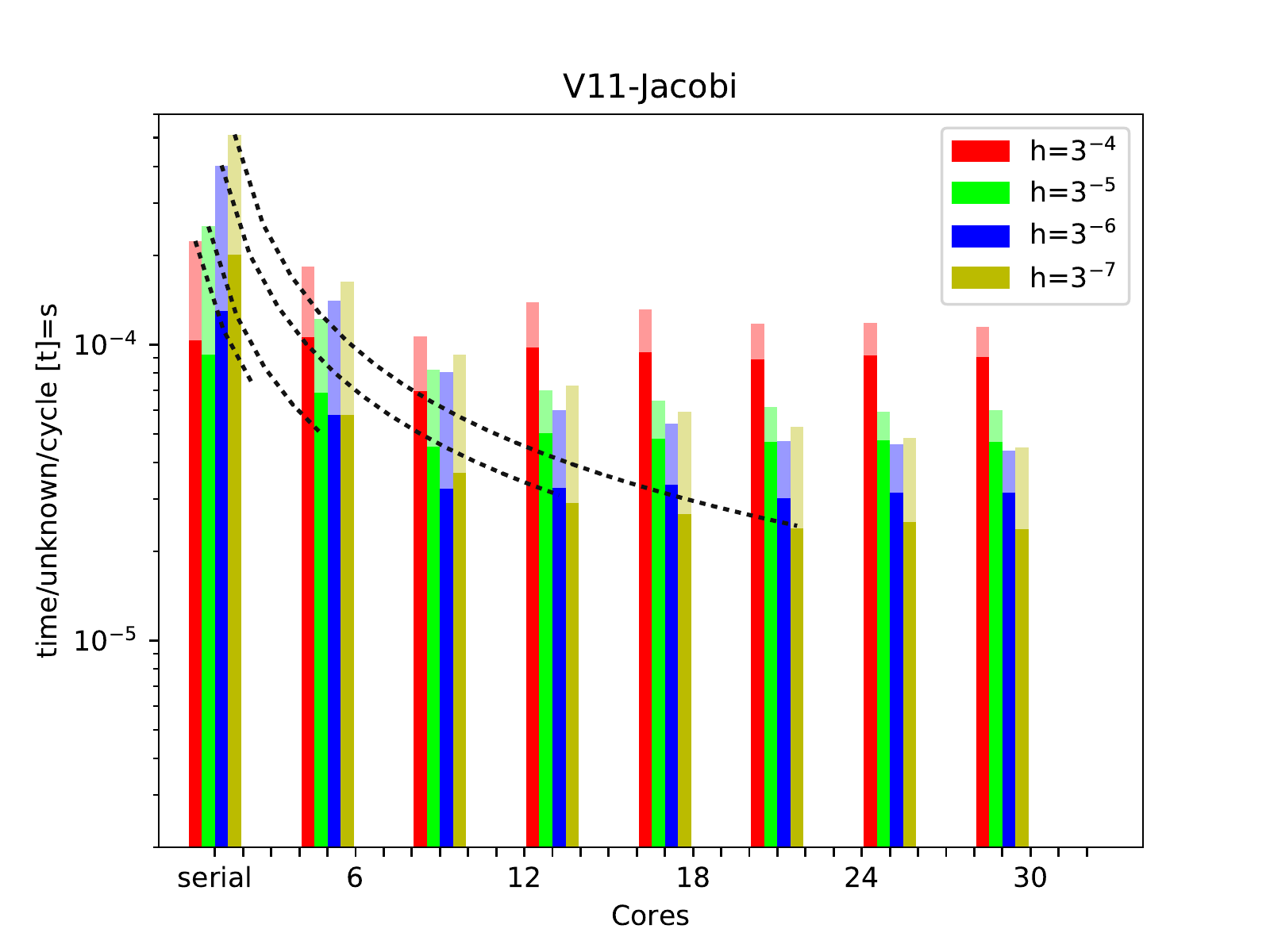}
    \includegraphics[width=0.4\textwidth]{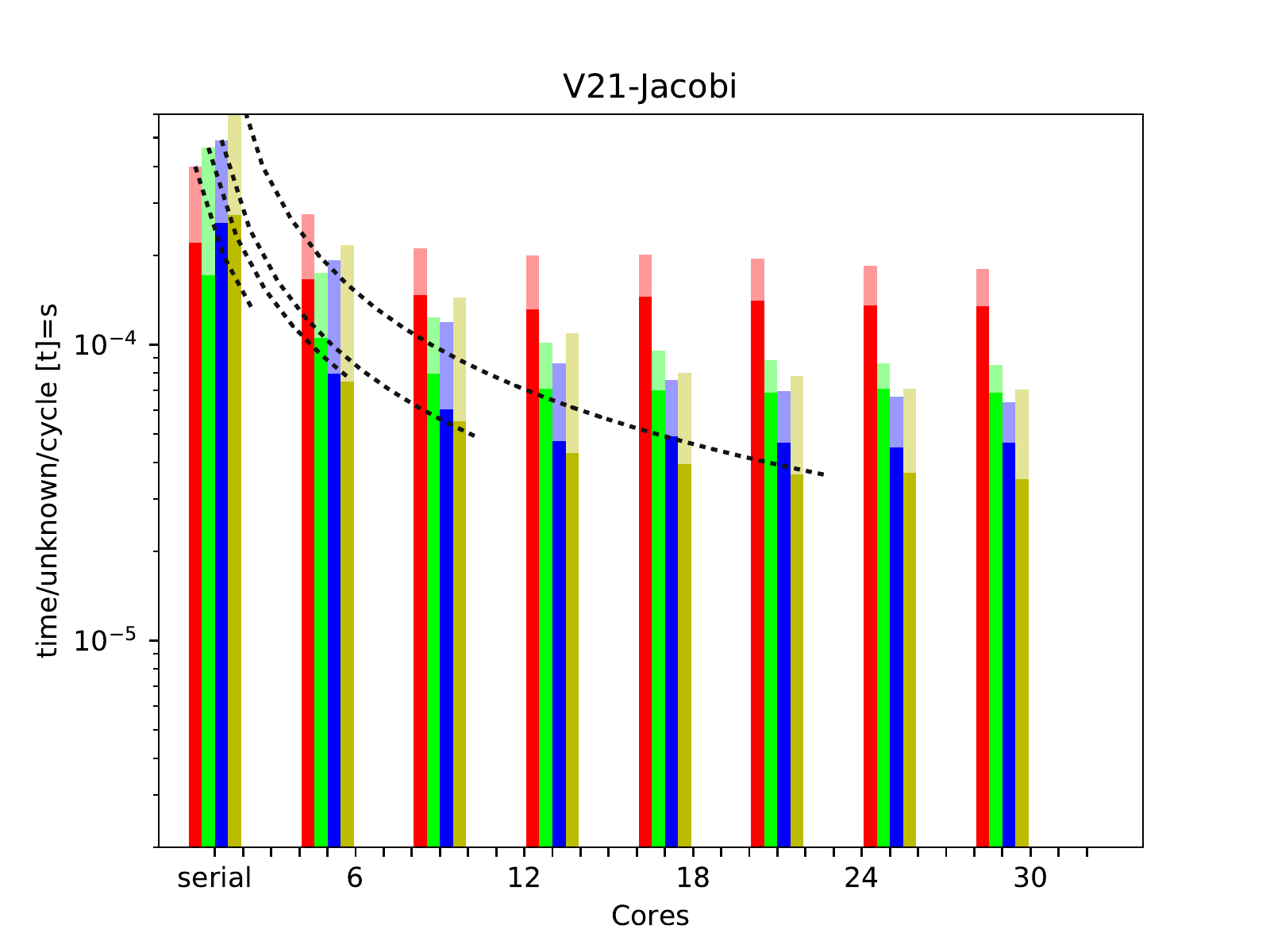}
  \end{center}
  \caption{
   \new{
    BoxMG on Broadwell for $d=3$ (top) and for the KNL for $d=2$ (bottom).
   }
    \label{figure:results:shared-memory-3d-knl}
  }
\end{figure}

%
%
We continue with $d=3$ experiments
(Fig.~\ref{figure:results:shared-memory-3d-knl}) and observe some changes
in the characteristics:
The finer the grid the smaller the cost per degree of freedom.
Administrative cost now does amortize while the \new{reduction of vertices per
coarsening by a factor of 27 is so significant that the coarse grids' runtime
behaviour has no major impact on the efficiency.} 
The qualitative scalability does not change dramatically, and we continue to see
strong scaling stagnation already for reasonably small core
counts.
\new{Stagnation in the plots however is reached sooner as the maximal grid
depths we are able to resolve on a single node are shallower than in the two-dimensional counterpart}.

\begin{observation}
  For $d=3$, the compression cost is not dominant anymore.
\end{observation}
  
\noindent
We finally rerun our experiments on the manycore architecture
\new{(Fig.~\ref{figure:results:shared-memory-3d-knl})}.
All of our statements qualitatively remain valid.
They however change quantitatively.
The relative compression cost on the manycore is lower and the code 
\new{scales to slightly higher core counts}.
\new{KNL's SNC-4 mode apparently is a well-suited hardware configuration
here which implies that we have to use at least four MPI ranks per node.} 
We summarize that the architecture seems to be even more sensitive to the
tasking overhead \new{as the
single-core performance difference does not directly relate to the difference in 
clock speed.
Yet, the chip benefits from our data compression more
significantly.
}

\begin{table}
  \tbl{
    \new{
    $d=2$ MPI experiments for various node counts (four MPI ranks
    per node, 24 cores per node) and minimal mesh sizes for V(1,1)-cycles (left)
    and V(2,1)-cycles (right).
    All speedups refer to the geometric multigrid variant running on one node
    with four ranks.
    Relative to the geometric multigrid baseline, we give the relative cost to
    compute and maintain the discretisation stencils of BoxMG and, on top of
    this, do the compression.
    }
    \label{table:mpi}
  } 
  {
   \begin{tabular}{l|l|rrr|rrr}
 \#dofs & nodes & $S_{geom}$ & stencil cost & comp.~cost 
                & $S_{geom}$ & stencil cost & comp.~cost \\
    \hline
 $5.91 \cdot 10^4$ & 
    1 & 1.00 & 1.06	& 1.53
      & 1.00 & 1.37 & 1.44
    \\
    \hline
 $5.31 \cdot 10^5$ & 
     1 & 1.14 & 1.14 & 1.57
       & 1.13 & 1.14 & 1.58
     \\
     & 
     8 & 7.74 & 2.59 & 1.16
       & 7.74 & 2.66 & 1.13
     \\
     &
     16 & 9.41 & 2.81 & 1.18
        & 9.41 & 2.84 & 1.16
     \\
     & 
     64	& 31.29 & 5.18 & 1.08
        & 33.35 & 5.22 & 1.09
     \\
    \hline
  $ 4.78 \cdot 10^6$ &          
     1 & 1.23 & 1.01 & 1.72
       & 1.22 & 1.00 & 1.68
     \\
     &
     8 & 9.20 & 1.43 & 1.44
       & 9.06 & 1.43 & 1.42
     \\
     & 
     16 & 11.40 & 1.55 & 1.44
        & 11.31 & 1.60 & 1.40
     \\
     & 
     64 & 61.33 & 3.19 & 1.21
        & 59.74 & 3.15 & 1.20
     \\
     \hline
   $ 4.31 \cdot 10^7$ & 
     8 & 9.78 & 1.06 & 1.78
       & 9.69 & 1.07 & 1.74
     \\
     & 
     16 & 12.11 & 1.25 & 1.76
        & 11.84 & 1.28 & 1.47
     \\
     & 
     64	& 54.91 & 1.38 & 1.65
        & 52.39 & 1.26 & 1.45
     \\
     \hline
   $ 3.87 \cdot 10^{8}$ &
     16 & 12.92 & 1.18 & 1.89 
        & 12.83 & 1.20 & 1.83 
     \\
     & 
     64 & 59.30 & 1.00 & 2.21
        & 62.22 & 1.27 & 1.53
\end{tabular}
    
  }
\end{table}

For our subsequent MPI experiments, we restrict to $(\mu
_{pre},\mu_{post})=(1,1)$ and $(\mu _{pre},\mu_{post})=(2,1)$.
Both are challenging choices in terms of scalability as the arithmetic work is
small compared to the inter-grid operator evaluations and grid level changes.
For the $(\mu _{pre},\mu_{post})=(2,1)$, we explicitly exchange all residuals
plus stencil contributions per smoothing step, i.e.~we anticipate the data flow
from a non-linear problem \new{and ignore the fact that the first smoothing 
step does not have to exchange partial stencils}.
For the decomposition, we naively apply graph partitioning on the start grid
chosen reasonably fine such that it can accommodate all MPI ranks.
Any dynamic load balancing is switched off.
We have validated that the partitioner yields perfectly balanced
subdomains for regular grids.
All cost per degree of freedom and, thus, all speedups
are normalized to the run with the smallest problem size.
\new{Following our shared memory results, we deploy four MPI ranks per node.}

%
%
We obtain a reasonable scalability $S_{geom}$ of the geometric baseline code
for both cycles once the problem sizes are sufficiently big (Table~\ref{table:mpi}):
When we increase the problem size, all grid administration overhead gets
amortized.
This effect materializes in classic weak speedup, and it also materializes in
speedups for serial runs \new{bigger than the minimum mesh size}.

%
\new{The merger of
algebraic multigrid's stencil storage into the geometric code increases
the serial runtime slightly.
Holding the stencil within the grid also reduces the scalability for the
majority of the setups, i.e.~the cost grows if we use multiple nodes.
}
More data has to be piped through the communication network. 
The bandwidth demands increase.
\new{This notably is problematic once we run into a strong scaling regime. 
For very large problems with higher node counts, the stencil administration and
communication penalty is not that significant.
}
%
Stencil compression increases the \new{runtime further, but this relative
increase decreases, for the majority of the setups, with growing node counts}.
\new{The scheme} releases stress from the communication network.


\section{Conclusion and outlook}
\label{sections:conclusion}


This paper proposes a, to the best of our knowledge, new combination of multigrid 
techniques and novel implementation concepts for quasi-matrix-free geometric-algebraic multigrid on 
dynamically adaptive grids.
First, we apply the BoxMG principle to additive and BPX-type solvers together
with HTMG and, thus, are able to support vertical integration and dynamically
adaptive grids without any constraints on the frequency of the grid refinement or
transition of refinement regions.
Second, we discuss an on-the-fly stencil compression that brings together the
robustness of BoxMG with the memory modesty of geometric rediscretization.
There are efficient matrix storage schemes for dynamically adaptive formats
\cite{King:16:DynamicCSR}, but our approach goes beyond that as it analyzes the
operators themselves.
It also reduces the amount of data exchanged between multiple ranks.
Finally, we sketch, as third methodological contribution, the impact of the
proposed algorithms on parallel programming.
The resulting family of solvers is a hybrid between algebraic and geometric
multigrid and a hybrid between matrix-free and stencil-holding techniques.
A third flavour of hybrid---purely algebraic solvers on coarse grids or
forests supplemented by geometric grid hierarchies on finer levels---would fit
to the proposed concepts.

While the realization idioms are elegant and the concept of mirroring 
the arising BoxMG equation systems to a reference configuration makes the higher dimensional implementation much less
tedious compared to setting up the equation systems straightforwardly, our experiments reveal that the
convergence speed for convection-dominated problems deserves additional
attention.
Next steps are deriving well-suited 
estimators that autonomously identify a good $\ell _{max}$ more
elegantly and implementing more sophisticated Petrov-Galerkin inter-grid
transfer operators.
The most important multigrid shortcoming of the present work is the
restriction to Jacobi and block Jacobi smoothers.
This restriction results from a single-touch single-traversal doctrine in
combination with the element-wise tree traversal.
As such, the present studies have academic character, and it is important in the
future to weaken the single-sweep paradigm
\new{if it renders it possible} to realize stronger
smoothers.
\new{Candidates for suitable smoothers are 2-sweep Krylov schemes \cite{Bader:08:MemorySierpinski} or red-black Gauss-Seidel with pipelining which
combines multiple sweeps \cite{Ghysels:13:ModelMG}.}
While \new{giving up on single touch harms implementational elegance},
it might even turn out to be favourable from a parallelization point of view to
run over the grid multiple times as long as the rank-local work increases faster than the exchanged data cardinality.

Finally, we emphasize a 
solver property that deserves particular attention:
the studied class of low order discretizations yields compact stencils with
relatively low arithmetic intensity.
Such stencil codes can, in the context of iterative solvers, significantly
benefit from careful tuning such as diamond tiling.
However, most tunings require invariant stencils in order to perform
\cite{Malas:2016:IntraTileStencils}.
Our work targets problems where stencil entries are not constant.
At the same time, it is able to compress data automatically in areas where
stencils are known a priori.
It therefore seems to be promising to inject state-of-the-art stencil techniques
for those regions where the compression pays off and to preserve the present
approach's robustness everywhere else.

With the obtained solver robustness our approach widens
the class of (sub)systems that can be solved with a spacetree-based multigrid
solver significantly, while it preserves the structuredness and low memory
footprint properties of geometric multigrid solvers.
The software concept thus can become an enabler to solve challenging
PDE problems on the grand scale where structuredness is important to optimize
and parallelize and memory (per core) is a precious resource.
To achieve this, state-of-the-art techniques from computer science and numerical
linear algebra had to be merged into one realization concept.

\ifthenelse{\boolean{arxiv}}{
  \section*{Acknowledgements}
  This work made use of the facilities of the Hamilton HPC Service of Durham
University.
Furthermore, we particularly have been benefitting from the support of the RSC
Group who granted us early access to their KNL machines.
Thanks are due to Hans-Joachim Bungartz for supervising and mentoring the underlying PhD theses.
It was Irad Yavneh's suggestion to use the idea of BoxMG in this context,
and first steps into this direction where made in
\cite{yavnehbendig12nonnsymbb}. Marion wants to thank him for hosting her
research stay at Technion which laid the foundations for the present work.
Operator and data compression is a workpackage in the ExaHyPE project and Tobias
thus appreciates support received from the European Union’s Horizon 2020
research and innovation programme under grant agreement No 671698 (ExaHyPE).
All underlying software is open source \cite{Software:Peano}.

  \bibliographystyle{siam}
  \bibliography{paper}
}{}

\ifthenelse{\boolean{sisc}}{
  \bibliographystyle{siamplain}
  { \footnotesize
  \bibliography{paper}
  }
}{}

\ifthenelse{\boolean{arxiv}}{
  \appendix
  \section{Spacetree construction, properties and traversal}

\label{sections:appendix:spacetrees}

All cells of a spacetree $\mathcal{T}$ with level $\ell $ span a
grid $\Omega _{h,\ell}$.
Such an $\Omega _{h,\ell}$ can be ragged, but all cells have exactly
the same size.
The union of all $\Omega _{h,\ell}$ yields an adaptive Cartesian mesh $\Omega
_{h}$.
The recursive construction scheme ensures that the grids 
$
  \Omega _{h,0} \subseteq \Omega _{h,1} \subseteq \Omega _{h,2} \subseteq \ldots
  \subseteq \Omega _{h}
$
are embedded into each other. 
These grids carry our shape function weights, i.e.~we solve
\begin{eqnarray*}
   A_\ell (u_\ell+e_\ell) & = & b_\ell + A_\ell u_\ell 
     =   Rr_{\ell+1} + A_\ell I u_{\ell+1} \\
   & = & R\left( b_{\ell+1}-A_{\ell+1}u_{\ell+1} \right) + RA_{\ell + 1} P I
   u_{\ell+1} 
     =  R\left( b_{\ell+1}-A_{\ell+1}(id - P I) 
   u_{\ell+1}\right) \\
   & =: & R\left( b_{\ell+1}-A_{\ell+1}\hat u_{\ell+1} \right) 
     =: R \hat r_{\ell+1}
\end{eqnarray*}
on them.

Various efficient strategies for storing and handling spacetrees,
i.e., traversing and providing adjacency information,  are
known---notably in combination with SFCs:
storage within (hash) maps \cite{Griebel:99:SFCAndMultigrid},
serialization/linearization of the whole tree \cite{Sundar:08:BalancedOctrees},
which allows us in our code \cite{Software:Peano} to encode a tree traversal
into a push-back automaton holding all adjacency data \cite{Weinzierl:11:Peano}, on-the-fly computation of
neighbours via Morton codes \cite{Burstedde:11:p4est}, and so forth.
All enlisted variants do not store adjacency explicitly and thus avoid memory
overhead.
In the present paper, we do not restrict ourselves to a
particular storage scheme but rely on the generic concept of {\em multiscale
element-wise traversals}.

\vspace{-0.1cm}
\begin{observation}
The classic depth-first tree traversal yields a
multiscale \linebreak element-wise traversal of the spacetree's multiscale grid. 
\end{observation} 

\vspace{-0.1cm}
\noindent
This observation \cite{Mehl:06:MG,Weinzierl:11:Peano} is important as
depth-first interweaves the traversal of multiple scales.
Furthermore, a cell is ``left'' if and only if all of its children have been
processed.
Such a vertical integration \cite{Adams:15:Segmental} ensures high temporal and spatial data access locality:
The probability that a vertex manipulated by one cell is required soon after
again by a neighbouring cell is high\new{---an effect amplified by the usage of
space-filling curves to order the children of any refined tree node}.
We obtain excellent cache behaviour \cite{Mehl:06:MG}.

\begin{figure}[htb]
 \begin{center}
   \includegraphics[width=0.6\textwidth]{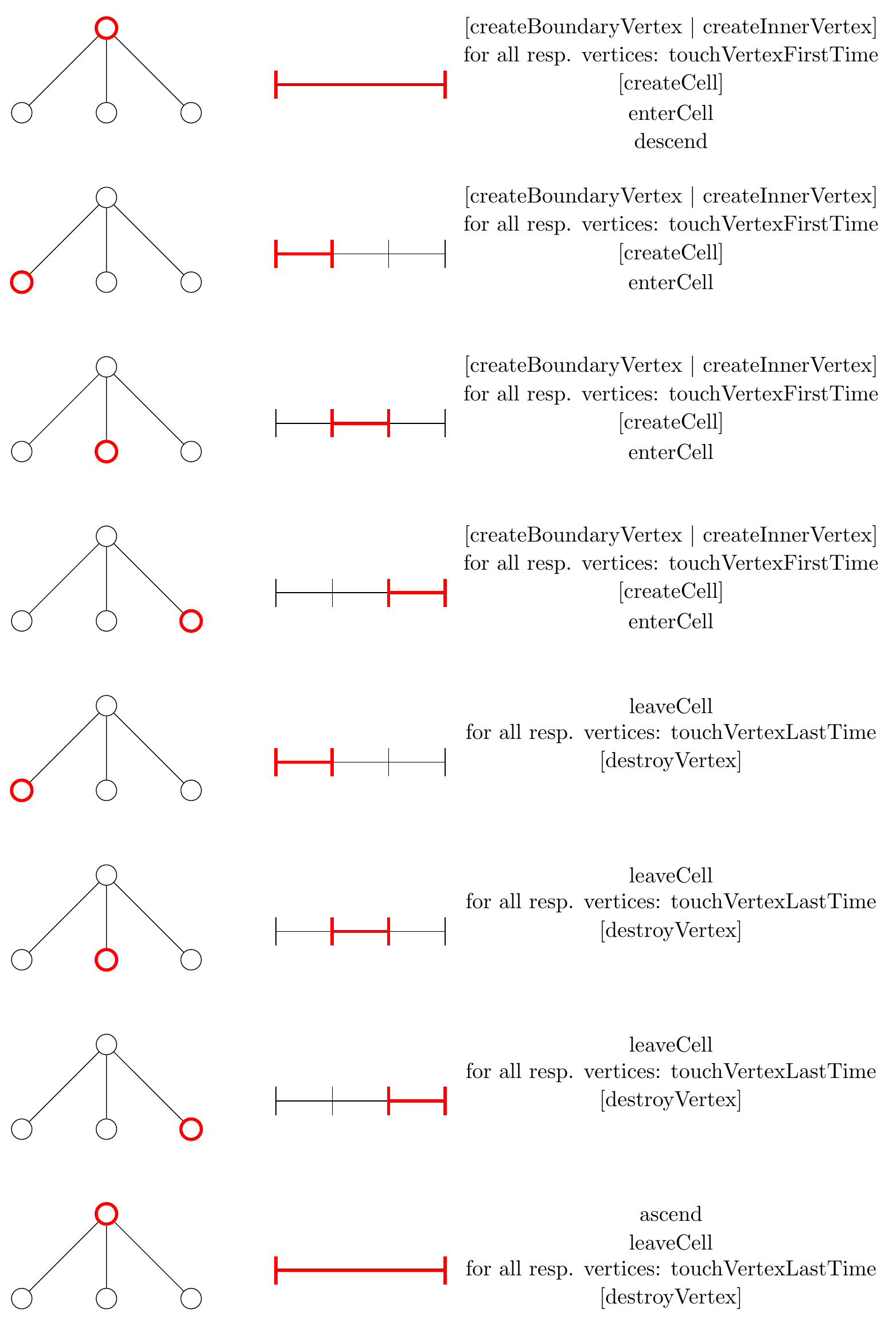}
   \caption{\new{A simple 1D example for the order of (sequential) events during
the spacetree traversal (from \cite{mweinzierl13diss}). The left column shows the tree, the middle column the grid, and the right column the respective events that are called.
The current active cell/tree node is highlighted in red.
At the beginning of the traversal \textit{beginIteration} and at the end \textit{endIteration} is called.
Parallel events and events concerning hanging vertices are omitted. Square
brackets stand for events which only occur during the setup/cleanup phase or in
an adaptive setting.}}
  \label{figure:events}
 \end{center}
\end{figure}

The resulting element-wise traversal facilitates solely point Jacobi smoothers
as outlined in Sect.~3.2.
Point Jacobi is a poor choice for many non-trivial parameter combinations in
(\ref{equation:pde}).
To facilitate more powerful smoothers without giving up data locality or single
touch \new{is our motivation to} generalize the tree traversal by a
$descend$ event \new{(Figure \ref{figure:events})}.
In a depth-first traversal code, such an operation makes a recursive step down
within the tree and loads all children of a node before it continues
recursively---a one-level recursion unrolling
\cite{Eckhardt:10:Blocking}.
Though this technique allows us to realize inter-grid transfer
operators, we stick to vertex-wise transfer operator
realizations. Yet, we use $descend$ to implement block smoothers.  

  \section{Review of additive multigrid realisation ideas}
\label{appendix:geoaddMG}

\new{
Geometric additive multigrid with rediscretization and one Jacobi
smoothing step per level fits into the multiscale element-wise traversal
once we shift the standard multigrid cycle by half a grid traversal (Algorithm
\ref{algorithm:geometric::additive}) and introduce a helper variable $d$ for the
correction terms:
we switch from a fine-to-coarse to a
coarse-to-fine grid level enumeration plus backtracking of the call stack.
All required matvec entries can be determined on-the-fly per cell.
This holds for $diag(A_\ell)$ extracting diagonal elements from $A_\ell$ too.
Let $\omega $ be a generic smoothing parameter.
In the present paper, we either use a constant $\omega $ or we damp $\omega$
exponentially, i.e., we use $\omega $ on the finest level, $\omega ^2$ on the
first correction level, $\omega ^3$ on the next level, and so forth.
Within the spacetree paradigm, the smoothing factor is decreased with the number of coinciding vertices on finer
levels: the coarser the level, the smaller the smoother impact.
Such a vertex count can be realized during the bottom-up steps
\cite{Reps:15:Helmholtz}.
}

\new{
Algorithm
\ref{algorithm:geometric::additive} works out-of-the-box for 
adaptive grids if we make the operator $R$ affect only refined vertices and 
set the nodal value $u$ in any hanging vertex to the $d$-linear
interpolant from the coarser levels.
Any grid region can be refined dynamically.
Textbook multigrid requires higher order interpolation for newly added
levels/vertices \cite{trottenberg01multigrid}.
We obtain reasonable convergence speed if we assign newly created vertices the
linear interpolant and then immediately apply one undamped Jacobi step.
Removing vertices works without any additional effort.
Textbook multigrid typically demands an exact coarse-grid solve on $\ell
_{max}$. 
We either skip exact coarse-grid solves---in this case, we run one
Jacobi step on the coarsest level and have to study the deterioration of the
convergence speed---or we run Jacobi sweeps on the coarsest grid until the
residual there underruns $10^{-12}$.
Better iterative or direct solvers are more reasonable in many cases.
See \cite{Bader:08:MemorySierpinski} for a  
(preconditioned) CG that fits to the present matrix-free
paradigma.
}

\new{
Lines 2--4 of Algorithm \ref{algorithm:geometric::additive} translate directly into activities that we perform whenever a
vertex is read for the first time during a grid traversal.
The residual accumulations in line 6 are performed as element-wise
operations. 
All remaining updates must be realized when a vertex is
touched for the last time during a multiscale grid traversal.
The recursion and the cell-wise updates are concurrent.
Their evaluation can be permuted.
A depth-first spacetree traversal for example intermixes cell
operations with vertex updates on finer levels.
Coarse-grid residual evaluations then work with inconsistent nodal
approximations $u$.
Since the fine-grid work updates coarse-grid values
during the computation (last branch), some coarser operator
accumulations rely on outdated unknown values from the previous traversal that
are then updated while the residual is computed. 
This data inconsistency can be eliminated by an
additional helper variable \cite{Reps:15:Helmholtz}.
}

\new{
Starting from the additive multigrid, we can write down a BPX variant
with a single-touch policy (Algorithm
\ref{algorithm:geometric::bpx}) if we introduce an additional helper variable
$i$ carrying an injection of the fine-grid updates from c-points.
A {\em c-point} is a grid point that also exists on the next coarser level.
Any refined vertex coincides with at least one c-point.
Different to additive multigrid, BPX automatically keeps all levels
consistent, as each unknown always is updated only on the coarsest grid
where a vertex exists.
Spatially coinciding vertices on finer levels hold copies
of coarse grid weights.
As a result, $\omega $ is level-independent. 
The present realization is introduced in \cite{Reps:15:Helmholtz}.
For both additive multigrid and BPX, the first grid sweep realizes only
a Jacobi smoother on the fine grid.
From the second traversal on, each grid sweep realizes one multilevel
update and anticipates operations from the follow-up iteration.
}

  \section{BoxMG realisation as one linear equation system solve}
\label{appendix:boxmgtensors}

The following section discusses the operator $P$ construction.
$R$ is constructed accordingly.
We may rewrite BoxMG's per patch operations into a matrix depending on a vector
$s \in \mathbb{R}^{4^d \cdot 3^d}$ that is applied to a vector $p \in
\mathbb{R}^{2^d \cdot 3^d}$.
Here, $s$ is a collection of the stencil entries of all vertices
within a $4^d$ patch.
The stencil entries are enumerated lexicographically.
The vector $p$ contains the prolongation operator's stencil entries of the $2^d$ affected coarse
grid vertices of one patch.
The stencil collapsing (cmp.~Sect.~\ref{sections:boxmg}) ensures that the
affected section of the $P$ vectors has only the cardinality $2^d \cdot 3^d$.

\paragraph{\new{$d=1$}}
We start to illustrate our notation and implementational techniques at hands of
a 1d setup where obviously BoxMG's collapse idea does not kick in.
All enumerations start with 0.
We study one patch.
Within a patch, we first study the impact of the patch on the prolongation
stencil of the left coarse grid vertex $p=[p_0, p_1, p_2, p_3, p_4]$.
BoxMG will determine the entries $p=[p_2, p_3, p_4]$, i.e.~only a subset of this
stencil is of interest.
We reiterate that it yields ones in the centre of the $P$ stencils,
i.e.~$p_{2}=1$, so we could reduce the $p$ input further. 
However, we enforce this central value manually and thus stick to three output
values.

BoxMG's per-patch method invocation accepts an input vector $s \in \mathbb{R}^{4 \cdot
3}$ holding all the stencils of the patch's vertices.
Different to the original BoxMG paper and the present manuscript's construction
of BoxMG step by step, our implementation computes all $p$ entries from one
equation system
\[
  C(s) p = f, \qquad p,f \in \mathbb{R}^{3^d}
\]
with a right-hand side 
\[
  f_i = \left\{
    \begin{array}{cl}
      1 & \mbox{for $i=1$ and} \\
      0 & \mbox{otherwise}
    \end{array}
  \right..
\]

\noindent
The artificial entry $1$ will ensure $p_2=1$.
\new{In classic multigrid manuscripts, the formula for prolongation entries
$p$ depends on the system matrix $A$ on the next finer level.
For our BoxMG implementation,} we use $C \cdot A$ as operator \new{of such a
formula} where the entries of the matrix $A$ are determined by the stencils $s$.
It is thus convenient to write $C(s)$ here.
We reiterate that $C$ depends linearly on $s$ and thus is a tensor of third
order.
However, we here write down the tensor as 2nd order tensor, i.e.~as matrix,
where the inner multiplication is explicitly expressed.
1d BoxMG then reads \new{as}

\begin{equation}
\underbrace{
\left(
\begin{array}{ccc}
      1       &       0 &       0  \\
      s_{1,0} & s_{1,1} & s_{1,2} \\
            0 & s_{2,0} & s_{2,1}  
\end{array}
\right)
}
_{C(s)}
\left(
\begin{array}{c}
p_{2} \\
p_{3} \\
p_{4} 
\end{array}
\right)
=
\left(
\begin{array}{c}
1 \\ 0 \\ 0
\end{array}
\right).
\label{equation:boxmg:1d-system-matrix}
\end{equation}

\noindent
Within the patch, we have to determine $P$
entries for the right coarse grid vertex, too.
Instead of setting up a separate equation system for the second coarse grid
vertex---a process that becomes tedious for higher dimensions---we mirror all
coarse grid computation problems to a reference configuration where the coarse
grid operator affected is tied to the left coarse grid vertex.
Before we trigger the equation system solve, 
we replace the entries in (\ref{equation:boxmg:1d-system-matrix})
with $s_{i,j} \mapsto s_{3-i,2-j} ,\ i \in \{0,1,2,3\} , j \in \{0,1,2\} $,
i.e.~we mirror them within the patch along the x-axis.
The indices $i\in \{0,3\}$ could be omitted here, but we use them below.
After the solve of (\ref{equation:boxmg:1d-system-matrix}), we mirror the
resulting $p$ entries back onto the $p$ entries of the right coarse grid vertex:
$p_i \mapsto p_{4-i}$.

For plain BoxMG, stencils in the corners of patches ($s_{0,x}$ and $s_{3,x}$) do 
not influence $P$.
They could be omitted.
We however keep them in our code as it allows us to work with a vector of $4^d$
stencil entries instead of $4^d-2^d$ which simplifies the implementation for
$d \geq 2$.

\begin{figure}[htb]
 \begin{center}
  \includegraphics[width=0.65\textwidth]{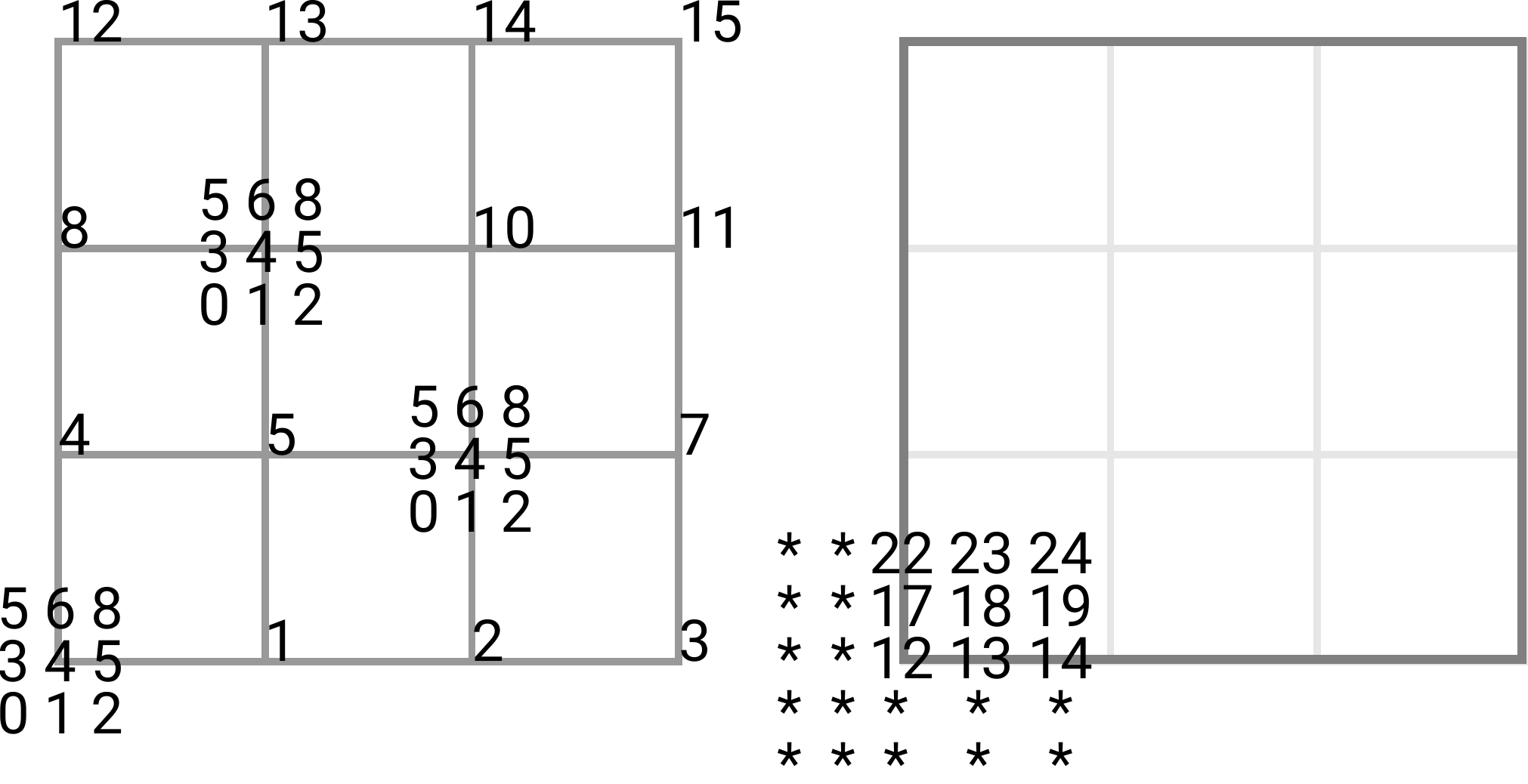}
 \end{center}
 \caption{
  \new{
  Right: Our BoxMG code computes nine entries of the bottom left vertex for
  $d=2$ within $descend$.
  Left: All fine grid vertices are enumerated lexicographically starting with 0
  at the bottom left (first index subscript). Each vertex carries a $3^2$
  stencil enumerated lexicographically, too.
  }
 }
 \label{figure:boxmg:enumeration}
\end{figure}

\paragraph{$d=2$}
For two dimensions, \new{our} per-patch \new{notation} can be written down
as
\[
 \left(
  \begin{array}{ccc}
   id & 0 & 0 \\
   \tilde A_{\gamma c} & \tilde A_{\gamma \gamma } & 0 \\
   A_{\iota c} & A_{\iota  \gamma } & A_{\iota \iota } 
  \end{array}
 \right) 
 \left(
  \begin{array}{c}
   P_c \\
   P_\gamma \\
   P_\iota
  \end{array}
 \right)
 u_{\ell -1} =: 
 C \cdot 
 \left(
  \begin{array}{ccc}
   A_{cc} & A_{c\gamma } & A_{c\iota } \\
   A_{\gamma c} & A_{\gamma \gamma } & A_{\gamma \iota } \\
   A_{\iota c} & A_{\iota  \gamma } & A_{\iota \iota } 
  \end{array}
 \right) 
 \left(
  \begin{array}{c}
   P_c \\
   P_\gamma \\
   P_\iota
  \end{array}
 \right)
 u_{\ell -1} = 
 \left(
  \begin{array}{c}
   1 \\
   0 \\
   0
  \end{array}
 \right),
\]
where $C$ is the collapsing operator acting on the PDE discretization $A$, and
$id \in \mathbb{R}^{2^d \times 2^d}$ is the identity.
\new{
We study only the bottom left vertex.
Our enumeration is bottom-up, left-to-right for both the vertices and the
stencil entries (Fig.~\ref{figure:boxmg:enumeration}).
So whenever we write $s_{i,j}$, $i$ is the vertex number and $j$ is the entry of
the vertex's stencil.
For the coarse grid vertex, we can only determine the nine entries 
$(
\underbrace{p_{12}}_{P_c},
\underbrace{p_{13},p_{14},p_{17}}_{{P_{\gamma_0}},{P_{\gamma_1}},{P_{\gamma_2}}},
\underbrace{p_{18},p_{19} }_{P_{\iota_0},P_{\iota_1}},
\underbrace{p_{22}}_{{P_\gamma}_3},
\underbrace{p_{23},p_{24} }_{P_{\iota _2},P_{\iota _3}})$. 
}

Since $C$ collapses all stencils along patch boundaries, BoxMG
translates the inter-grid transfer operator computation into a set of  
small, decoupled linear equation system solves;
one solve per patch.
These solves decompose further: For each of the $2^d$ vertices, we have to
determine those $P$ entries that coincide with the patch.
These computations are independent of each other.

With a reference configuration at hand, the $C$ matrix for $d=2$ and the
reference configuration reads as

{\tiny
\[
\left(
\begin{array}{ccccccccc}
1 & 0 & 0 & 0 & 0 & 0 & 0 & 0 & 0 \\
s_{1,0}+s_{1,3}+s_{1,6} & 
     s_{1,1}+s_{1,4}+s_{1,7} &
        s_{1,2}+s_{1,5}+s_{1,8} & 0 & 0 & 0 & 0 & 0 & 0 \\
0 & s_{2,0}+s_{2,3}+s_{2,6} 
                  & s_{2,1}+s_{2,4}+s_{2,7} & 0 & 0 & 0 & 0 & 0 & 0\\        
s_{4,0}+s_{4,1}+s_{4,2} & 0 & 0 
          & s_{4,3}+s_{4,4}+s_{4,5} & 0 & 0 
                      & s_{4,6}+s_{4,7}+s_{4,8} & 0 & 0 \\
0 & 0 & 0 & s_{8,0}+s_{8,1}+s_{8,2} & 0 & 0 
                  & s_{8,3}+s_{8,4}+s_{8,5} & 0 & 0 \\
s_{5,0} & s_{5,1} & s_{5,2} s_{5,3} & s_{5,4} & s_{5,5} & s_{5,6} & s_{5,7} &
s_{5,8} \\
0 & s_{6,0} & s_{6,1} & 0 & s_{6,3} & s_{6,4} & 0 & s_{6,6} & s_{6,7} \\
0 & 0 & 0 & s_{9,0} & s_{9,1} & s_{9,2} & s_{9,3} & s_{9,4} & s_{9,5} \\
0 & 0 & 0 & 0 & s_{10,0} & s_{10,1} & 0 & s_{10,3} & s_{10,4}
\end{array}
\right)
\]
}

\paragraph{\new{$d=3$}} 
The $C$ matrix for $d=3$ is too big to write it down
here.
\new{All routines are however contained within one of the toolboxes (function
collections) available for the underlying software Peano \cite{Software:Peano}
from the project's repository.}

  \section{Remarks on the shared memory parallelization}

Patch-based strategies
\cite{Feichtinger:11:Walberla,Ghysels:13:ModelMG,Weinzierl:14:BlockFusion},
where patches of regular grids are embedded into cells, have been applied successfully for spacetrees and facilitate loop parallelism.
Such approaches even can be generalized in a multiscale way, where whole regions
are tessellated by a cascade of regular grids \cite{Gmeiner:14:HHG,Gmeiner:15:HHG}.
Alternatively, we may fix the grid, cut the linearized tree into chunks and
distribute those among threads \cite{Schreiber:13:Cluster}.

To obtain reasonable peak performance, such optimizations might become
necessary.
We do not study them here as they impose grid regularity constraints.
Instead, we focus on a task-based \new{parallelisation} formalism.
Combinations of both techniques seem to be promising to obtain high performance
in practice.
Furthermore, we note that the mapping of multigrid element activities onto tasks
yields a high theoretical concurrency but also yields high task management
overhead.
To reduce this overhead and, hence, to increase the arithmetic intensity,
our tasks have to be merged into bigger task assemblies \cite{Schreiber:13:Cluster,Weinzierl:14:BlockFusion}.
The exact choice of the size of such mergers is a non-trivial task
\cite{Charrier:17:Autotuning}.

  \section{Remarks on block smoothers}

\new{
Block Jacobi can not be used in the
coarsening step in the multiplicative algorithm variant as 
we fuse the computation of the
correction's right-hand side with a smoothing step on the new coarse grid. 
This is possible as the right-hand side on the correction grid is not required
to evaluate the element-wise operators. 
However, a block smoother on the new coarse grid would require such information.
}

\new{
This would not hold if we plugged into an $ascend$ operation that integrates
into the backtracking steps within the spacetree.
Only offering $ascend$ however would create a twin problem:
block smoothing throughout the coarsening would be possible but block smoothing
throughout a prolongation step would become impossible.
Combinations of $descend$ and $ascend$ solve the problem but sacrifice algorithmic simplicity. 
}

\new{
BoxMG's stencil collapsing
seems to be a natural candidate to realize better block smoothers:
If we apply collapsing on $\gamma $-points, these points can be subject to a
Gau\ss -Seidel smoother along the $\iota $-points of a patch.
Yet, we were not able to identify any significant convergence speedup in our
numerical experiments for such a smoother variant. This might be different for harder problems.
}

}{}

\ifthenelse{\boolean{toms}}{
  \pagebreak
  \section*{APPENDIX}
  \setcounter{section}{0}

  \section{Remarks on the shared memory parallelization}

Patch-based strategies
\cite{Feichtinger:11:Walberla,Ghysels:13:ModelMG,Weinzierl:14:BlockFusion},
where patches of regular grids are embedded into cells, have been applied successfully for spacetrees and facilitate loop parallelism.
Such approaches even can be generalized in a multiscale way, where whole regions
are tessellated by a cascade of regular grids \cite{Gmeiner:14:HHG,Gmeiner:15:HHG}.
Alternatively, we may fix the grid, cut the linearized tree into chunks and
distribute those among threads \cite{Schreiber:13:Cluster}.

To obtain reasonable peak performance, such optimizations might become
necessary.
We do not study them here as they impose grid regularity constraints.
Instead, we focus on a task-based \new{parallelisation} formalism.
Combinations of both techniques seem to be promising to obtain high performance
in practice.
Furthermore, we note that the mapping of multigrid element activities onto tasks
yields a high theoretical concurrency but also yields high task management
overhead.
To reduce this overhead and, hence, to increase the arithmetic intensity,
our tasks have to be merged into bigger task assemblies \cite{Schreiber:13:Cluster,Weinzierl:14:BlockFusion}.
The exact choice of the size of such mergers is a non-trivial task
\cite{Charrier:17:Autotuning}.

  \section{Remarks on block smoothers}

\new{
Block Jacobi can not be used in the
coarsening step in the multiplicative algorithm variant as 
we fuse the computation of the
correction's right-hand side with a smoothing step on the new coarse grid. 
This is possible as the right-hand side on the correction grid is not required
to evaluate the element-wise operators. 
However, a block smoother on the new coarse grid would require such information.
}

\new{
This would not hold if we plugged into an $ascend$ operation that integrates
into the backtracking steps within the spacetree.
Only offering $ascend$ however would create a twin problem:
block smoothing throughout the coarsening would be possible but block smoothing
throughout a prolongation step would become impossible.
Combinations of $descend$ and $ascend$ solve the problem but sacrifice algorithmic simplicity. 
}

\new{
BoxMG's stencil collapsing
seems to be a natural candidate to realize better block smoothers:
If we apply collapsing on $\gamma $-points, these points can be subject to a
Gau\ss -Seidel smoother along the $\iota $-points of a patch.
Yet, we were not able to identify any significant convergence speedup in our
numerical experiments for such a smoother variant. This might be different for harder problems.
}

  \acks
  This work made use of the facilities of the Hamilton HPC Service of Durham
University.
Furthermore, we particularly have been benefitting from the support of the RSC
Group who granted us early access to their KNL machines.
Thanks are due to Hans-Joachim Bungartz for supervising and mentoring the underlying PhD theses.
It was Irad Yavneh's suggestion to use the idea of BoxMG in this context,
and first steps into this direction where made in
\cite{yavnehbendig12nonnsymbb}. Marion wants to thank him for hosting her
research stay at Technion which laid the foundations for the present work.
Operator and data compression is a workpackage in the ExaHyPE project and Tobias
thus appreciates support received from the European Union’s Horizon 2020
research and innovation programme under grant agreement No 671698 (ExaHyPE).
All underlying software is open source \cite{Software:Peano}.

}{}

\ifthenelse{\boolean{toms}}{
 \bibliographystyle{ACM-Reference-Format-Journals}
 \bibliography{paper}
}{}

\end{document}